\documentclass {amsart}
\usepackage {amsrefs}
\usepackage [mathscr]{eucal}
\newtheorem {theorem}{Theorem}[section]

\newtheorem {lemma}[theorem]{Lemma}

\numberwithin {equation}{section}

\DeclareMathOperator*{\osc}{osc}

\begin {document}

\title [Gradient estimates]{Gradient estimates for oblique derivative problems via the maximum principle}

\author {Gary\ M.~ Lieberman}
\address {Department of Mathematics \\
Iowa State University\\
Ames, Iowa 50011}
\date {\today}
\maketitle
\newcommand \clOmega{\overline\Omega}

\section* {Introduction}
In this work, we study the gradient estimate for solutions of the boundary value problem
\begin {subequations} \label {EODP}
\begin {gather}
a^{ij}(x,u,Du)D_{ij}u +a(x,u,Du) =0 \text { in } \Omega, \label {EPDE} \\
b(x,u,Du)= 0 \text { on } \partial\Omega, \label {EBC}
\end {gather}
\end {subequations}
where $\Omega$ is a domain in $\mathbb R^n$ for some positive integer $n$ and we follow the summation convention.  Our two primary hypotheses are that the matrix $[a^{ij}]$ is positive definite and that the vector derivative $b_p(x,z,p)= \partial b(x,z,p)/\partial p$ satisfies the condition $b_p\cdot \gamma >0$ on $\partial\Omega$ for $\gamma$ the exterior unit normal to $\partial\Omega$. Such problems have been studied for a long time and are the topic of the book \cite {ODPbook}.  It is well-known that the key step in proving existence of solutions to this problem is a bound on the gradient of the solution, and, here, we study the gradient estimate under more general hypotheses than in previous works.  Specifically, we obtain a gradient bound under conditions on the differential equation that are based on those in \cite {Serrin}.  Unfortunately, we are unable to include some important special cases for reasons that will be described in more detail later.  Our model equation is the so-called false mean curvature equation, in which $a^{ij}=\delta^{ij}+p_ip_j$, where 
\[
\delta^{ij}=\begin {cases} 1 &\text { if } i=j, \\
0 &\text { if } i\neq j. \end {cases}
\]
The oblique derivative problem for this equation was first studied in \cite {NUE} and our results improve the ones there.  It is our hope that the method described here can be extended to other oblique derivative problems, such as the capillary problem
\begin {gather*}
\left(\delta^{ij} - \frac {D_iuD_ju}{1+|Du|^2}\right)D_{ij}u +a(x,u,Du)= 0 \text { in } \Omega, \\
\frac {Du\cdot\gamma}{(1+|Du|^2)^{1/2}} +\psi(x,u)=0 \text { on } \partial\Omega
\end {gather*}
under suitable conditions on $a$ and $\psi$ but we have not seen how to do so, yet.  Of course, this problem can be completely analyzed by various other methods, and we refer the interested reader to the Notes of Chapter 10 in \cite {ODPbook} as well as \cite {maxcap} and \cite {XuMa} for details on these other methods.

We also analyze the corresponding parabolic problem
\begin {subequations} \label {EPODP}
\begin {gather}
-u_t+a^{ij}(X,u,Du)D_{ij}u +a(X,u,Du) =0 \text { in } \Omega, \label {EPPDE} \\
b(X,u,Du)= 0 \text { on } S\Omega, \label {EPBC} \\
u=u_0 \text { on } B\Omega
\end {gather}
\end {subequations}
for a fairly general space-time domain $\Omega$, where we have followed the notation in \cite {LBook}.  More specifically, we first write $\mathcal P\Omega$ for the parabolic boundary of $\Omega$, that is, $\mathcal P\Omega$ is the set of all points $X_0=(x_0,t_0)$ in the topological boundary of $\Omega$ such that, for any $R>0$, the cylinder
\[
Q(X_0,R) =\{X\in\mathbb R^{n+1}: |x-x_0|<R, t_0-R^2<t<t_0\}
\]
contains at least one point not in $\Omega$. Then $S\Omega$ denotes the set of all $X_0\in \mathcal P\Omega$ such that, for any $R>0$, $Q(X_0,R)$ contains at least one point in $\Omega$.  Finally $B\Omega= \mathcal P\Omega\setminus S\Omega$.  The currently available situation for this problem is quite limited compared to that for the elliptic problem.  Gradient estimates are only known in three cases.  The first case is when the equation is uniformly parabolic in the sense that the eigenvalues of the matrix $[a^{ij}(X,z,p)]$ are bounded from above and below by positive constants for all $(X,z,p) \in \Omega\times \mathbb R\times \mathbb R^n$. In this case, gradient estimates appear in \cite {UR3} and Section 13.3 of \cite {LBook}. The second case is that of the conormal problem, which means that there is a vector-valued function $A$ such that $a^{ij}(X,z,p)= \partial A^i(X,z,p)/\partial p_j$ and
\[
b(X,z,p)= A(X,z,p)\cdot \gamma +\psi(X,z)
\]
for some scalar function $\psi$. In this case, gradient estimates appear in \cite {IUMJ88} (see Section 13.2 of \cite {LBook} for further discussion of this work).  Finally, Huisken \cite {H} and, more recently, Mizuno and Takasao \cite {Mizuno}, proved gradient estimates for  the problem
\begin {gather*}
- u_t+\left( \delta^{ij}- \frac {D_iuD_ju}{1+|Du|^2}\right)D_{ij}u =f(X,u,Du) \text { in } \Omega, \\
\gamma \cdot Du=0 \text { on } S\Omega, \\
u=u_0 \text { on } B\Omega
\end {gather*}
when $\Omega = \omega\times (0,T)$ for some domain $\omega\subset \mathbb R^n$ and $T>0$ as part of a more detailed study of this problem.  We are not concerned here with the complete analysis, but we do point out that Huisken only studied the case when $f(X,z,p)\equiv0$ to analyze the long time behavior of the solution, while Mizuno and Takasao assume an integrability condition on $f$ and only prove their gradient bound for small time.

Our present results expand on the first and third cases in this list (although our results do not really include those of Mizuno and Takasao).  In fact, our argument, when specialized to uniformly parabolic equations, is essentially identical to the one in Section 13.3 of \cite {LBook}.  More importantly, we obtain gradient bounds for a number of parabolic problems not previously accessible.

Our scheme is straightforward but, sadly, rather heavy on computation.  We begin with some basic assumptions and notation in Section \ref {S1} and some simple properties of monotone functions in Section \ref {Smono}.  In Section \ref {SN}, we introduce an auxiliary function which is crucial in our gradient estimates.  To apply our maximum principle argument, we present some preliminary calculations in Section \ref {Sp} which are then used in Sections \ref {S1GB} and \ref {Sso} to derive our gradient estimates for elliptic problems. Examples to illustrate these results are given in Section \ref {Se}.  We then turn to the parabolic problem, giving gradient estimates in Sections \ref {SpG} and \ref {SLp}.  We close with examples of our gradient for parabolic problems in Section \ref {Spe}.

Before beginning, we note an important limitation to our current gradient estimates.  Our maximum principle argument is based rather heavily on the corresponding analysis of Serrin \cite {Serrin} for interior gradient, but he introduced a decomposition
\[
a^{ij}(x,z,p)= a_*^{ij}(x,z,p) +f^i(x,z,p)p_j+f^j(x,z,p)p_i
\]
for some functions $a_*^{ij}$ and $f^i$.  This decomposition gives a number of striking results involving the structure of the coefficient $a$.  In particular, interior gradient estimates for solutions of the equation
\[
(\delta^{ij}+D_iuD_ju)D_{ij}u +|Du|^4=0
\]
are derived in \cite {Serrin}.  Such a decomposition is not available in our case since $p_i$ and $p_j$ would have to be replaced with more complicated functions of $p$ determined in a convoluted way from the boundary condition.  Hence we can only obtain gradient estimates for solutions of the
\[
(\delta^{ij}+D_iuD_ju)D_{ij}u +|Du|^q=0
\]
with $q\in (1,4)$.

\section {Basic assumptions and notation} \label {S1}

First, we say that  $\partial\Omega\in C^3$ if $\partial\Omega$ can be written as the level set of a $C^3$ function $f$ such that $Df$ doesn't vanish on $\partial\Omega$. We use $d$ to denote the distance function to $\partial\Omega$, and we denote by $\Omega_r$, for any positive number $r$, the subset of $\Omega$ on which $d<r$.
We then recall from Lemma 5.18 of \cite {ODPbook} (see also Section 10.3 of \cite {LBook}) that if $\partial\Omega\in C^3$, then there is a positive constant $R_0$, determined only by $\Omega$, such that $d\in C^{3}(\partial\Omega\cup\Omega_{R_0})$.
We set $\gamma=Dd$ in $\Omega_{R_0}$ and note that $\gamma$ is a $C^2$ unit vector. Moreover, for any vector $\xi$, if we define the vector $\tilde \xi$ by
\begin {equation} \label {E1.1}
\tilde \xi_i=D_i\gamma^k\xi_k,
\end {equation}
then $|\tilde \xi|\le 2|\xi|/R_0$ in $\Omega_{R_0/2}$ by virtue of Lemma 14.17 of \cite {GT}.

To describe our operator more easily, we define $v=(1+|p|^2)^{1/2}$ and $\nu = p/v$. We define the underlying matrix $[g^{ij}]$ by 
\[
g^{ij} = \delta^{ij}-\nu_i\nu_j,
\]
where $\delta^{ij}$ is the Kronecker $\delta$, that is, $\delta^{ij}=1$ if $i=j$ and $0$ otherwise.
We also set
\[
c^{ij}=\delta^{ij}-\gamma^i\gamma^j,
\]
and, for any $n$-dimensional vector $p$, we write $p'$ for the vector with components $p'_i=c^{ij}p_j$.
We also set
\[
v'=(1+|p'|^2)^{1/2}.
\]

We also adopt the following convention concerning derivatives of functions. First, we use subscripts to denote derivatives with respect to $x$, $z$, and $p$. So
\[
f_x= \frac {\partial f}{\partial x},\ f_z =\frac {\partial f}{\partial z}, \ f_p=\frac {\partial f}{\partial p}.
\]
We also use subscript indices to denote derivatives with respect to the corresponding coordinate of $x$:
\[
f_i = \frac {\partial f}{\partial x^i}
\]
and we use superscript indices to denote derivatives with respect to the corresponding component of $p$:
\[
f^i = \frac {\partial f}{\partial p_i}.
\]

\section {Some properties of monotone functions} \label {Smono}

Before we study the special function used in the remainder of this work, we present some properties of monotone functions.  Specifically, we show that a decreasing function is essentially equivalent to a smooth decreasing function that satisfies some simple differential inequalities.  The property of decreasing functions is essentially contained in Lemma 1.1 of \cite {CPDE16}.

\begin {lemma} \label {L1.1}
Let $L\ge0$. Then, for  any bounded, positive, decreasing function $f$, defined on $(L,\infty)$  there is a function $F\in C^2(L,\infty)$ with
\begin {subequations} \label {L1.1E}
\begin {gather}
f(x) \le F(x), \label {L1.1F}\\
-\frac {F(x)}x\le F'(x) \le 0, \label {L1.1F1}\\
- \frac {F(x)}{x^2}\le F''(x) \le  \frac {2F(x)}{x^2} \label {L1.1F2}, \\
\intertext {for all $x\in (L,\infty)$ and}
\lim_{x\to\infty} F(x)= \lim_{x\to\infty} f(x). \label {L1.1Flim}
\end {gather}
\end {subequations}
\end {lemma}
\begin {proof} By extending $f$ to be constant on the interval $[0,L]$, we may assume that $L=0$.
As in the proof of Lemma 1.1 from \cite {CPDE16}, we set
\[
g(x)= \frac 1x\int_0^x  f(y)\, dy,
\]
and we note that $g$ is continuous.  Since $f$ is decreasing, we have $g(x) \ge f(x)$.  Next, a simple integration by parts shows that
\begin {equation}
\int_\varepsilon^x\frac {g(y)}y\, dy = g(\varepsilon)-g(x) + \int_\varepsilon^x \frac {f(y)}y\, dy
\end {equation} \label {EDini}
for any positive $x$ and $\varepsilon$. It follows that
\[
g(x)= g(\varepsilon) +\int_\varepsilon^x \frac 1y(f(y)-g(y))\, dy,
\]
and hence $g$ is decreasing.

Next, we set 
\[
h(x)=\frac 1x\int_0^x g(y)\, dy
\]
and 
\[
F(x)= \frac 1x\int_0^x h(y)\, dy.
\]
Then $h$ is $C^1$ and $F$ is $C^2$, and the preceding arguments show that $h$ and $F$ are decreasing with $F(x)\ge h(x)\ge g(x)\ge f(x)$. In addition,
\[
xF'(x)+F(x)= h(x)\ge 0,
\]
which implies the first inequality of \eqref {L1.1F1}.

We now compute
\[
F'(x)= \frac {h(x)-F(x)}x,
\]
so
\[
F''(x) = \frac {F(x)-h(x)}{x^2}+ \frac {h'(x)-F'(x)}x = \frac {2(F(x)-h(x))}{x^2} + \frac {h'(x)}x.
\]
Since $F(x)\ge h(x)$ and $h'(x)\ge - h(x)/x\ge -F(x)/x$, we infer the first inequality of \eqref {L1.1F2}. The second inequality of \eqref {L1.1F2} follows because $h\ge 0\ge h'$.

Finally, \eqref {L1.1Flim} follows from repeated application of  l$'$H{\^o}spital's Rule.
\end {proof}

For ease of notation, we say that a function $F$ is {\it $*$-decreasing} if it satisfies \eqref {L1.1F1} and \eqref {L1.1F2}.  This concept will be useful in later sections.

We also make a corresponding observation for Dini functions, recalling that a function $f$ is \emph {Dini} if $f$ is increasing and 
\[
\int_0^L \frac {f(y)}y\,dy<\infty
\]
for some positive $L$. (In particular, $f$ must be nonnegative with $\lim_{x\to0} f(x)=0$.)  This observation will not be used here, but it may be helpful in other investigations.

\begin {lemma} \label {LDini}
If $f$ is Dini, then there is a $C^1$ Dini function $F$ such that 
\begin {subequations} \label {ELDini}
\begin {gather}
F(x) \ge f(x), \label {ELDini1}\\
F(x) \le 4f(4x), \label {ELDini4} \\
F'(x) \le \frac {F(x)}x, \label {ELDinig}
\end {gather}
\end {subequations}
for all sufficiently small positive $x$.
\end {lemma}
\begin {proof}
We define $g$ and $h$ as in the proof of Lemma \ref {L1.1} and set  $F(x) = 2h(2x)$. 
Then
\[
g(x) \ge \frac 1x\int_{x/2}^x f(x/2)\, dy = \frac {f(x/2)}2.
\]
Similarly, 
\[
h(x) \ge \frac {g(x/2)}2\ge \frac {f(x/4)}4,
\]
which yields \eqref {ELDini1}.
This time, the monotonicity of $f$ implies that $g$ and $h$ are increasing with $h(x) \le g(x) \le f(x)$, and this inequality yields \eqref {ELDini4}. In addition, \eqref {EDini} and the monotone convergence theorem imply that $g$ is Dini.  It follows that $h$ is Dini and hence so is $F$.  The proof of \eqref {ELDinig} follows that of the first inequality of \eqref {L1.1F1}.
\end {proof}

\section {The $N$ function} \label {SN}

In this section, we introduce a function $N$, determined by the boundary condition, which is key to deriving our estimates.
Such a function was first presented in \cite {Dong} and our function takes advantage of a modification due to the present author \cite {IUMJ88} that allows us to study boundary conditions under weaker hypotheses than in \cite {Dong}.
In the previous applications of this function, a positive parameter $\varepsilon$ is introduced and then fixed at a particular value.  For our purposes, it will be very important to keep track of this parameter throughout our work, so we rewrite the results and present the proof with a careful accounting of this parameter.

To state the properties of this function in a useful way, we introduce one additional bit of terminology that will be useful.
For $\tau\ge1$, we write $\Sigma(\tau)$ for the subset of $\partial\Omega\times\mathbb R\times\mathbb R^n$ on which $v\ge\tau$.

\begin {theorem} \label {TN}
Let $\partial\Omega\in C^2$.  Let $\tau_0\ge1$ be given and let $b\in C^1(\Sigma(\tau_0))$ with $b_p\cdot\gamma >0$ on $\Sigma_0(\tau_0)$, the subset of $\Sigma(\tau_0)$ on which $b=0$.
Suppose
\begin {equation} \label {E10.32}
\lim_{t\to\infty} b(x,z,p-\tau\gamma(x))<0<\lim_{t\to\infty} b(x,z,p+\tau\gamma(x))
\end {equation}
for all $(x,z,p)\in \Sigma(\tau_0)$, and that there are positive constants $\beta_0$ and $c_0$ such that
\begin {equation} \label {E10.33}
|b_p| \le \beta_0b_p\cdot \gamma, \quad
|p\cdot\gamma| \le c_0v'
\end {equation}
on $\Sigma_0(\tau_0)$.
Suppose also that there are a $*$-decreasing function $\varepsilon_x$ and a nonnegative constant $\beta_1$ such that
\begin {subequations} \label {E10.35}
\begin {gather}
|b_x| \le \varepsilon_x(v)v^2b_p\cdot\gamma \\
|b_z| \le \beta_1vb_p\cdot\gamma
\end {gather}
\end {subequations}
on $\Sigma_0(\tau_0)$.
Then there is a positive constant $\varepsilon_0(\beta_0,c_0)$, along with  a $C^2(\overline{\Omega_{R_0/4}}\times\mathbb R\times \mathbb R^n\times (0,\varepsilon_0))$ function $N$ such that
\begin {subequations} \label {E10.37}
\begin {gather}
\frac 12\le N_p\cdot\gamma \le \frac 32, \label {E10.37a} \\
|N_p| \le 3\beta_0, \label {E10.37b} \\
|N_z| \le 4(1+c_0)\beta_1v_\varepsilon, \label {E10.37c} \\
|N_x| \le 55(1+c_0)\frac {\beta_0}{R_0}v +12(1+c_0)^2\varepsilon_x(v)v^2 , \label {E10.37d} \\
|N-N_p\cdot p| \le 6\beta_0(1+c_0)^2v_\varepsilon, \label {E10.372dot}
\end {gather}
\end {subequations}
on $\Omega_{R_0/4}\times\mathbb R\times\mathbb R^n\times(0,\varepsilon_0)$, where the arguments of $N$ are $(x,z,p,\varepsilon)$ and $v_\varepsilon = \sqrt{1+|p'|^2+\varepsilon|p\cdot\gamma|^2}$. Moreover, $N=0$ on $\Sigma_0( (1+c_0^2)^{1/2}\tau_0) \times(0,\varepsilon_0)$. If we set
\begin {equation} \label {E10.38}
W=|p'|^2+\varepsilon N^2, \ \nu_1= p'+\varepsilon NN_p
\end {equation}
for some $\varepsilon\in(0,\varepsilon_0)$, then
\begin {subequations} \label {E10.390}
\begin {gather}
 (1-8(1+c_0)^2\varepsilon^{1/2})v_\varepsilon^2 \le 1+W \le (1+8(1+c_0)^2\varepsilon^{1/2})v_\varepsilon^2,  \label {E10.39}\\
|\nu_1| \le 2v_\varepsilon, \label {E10.41a} \\
|W-p\cdot\nu_1| \le 12(1+c_0)^3\beta_0\varepsilon^{1/2}v_\varepsilon^2. \label {EWnu1}
\end {gather}
\end {subequations}
Also, there is a nonnegative constant $c_1$, determined only by $\beta_0$, $c_0$ and $n$ such that 
\begin {subequations} \label {E10.372}
\begin {gather}
|NN^{km}\xi_k\xi_m| \le \frac{c_1}{\varepsilon^{1/2}}|\xi'|^2+ \frac 14|\xi\cdot\gamma|^2 , \label {E10.372pp} \\
(1- c_1\varepsilon^{1/2})|\xi'|^2+ \frac 12\varepsilon|\xi\cdot\gamma|^2 \le \left[c^{km}+\varepsilon NN^{km}+\varepsilon N^kN^m\right]\xi_k\xi_m,  \label {E10.43}\\
 |NN_{pz}\cdot\xi| \le c_1\left(\frac 1{\varepsilon^{1/2}}|\xi'|+|\xi\cdot\gamma|\right) \beta_1v_\varepsilon, \label {E10.372pz} \\
|NN^k_x\xi_k| \le c_1\left(\frac 1{\varepsilon^{1/2}}|\xi'|+|\xi\cdot\gamma|\right) \left[\frac v{R_0}+\varepsilon_x(v)v^2\right] \label {E10.372px} \\
\intertext {for any $\xi\in \mathbb R^n$,}
|NN_{zz}|  \le c_1\beta_1^2v_\varepsilon^2, \label {E10.372zz} \\
|NN_{xz}| \le \frac {c_1}{\varepsilon^{1/2}}\beta_1v_\varepsilon\left[\frac v{R_0}+\varepsilon_x(v)v^2\right] . \label {E10.372xz}\\
\intertext {Finally, if $\partial\Omega\in C^3$, then there is a nonnegative constant $c_2$, determined only by $\Omega$, such that}
|NN_{xx}|\le \frac {c_2}{\varepsilon^{1/2}}\left[\frac v{R_0}+\varepsilon_x(v)v^2\right] ^2+c_2v. \label {E10.372xx}
\end {gather}
\end {subequations}
\end {theorem}
\begin {proof} Although the results are mostly contained in Lemma 10.8 of \cite {ODPbook} and Lemma 4.3 of \cite {NA119}, we sketch the proof here because there are some points that are not immediate consequences of the arguments there.  Specifically, Lemma 10.8 of \cite {ODPbook} only studies second derivatives of $N$ with respect to $p$ (and the argument there will not give the existence of the other second derivatives) while Lemma 4.3 of \cite {NA119} assumes that $\Omega$ is the upper half-plane. In addition, the exact form of the estimates in terms of  $\varepsilon$ was not studied in either of those works.

In our proof, we shall assume that $R_0$ is finite and $\varepsilon_x$ and $\beta_1$ are positive. The other cases are proved by similar but (sometimes) simpler arguments.
We start by noting (exactly as in Lemma 10.7 of \cite {ODPbook}) that there is a function $g$, defined on $\Omega^*_{R_0}\times\mathbb R\times \mathbb R^n$ (where $\Omega^*_{R_0}$ is the set of all $x\in\mathbb R^n$ with $\inf\{|x-y|:y\in\partial\Omega\}<R_0$ including points outside $\Omega$), such that $p\cdot\gamma=g(x,z,p)$ if and only if $b(x,z,p)=0$. Specifically, we have
\begin {equation} \label {Egdef}
0= b(x,z,p+\gamma [g(x,z,p)-p\cdot\gamma]).
\end {equation}
It follows that  $g(x,z,p)= g(x,z,q)$ whenever $p\cdot\gamma = q\cdot\gamma$ and that
\begin {subequations}
\begin {gather}
|g(x,z,p)|\le c_0v',  \label {Egvprime} \\
|g_x| \le \theta_x(v)v^2, \quad |g_z|\le \beta_1v', \quad |g_p| \le \beta_0, \label {Egderivs}
\end {gather}
\end {subequations}
where the function $\theta_x$ is defined by
\[
\theta_x(\sigma)= \frac {4(1+c_0)\beta_0}{R_0\sigma} +(1+c_0)^2\varepsilon_x(\sigma).
\]
Note that $\theta_x$ is $*$-decreasing.

We now let $\varphi$ be a nonnegative $C^\infty(\mathbb R^{2n+1})$ function with support in the unit ball and 
\[
\int \varphi(Y)\,dY=1,
\]
where here and in the remainder of this proof, all integrals are taken over $\mathbb R^{2n+1}$.
We also write $Y=(y,w,q)$ with $y\in \mathbb R^n$, $w\in\mathbb R$ and $q\in \mathbb R^n$ and we set $\tilde v = \sqrt{1+|p|^2+s^2}$ and $\tilde v' = \sqrt{1+|p'|^2+s^2}$. 
With $K= 1/(6\beta_0)$ and $\varepsilon'$ a positive constant to be determined, we introduce the following abbreviations:
\begin {align*}
x^* &= x+ \frac {\varepsilon'sy}{\theta_x(\tilde v)\tilde v^2}, \\
z^* &= z+ \frac {\varepsilon'sw}{\beta_1\tilde v'}, \\
p^*&= p+ Ksq, \\
X^*&= (x^*, z^*,p^*),
\end {align*}
and we define the function $\tilde g$ by
\[
\tilde g(x,z,p,s) = \int g(X^*)\varphi(Y)\, dY.
\]
An elementary calculation shows that 
\[
\frac {1}{\theta_x(\tilde v)\tilde v}\le  \frac {R_0}{4},
\]
 so if  $\varepsilon \le1$, then $\tilde g$ is defined for all $(x,z,p,s)\in \Omega_{R_0/4}\times \mathbb R\times \mathbb R^n\times \mathbb R$. 

To proceed, we note that, even though $g$ depends on $p'$ rather than $p$, $\tilde g$ may also depend on $p\cdot\gamma$ because $\gamma$ changes with $x$. (This situation does not arise in \cite {NA119} and it is not relevant for the argument in \cite {ODPbook}.)
In particular, we have
\[
\gamma(x)- \gamma ( x^*)= D_i\gamma(x^{**}) \frac {\varepsilon' sy^i}{\theta_x(\tilde v)\tilde v^2}
\]
for some point $x^{**}$ on the line segment between $x$ and  $x^*$ and hence 
\begin {equation} \label {Egammadiff}
\left|\gamma(x)- \gamma (x^*)\right| \le \frac {\varepsilon'|s|}{\tilde v}.
\end {equation}
Using this inequality, we can estimate the first derivatives of $\tilde g$. We start by introducing two more functions $h_1$ and $h_2$, defined by
\[
h_1(\sigma) = \frac {\varepsilon'}{\theta_x(\sigma)\sigma^2}, \quad h_2(\sigma) = \frac {\varepsilon'}{\beta_1\sigma}.
\]
A simple computation shows that
\[
\tilde g_s(x,z,p,s) = I_1+I_2+I_3,
\]
with
\begin {align*}
I_1 &= \int g_i(X^*)y^i \left[ h_1(\tilde v) -\frac {s^2h_1'(\tilde v)}{\tilde v}\right] \varphi(Y)\, dY, \\
I_2 &= \int g_z(X^*)w\ \left[ h_2(\tilde v') -\frac {s^2h_2'(\tilde v')}{\tilde v'}\right]\varphi(Y)\, dY\\
I_3 &= \int g^k(X^*) Kq_k \varphi(Y)\, dY.
\end {align*}
To estimate  $I_1$ and $I_2$ (and to estimate many of our later integrals), we begin by using \eqref {Egammadiff} and noting that $(p+Ksq)'_m = c^{km}(x^*)(p+Ksq)_k$ to infer that
\[
\left|\left(p+Ksq\right)'\right| \le |p'| + |s| + (|p|+|s|) \frac {2\varepsilon'|s|}{\tilde v}
\]
and hence
\[
\left( 1+ \left|\left(p+Ksq\right)'\right|^2\right)^{1/2}
 \le\sqrt2(1+2\varepsilon')|\tilde v'|.
\]
If $\varepsilon'\le \frac 1{5}$ (so that $1+2\varepsilon^{1/2}\le \sqrt2$), then it follows that
\begin {subequations} 
\begin {gather}
|g_x(X^*)| \le 2\theta_x(\tilde v)\tilde v^2,  \label {Egi} \\
|g_z(X^*)| \le 2\beta_1\tilde v'. \label {Egz}
\end {gather}
\end {subequations}
Next, we invoke \eqref {L1.1F1} to see that
\[
\left| h_1(\tilde v) -\frac {s^2h_1'(\tilde v)}{\tilde v}\right| \le 1
\]
for any $s$, and hence $|I_1| \le 4\varepsilon'$, and a similar argument gives $|I_2| \le 2\varepsilon'$. The choice of $K$ implies that $|I_3| \le 1/6$, so 
\begin {equation} \label {Etildegs}
|\tilde g_s| \le \frac 12
\end {equation}
provided $\varepsilon' \le 1/18$.
We now compute and estimate the other derivatives of $\tilde g$.  First,
\[
\tilde g_i(x,z,p,s) = \int g_i(X^*)\varphi(Y)\, dY+ \int g_z(X^*)sh_2'(\tilde v') \frac {D_i(c^{km})p_kp_m}{\tilde v'}\varphi(Y)\, dY,
\]
so \eqref {Egi} and \eqref {Egz} imply
\[
|\tilde g_x| \le 3\theta_x(\tilde v)\tilde v^2.
\]
Similarly,
\[
\tilde g_z(x,z,p,s) = \int g_z(X^*)\varphi(Y)\, dY,
\]
so
\[
|\tilde g_z| \le 2\beta_1\tilde v'.
\]

Finally, we compute
\[
\tilde g_p(x,z,p,s) = J_1+J_2 +J_3
\]
with
\begin {align*}
J_1 &= -p\int g_i(X^*)y^i \frac {sh_1'(\tilde v)}{\tilde v} \varphi(Y)\, dY, \\
J_2 &=-p' \int g_z(X^*) w\frac {sh_2'(\tilde v')}{\tilde v'}\varphi(Y)\, dY, \\
J_3 &= \int g_p(X^*)\varphi(Y)\, dY.
\end {align*}
It's easy to check that $|J_1|+|J_2| \le 4\varepsilon'|s|/\tilde v$. 
The analysis of $J_3$ is more subtle.  For any vector $\xi$, we have
\[
\xi_k\int g^k(X^*)\varphi(Y)\, dY= \xi'_k\int g^k(X^*) \varphi(Y)\, dY + \xi\cdot \gamma \int g^k(X^*)[\gamma_k(x)-\gamma_k(x^*)]\, dY
\]
because $g_p\cdot\gamma =0$. It follows that
\begin {equation} \label {Etildegp}
|\tilde g^k\xi_k| \le \beta_0(1+4\varepsilon')|\xi'| + 5\beta_0 \varepsilon' \frac {|s|}{\tilde v}|\xi\cdot\gamma|.
\end {equation}
because $\beta_0\ge 1$.
Due to \eqref {Etildegs}, the equation
\[
p\cdot\gamma - \tilde g(x,z,p,\varepsilon^{1/2}N)=N
\]
defines $N$ as a function of $(x,z,p)$ and $\varepsilon$. 
By writing
\[
N= p\cdot\gamma- g(x,z,p)+[g(x,z,p)-\tilde g(x,z,p,\varepsilon^{1/2}N)],
\]
we infer that
\[
|N| \le |p\cdot\gamma-g(x,z,p)| +\frac 12\varepsilon^{1/2}|N|,
\]
and hence
\begin {equation} \label {EgN}
|N(x,z,p;\varepsilon)| \le \frac {1}{1-\varepsilon^{1/2}}(|p\cdot\gamma|+|g(x,z,p)|).
\end {equation}
It then follows from \eqref {Egvprime} that
\begin {equation} \label {ENv}
|N| \le 2(1+c_0)v
\end {equation}
if $\varepsilon_0\le 1/4$.  Since 
\[
\tilde v'= \sqrt{(v')^2+\varepsilon N^2},
\]
we have $\tilde v' \le 2(1+c_0)v_\varepsilon$.
We then infer that
\[
N_p\cdot\gamma -1 = \frac {-\tilde g_p\cdot\gamma-\tilde g_s\varepsilon^{1/2}}{1+\tilde g_s\varepsilon^{1/2}}
\]
so our estimates for the derivatives of $\tilde g$ imply \eqref {E10.37a}, \eqref {E10.37b}, \eqref {E10.37c}, and \eqref {E10.37d} provided $\varepsilon_0\le (10(1+c_0))^{-2}$ and $\varepsilon'\le 1/(20\beta_0)$.

Differentiating \eqref {Egdef} with respect to $p$ shows that
\[
g_p- \gamma =- \frac {b_p}{b_p\cdot\gamma},
\]
where $b_p$ is evaluated at $(x,z,p'+\gamma g(x,z,p))$ since $p'+\gamma g= p+\gamma(g-p\cdot\gamma)$, and hence, after taking the dot product of this equation with $p'+\gamma g$, we find from \eqref {E10.33} that 
\[
|g_p\cdot p-g| \le \beta_0(1+|p'|^2+g(z,x,p)^2)^{1/2}.
\]
Another application of \eqref {Egvprime} shows that
\begin {equation} \label {Egp}
|g_p\cdot p-g| \le (1+c_0)\beta_0v'.
\end {equation}

Next, we recall that
\[
\tilde g_p(x,z,p,s) = J_1+J_2+J_3, \quad \tilde g_s(x,z,p,s) =I_1+I_2+I_3
\]
with $|J_1|+|J_2| \le 4\varepsilon'|s|/\tilde v$ and $|I_1|+|I_2| \le 6\varepsilon'$. In addition,
\[
J_3=\int g_p(X^*)\varphi(Y)\, dY, \quad I_3= \int g_p(X^*)\cdot(Kq)\, dY,
\]
so
\begin {align*}
|\tilde g_p\cdot p+\tilde  g_ss-\tilde g| &\le 10\varepsilon'|s| +\int |g_p(X^*)\cdot p^*-g(X^*)|\varphi(Y)\, dY \\
&\le (1+c_0)\beta_0\tilde v'+ 10\varepsilon'|s|,
\end {align*}
with $\tilde g$ and its derivatives evaluated at $(x,z,p,s)$, by \eqref {Egp}. A direct computation gives
\[
N-N_p\cdot p= \frac {-\tilde g +\tilde g_p\cdot p +\tilde g_s \varepsilon^{1/2}N}{1- \tilde g_s\varepsilon^{1/2}}
\]
with $\tilde g$ and its derivatives now evaluated at $(x,z,p,\varepsilon^{1/2}N)$, so
\[
|N-N_p\cdot p| \le 2[(1+c_0)\beta_0(v+\varepsilon^{1/2}|N|) + 10\varepsilon'\varepsilon^{1/2} |N|] .
\]
An application of \eqref {ENv} yields \eqref {E10.372dot}.

A simple modification of the argument leading to \eqref {EgN} shows that
\[
|N| \ge \frac 1{1+\varepsilon^{1/2}}|p\cdot\gamma-g|.
\]
The Cauchy-Schwarz inequality then gives
\begin {align*}
N^2 &\ge \frac {1-\varepsilon^{1/2}}{(1+\varepsilon^{1/2})^2}(p\cdot\gamma)^2 - \frac {1+\varepsilon^{1/2}}{(1+\varepsilon^{1/2})^2}g(x,z,p)^2 \\
&\ge (1-3\varepsilon^{1/2})(p\cdot\gamma)^2 -\frac {2}{\varepsilon^{1/2}}g(x,z,p)^2 
\end {align*}
if $\varepsilon\le 1$.
A similar argument, using \eqref {EgN} shows that
\[
N^2 \le (1+4\varepsilon^{1/2})(p\cdot\gamma)^2+ \frac 8{\varepsilon^{1/2}}g(x,z,p)^2
\]
if $\varepsilon \le 1/\sqrt 2$.
Using these inequalities along with \eqref {Egvprime} gives \eqref {E10.39}.

Estimates \eqref {E10.37b},  \eqref {ENv}, and the Cauchy-Schwarz inequailty imply \eqref {E10.41a} provided $\varepsilon_0\le 1/(36\beta_0^2(1+c_0)^2)$.

Since $|W-p\cdot \nu_1| = \varepsilon|N||N-N_p\cdot p|$, \eqref {EWnu1} follows from \eqref {ENv}, \eqref {E10.372dot}, and the observation that $\varepsilon^{1/2}v\le v_\varepsilon$.

To estimate the second derivatives of $N$, we use integration by parts repeatedly.  First, for the second derivatives with respect to $p$ to obtain
\begin {align*}
J_1 &= p\int g(X^*) \frac {h_1'(\tilde v)}{h_1(\tilde v)\tilde v} \frac {\partial}{\partial y^i}(y^i\varphi(Y))\, dY , \\
J_2 &= -p'\int g(X^*) \frac 1{(\tilde v')^2}\frac {\partial}{\partial w}(w\varphi(Y))\, dY, \\
J_3 &= - \frac 1{Ks}\int g(X^*)\varphi_q(Y)\, dY.
\end {align*}
Straightforward computation then shows that, for any vector $\xi$, we have
\[
\tilde g^{km}\xi_k\xi_m
= J_{11}+J_{12} +J_{13}+ J_{14}+J_{15}+ J_{21}+J_{22}+J_{23}+J_{24} +J_{25}
+J_{31}+J_{32} +J_{33},
\]
where $J_{11}, \dots, J_{15}$ come from differentiating $J_1$:
\begin {align*}
J_{11} &=-|\xi|^2 \int g(X^*) \frac {h_1'(\tilde v)}{h_1(\tilde v)\tilde v}\frac {\partial}{\partial y^i}(y^i\varphi(Y))\, dY , \\
J_{12} &= -(p\cdot \xi)^2 \int g_j(X^*)y^jh_1'(\tilde v) \frac {sh_1'(\tilde v)}{h_1(\tilde v)\tilde v^2}\frac {\partial}{\partial y^i}(y^i\varphi(Y))\, dY , \\
J_{13} &= -(p\cdot \xi)(p'\cdot\xi')\int g_z(X^*)wh_2'(\tilde v') \frac {sh_1'(\tilde v)}{h_1(\tilde v)\tilde v\tilde v'} \frac {\partial}{\partial y^i}(y^i\varphi(Y))\, dY, \\
J_{14} &= -(p\cdot\xi)\int g_p(X^*)\cdot\xi \frac {h_1'(\tilde v)}{h_1(\tilde v)\tilde v}\frac {\partial}{\partial y^i}(y^i\varphi(Y))\, dY , \\
J_{15}&=-(p\cdot\xi)^2\int g(X^*)\left[\frac{h_1''(\tilde v)}{h_1(\tilde v)\tilde v^2}- \frac {h_1'(\tilde v)(h_1'(\tilde v)\tilde v-h_1(\tilde v))}{h_1(\tilde v)^2\tilde v^3}\right]\frac {\partial}{\partial y^i}(y^i\varphi(Y))\, dY ;
\end {align*}
 $J_{21}, \dots, J_{25}$ come from differentiating $J_2$:
\begin {align*}
J_{21} &=-|\xi'|^2 \int g(X^*) \frac {1}{(\tilde v')^2}\frac {\partial}{\partial w}(w\varphi(Y))\, dY , \\
J_{22} &= -(p\cdot \xi)(p'\cdot \xi') \int g_i(X^*)y^ih_1'(\tilde v)\frac s{(\tilde v')^2\tilde v}\frac {\partial}{\partial w}(w\varphi(Y))\, dY , \\
J_{23} &= - (p'\cdot\xi')^2 \int g_z(X^*) wh_2'(\tilde v')\frac s {(\tilde v')^3}\frac {\partial}{\partial w}(w\varphi(Y))\, dY, \\
J_{24} &= (p'\cdot\xi')\int g_p(X^*)\cdot\xi \frac 1{(\tilde v')^2}\frac {\partial}{\partial w}(w\varphi(Y))\, dY , \\
J_{25}&=2(p'\cdot\xi')^2\int g(X^*)\frac 1{(\tilde v')^4}\frac {\partial}{\partial w}(w\varphi(Y))\, dY;
\end {align*}
and $J_{31}$, $J_{32}$, and $J_{33}$ come from differentiating $J_3$:
\begin {align*}
J_{31} &= -\frac {p\cdot \xi}{Ks\tilde v} \int g_i(X^*)y^ih_1'(\tilde v) s \frac {\partial \varphi(Y)}{\partial q_k}\xi_k\, dY, \\
J_{32} &= - \frac 1{Ks} \int g_z(X^*) wh_2'(\tilde v') s \frac {p'\cdot \xi'}{\tilde v'} \frac {\partial\varphi(Y)}{\partial q_k}\xi_k\, dY, \\
J_{33} &= -\frac 1{Ks} \int g_p(X^*)\cdot\xi \frac {\partial \varphi(Y)}{\partial q_k}\xi_k\, dY.
\end {align*}

To estimate $J_{11}$, we first integrate by parts and then rewrite the derivative with respect to $y$.  In this way, we see that
\begin {align*}
J_{11} &=-|\xi|^2 \int \frac {\partial g(X^*)}{\partial y^i}y^i \frac {h_1'(\tilde v)}{h_1(\tilde v)\tilde v}\varphi(Y)\, dY \\
&=-|\xi|^2\int  g_i(X^*)y^i \frac {sh_1'(\tilde v)}{\tilde v}\varphi(Y)\, dY.
\end {align*}
Therefore $|J_{11}| \le 2\varepsilon'|\xi|^2/\tilde v$.  In a similar fashion, we find that
\[
|J_{11}|+|J_{15}| \le C(n,c_0,\beta_0)\varepsilon'\frac {|\xi|^2}{\tilde v}.
\]

It is straightforward to estimate $J_{12}$,  $J_{13}$, and $J_{31}$.  The resultant inequality is
\[
|J_{12}|+|J_{13}|+|J_{31}| \le C(n,c_0,\beta_0) \varepsilon'\frac {|\xi|^2}{\tilde v}.
\]

To estimate $J_{14}$, we use the decomposition
\[
J_{14} = J_{14a}+J_{14b}
\]
with
\begin {align*}
J_{14a} &= (p\cdot\xi) \xi'_k \int g^k(X^*) \frac {h_1'(\tilde v)}{\tilde vh_1(\tilde v)} \frac {\partial}{\partial y^i} (y^i\varphi(Y))\, dY \\
J_{14b} &= (p\cdot\xi) (\xi\cdot \gamma) \int g^k(X^*)[\gamma_k(x)-\gamma_k(x^*)] \frac {h_1'(\tilde v)}{\tilde vh_1(\tilde v)} \frac {\partial}{\partial y^i} (y^i\varphi(Y))\, dY,
\end {align*}
recalling that $g_p(X^*)\cdot\gamma(x^*)=0$.
We then have
\[
|J_{14a}| \le C(n,\beta_0) \frac {|\xi||\xi'|}{\tilde v}
\]
and, by virtue of \eqref {Egammadiff},
\[
|J_{14b}| \le C(n,\beta_0)\varepsilon' \frac {|\xi||\xi\cdot\gamma|}{\tilde v}.
\]
It follows that 
\[
|J_{14}| \le C(n,\beta_0) \frac {|\xi|}{\tilde v}( |\xi'|+\varepsilon'|\xi\cdot\gamma|).
\]

The estimates for $J_{21}$ and $J_{25}$ follow the same idea as for $J_{11}$, yielding
\[
|J_{21}| +|J_{25}| \le C(n,\beta_0) \varepsilon' \frac{|\xi||\xi'|}{\tilde v'}.
\]

Straightforward calculation gives the following estimates for $J_{22}$ and $J_{23}$:
\begin {align*}
|J_{22}|  &\le  C(n,\beta_0)\varepsilon'\frac {|\xi||\xi'|}{\tilde v}, \\
|J_{23}| &\le C(n,\beta_0) \frac {|\xi'|^2}{\tilde v'}.
\end {align*}

By using the same decomposition as for $J_{14}$, we find that
\[
|J_{24}| \le C(n,\beta_0) \left(\frac {|\xi'|^2}{\tilde v'} + \frac { \varepsilon'|\xi'||\xi\cdot\gamma|}{\tilde v}\right).
\]
A simple integration by parts shows that $J_{31}=J_{14}$  and hence
\[
|J_{32}| \le C(n,\beta_0) \frac {|\xi|}{\tilde v}( |\xi'|+\varepsilon'|\xi\cdot\gamma|).
\]
Integration by parts also shows that $J_{32} = J_{23}$, so 
\[
J_{32} \le C(n,\beta_0) \frac {|\xi'|^2}{\tilde v'}.
\]

The estimate for $J_{33}$ is somewhat more complex.  First, we write
\[
J_{33}= J_{33a}+J_{33b}
\]
with
\begin {align*}
J_{33a} &= -\frac {1}{Ks}\xi'_m\int g^m(X^*)\frac {\partial \varphi(Y)}{\partial q_k} \xi_k\, dY, \\
J_{33b} &=-\frac 1{Ks}(\xi\cdot \gamma)\int g^m(X^*)[\gamma_m(x)-\gamma_m(x^*)] \frac {\partial \varphi(Y)}{\partial q_k} \xi_k\, dY.
\end {align*}
Then we write
\[
J_{33a} = J_{33c}+ J_{33d}
\]
with
\begin {align*}
J_{33c} &= -\frac {1}{Ks}\xi'_m\xi'_k\int g^m(X^*)\frac {\partial \varphi(Y)}{\partial q_k}\, dY, \\
J_{33d} &=-\frac 1{Ks}\xi'_m(\xi\cdot \gamma)\int g^m(X^*) \frac {\partial \varphi(Y)}{\partial q_k}\gamma_k(x) \, dY.
\end {align*}
It follows that $|J_{33c}|\le C(n,\beta_0)|\xi'|^2/|s|$ and an integration by parts yields
\[
J_{33d} =-\frac 1{Ks\tilde v}\xi'_m(\xi\cdot \gamma)\int g^k(X^*)[\gamma_k(x)-\gamma_k(x^*)] \frac {\partial \varphi(Y)}{\partial q_m}\, dY.
\]
It follows that
\[
|J_{33d}| \le C(n,\beta_0)\varepsilon'\frac {|\xi'||\xi\cdot\gamma|}{\tilde v}.
\]
Arguing as before, we also have
\[
|J_{33b}|\le C(n,\beta_0)\varepsilon\frac {|\xi||\xi\cdot\gamma|}{\tilde v}
\]
and therefore
\[
|\tilde g^{km}\xi_k\xi_m| \le C(n,c_0,\beta_0)|\left(\varepsilon'\frac {|\xi|^2}{\tilde v} + \frac {|\xi||\xi'|}{\tilde v} +\varepsilon'\frac {|\xi||\xi\cdot\gamma|}{\tilde v} + \frac{|\xi'|^2}{|s|}\right).
\]
By using the Cauchy-Schwarz inequality along with the inequalities $\varepsilon'\le 1$, and $|s|\le \tilde v'\le \tilde v$, and $v\le \tilde v$, we conclude that
\begin {equation} \label {Etildegpp}
|\tilde g^{km}\xi_k\xi_m| \le  C(n,\beta_0)\left(  \frac {|\xi'|^2}{\varepsilon'|s|} + \frac {\varepsilon'(\xi\cdot\gamma)^2}{ v}\right).
\end {equation} 

To estimate the derivative $g_{sp}$, we integrate $I_1$, $I_2$, and $I_3$ by parts and then differentiate the resultant integrals to obtain
\[
\tilde g^k_s\xi_k = I_{11}+ I_{12}+ I_{13}+I_{14}+ I_{21}+I_{22}+I_{23}+I_{24}+I_{31} +I_{32}+I_{33}
\]
for any vector $\xi$ with
\begin {align*}
 I_{11} &= -\frac {p\cdot\xi}{\tilde v}\int g_j(X^*)y^j h_1'(\tilde v)\left[1- \frac {s^2h_1'(\tilde v)}{h_1(\tilde v)\tilde v} \right] \frac {\partial}{\partial y^i}(y^i\varphi(Y))\, dY, \\
I_{12} &= -\frac {p'\cdot\xi'}{\tilde v'}\int g_z(X^*)wh_2'(\tilde v')\left[1- \frac {s^2h_1'(\tilde v)}{h_1(\tilde v)\tilde v} \right] \frac {\partial}{\partial y^i}(y^i\varphi(Y))\, dY, \\
I_{13} &=-K\int g^k(X^*)\xi_k\left[1- \frac {s^2h_1'(\tilde v)}{h_1(\tilde v)\tilde v}\right] \frac {\partial}{\partial y^i}(y^i\varphi(Y))\, dY, \\
I_{14} &= -\frac{p\cdot\xi}{\tilde v}\int g(X^*) s \frac {h_1''(\tilde v) h_1(\tilde v)\tilde v-(h_1'(\tilde v))^2\tilde v- h_1'(\tilde v)h_1(\tilde v)}{h_1(\tilde v)^2\tilde v^2}
\frac {\partial}{\partial y^i}(y^i\varphi(Y))\, dY;
\end {align*}
and
\begin {align*}
 I_{21} &= -\frac {p'\cdot\xi'}{\tilde v'}\int g_j(X^*)y^j h_1'(\tilde v)\left[1+ \frac {s^2}{(\tilde v')^2} \right] \frac {\partial}{\partial w}(w\varphi(Y))\, dY, \\
I_{22} &= -\frac {p'\cdot\xi'}{\tilde v'}\int g_z(X^*)wh_2'(\tilde v)\left[1+ \frac {s^2}{(\tilde v')^2} \right] \frac {\partial}{\partial w}(w\varphi(Y))\, dY, \\
I_{23} &= -\frac 1s\int g^k(X^*)\xi_k\left[1+ \frac {s^2}{(\tilde v')^2}\right] \frac {\partial}{\partial w}(y^i\varphi(Y))\, dY, \\
I_{24} &=\frac {2s(p'\cdot\xi')}{(\tilde v')^3}\int g(X^*) 
\frac {\partial}{\partial w}(w\varphi(Y))\, dY;
\end {align*} 
and
\begin {align*}
I_{31} &= - \frac {p\cdot\xi}{\tilde v}\int g_j(X^*)y^j h_1'(\tilde v)\frac {\partial}{\partial q_k}(q^k\varphi (Y)) \, dY, \\
I_{32} &= -\frac {p'\cdot\xi'}\int g_z(X^*)w h_2'(\tilde v) \frac {\partial}{\partial q_k}(q_k\varphi (Y)) \, dY, \\
I_{33} &= - \frac 1s\int g^m(X^*)\xi_m \frac {\partial}{\partial q_k}(q_k\varphi (Y)) \, dY.
\end {align*}
The integrals $I_{11}$, $I_{12}$, $I_{14}$, $I_{21}$, $I_{22}$, $I_{24}$, $I_{31}$, and $I_{32}$ are estimated directly, and the integrals $I_{13}$, $I_{23}$, and $I_{33}$ are estimated using the same decomposition as for $J_{13}$.
We therefore obtain
\[
|\tilde g_s^k\xi_k| \le C\left( \frac {|\xi'| }s+ \frac {|\xi\cdot\gamma|}{v}\right)
\]
if $\varepsilon\le 1$.

The other second derivatives involving $p$ are straightforward.
To estimate $g_{px}$, we  compute, using integration by parts
\[
\tilde g_i =- \frac 1{sh_1(\tilde v)}\int g(X^*)\frac {\partial }{\partial y^i}\varphi(Y)\, dY+ \frac 1{(\tilde v')^2}\int g(X^*)D_i(c^{km})p_kp_m\frac {\partial} {\partial w}(w\varphi(Y))\, dY,
\]
and hence
\[
\tilde g^k_i\xi_k\eta^i = J_{41}+ J_{42}+J_{43}+ J_{44} +J_{51} + J_{52} +J_{53}+ J_{54} +J_{55},
\]
with 
\begin {align*}
J_{41} &= \frac {h_1'(\tilde v)p\cdot \xi}{sh_1(\tilde v)^2\tilde v}\int g(X^*)\eta\cdot \frac {\partial \varphi(Y)}{\partial y}\, dY, \\
J_{42}&= - \frac {h_1'(\tilde v)p\cdot\xi}{h_1(\tilde v)\tilde v}\int 
g_j(X^*)y^j \eta\cdot \frac {\partial \varphi(Y)}{\partial y}\, dY, \\
J_{43} &= - \frac {h_2'(\tilde v')p'\cdot \xi'}{h_1(\tilde v)\tilde v'}\int
g_z(X^*)\eta\cdot \frac {\partial \varphi(Y^*)}{\partial y}\, dY, \\
J_{44} &= - \frac 1{sh_1(\tilde v)}\int g^k(Y^*)\xi_k\eta\cdot \frac {\partial \varphi(Y)}{\partial y}\, dY;
\end {align*} 
and
\begin {align*}
J_{51} &= -\frac {2p\cdot \xi}{(\tilde v')^3\tilde v} \int g(X^*)
\eta^iD_i(c^{km})p_kp_m\frac {\partial}{\partial w}(w\varphi(Y))\, dY, \\
J_{52} &= \frac {sp\cdot\xi}{(\tilde v')^2\tilde v} \int g_j(X^*)y^jh_1'(\tilde v)\eta^iD_i(c^{km})p_kp_m\frac {\partial}{\partial w}(w\varphi(Y))\, dY, \\
J_{53} &=\frac {sp'\cdot\xi'}{(\tilde v')^3} \int g_z(X^*) h_2'(\tilde v')\eta^iD_i(c^{km})p_kp_m\frac {\partial}{\partial w}(w\varphi(Y))\, dY, \\
J_{54} &= \frac 1{(\tilde v')^2} \int g^r(X^*)\xi_r\eta^iD_i(c^{km})p_kp_m \frac {\partial}{\partial w}(w\varphi(Y^*))\, dY, \\
J_{55} &= -\frac 2{(\tilde v')^2} \int g(X^*)\eta^iD_i\gamma^{k}p\cdot\gamma\xi_k\frac {\partial}{\partial w}(w\varphi(Y^*))\, dY, \\
J_{56} &= -\frac 2{(\tilde v')^2} \int g(X^*)\eta^iD_i\gamma^{m})p_m\xi\cdot\gamma\frac {\partial}{\partial w}(w\varphi(Y^*))\, dY.
\end {align*}
After integrating $J_{41}$ by parts, direct estimation shows that
\[
|\tilde g^k_i\xi_k\eta^i| \le C\theta_x(\tilde v)\tilde v^2\frac 1s|\eta||\xi|.
\]
Similar arguments give
\[
|\tilde g^k_z\xi_k| \le \frac {C}{\varepsilon'}\beta_1\tilde v'\left( \frac {|\xi'|}{s}+ \frac {|\xi\cdot \gamma|}{ v} \right).
\]

The remaining second derivatives are estimated using similar arguments.  
After integrating $I_1$, $I_2$, and $I_3$ by parts, we find that
\[
|\tilde g_{ss}| \le \frac {C}s, \quad |\tilde g_{sx}| \le \frac {C\theta_x(\tilde v)\tilde v^2}s, \quad |\tilde g_{sz}| \le \frac {C\beta_1\tilde v'}s.
\]

In the same vein, we have
\begin {gather*}
|\tilde g_{xx}| \le \frac {C\theta_x(\tilde v)^2\tilde v^4}{\varepsilon'|s|}, \quad |\tilde g_{xz}| \le \frac {C\beta_1\theta_x(\tilde v)\tilde v^2\tilde v'}{\varepsilon'|s|}, \\
 |\tilde g_{zz}| \le \frac {C\beta_1^2(\tilde v')^2}{\varepsilon'|s|}.
\end {gather*}
Here, the constant $C$ in the estimate of $\tilde g_{xx}$ also depends on $\Omega$, specifically, on the $C^3$ nature of $\partial\Omega$.

From our estimates for the derivatives of $\tilde g$, we obtain estimates  for the second derivatives of $N$.

First, by direct computation,
\[
N^{km}= \frac {-\tilde g^{km}-\varepsilon^{1/2}\tilde g^k_sN^m-\varepsilon^{1/2}\tilde g^m_sN^k-\varepsilon\tilde g_{ss}\rho^k\rho^m}{1+\tilde g_s\varepsilon^{1/2}},
\]
with $\tilde g$ and its derivatives now evaluated at $(x,z,p,\varepsilon^{1/2}N)$.  Recalling also \eqref {ENv}, we infer that there is a positive constant $c_3$, determined only by $\beta_0$, $c_0$, and $n$, such that
\[
|NN^{km}\xi_k\xi_m| \le c_3\left( \frac {|\xi'|^2}{\varepsilon'\varepsilon^{1/2}}+ \varepsilon'(\xi\cdot\gamma)^2 \right).
\]
If we take 
\[
\varepsilon'=\min\left\{ \frac 1{20\beta_0}, \frac 1{4c_3}\right\},
\]
we have \eqref {E10.372pp} provided $c_1\ge 2c_3/\varepsilon'$.

Next, we recall that
\[
N^kN^m\xi_k\xi_m = \frac {(\xi\cdot\gamma-\varepsilon^{1/2}\tilde g_p\cdot \xi)^2}{(1+\varepsilon^{1/2}\tilde g_s)^2}.
\]
Since $|\tilde g_s|\le \frac 12$, it follows that  
\[
\frac 1{(1+\varepsilon^{1/2}\tilde g_s)^2}\ge 1-2\varepsilon.
\]
Also
\[
(\xi\cdot\gamma-\varepsilon^{1/2}\tilde g_p\cdot \xi)^2\ge (\xi\cdot\gamma)^2-2\varepsilon^{1/2}|\xi\cdot\gamma||\tilde g_p\cdot\xi|,
\]
The Cauchy-Schwarz inequality and \eqref {Etildegp} then give
\begin {align*}
N^kN^m\xi_k\xi_m &\ge (1-2\varepsilon)[(1-20\beta_0\varepsilon)(\xi\cdot\gamma)^2-10\beta_0|\xi'|^2] \\
&\ge (1-22\beta_0\varepsilon)(\xi\cdot\gamma)^2- 22\beta_0|\xi'|^2.
\end {align*}
In combination with \eqref {E10.372pp}, this inequality implies \eqref {E10.43}.
Estimates \eqref {E10.372pz}, \eqref {E10.372px}, \eqref {E10.372zz}, \eqref {E10.372xz}, and \eqref {E10.372xx} are proved by similar arguments.
\end {proof}

\section {Some preliminary calculations} \label {Sp}

The basic idea in the proof of the gradient estimate is to examine a quadratic function of the gradient of the solution.  As first seen in \cite {IUMJ88}, the function is $c^{km}D_kuD_mu+\varepsilon N(x,u,Du;\varepsilon)^2$ for a suitably small $\varepsilon$, and this function has been used in several special circumstances as well. Here, we want to introduce a suitable change of variables to more closely mimic the gradient estimates in \cite {Serrin} (and subsequently in Chapter 15 of \cite {GT}).  
The combination of the more complicated quadratic function and the more general structure conditions leads to messier calculations, so we start here with some basic calculations which will be used in the next few sections to prove our gradient bound.

We begin by introducing an increasing $C^3$ function $\Psi$, defined on some interval which includes the range of $u$, and we write $\psi$ for the inverse to $\Psi$. To simplify the writing, we also use two standard bits of notation.  We set $\bar u = \Psi\circ u$ and $\omega= \psi''/(\psi')^2$.

We also define $w_1$ by
\begin {equation} \label {Ew1def}
w_1 = c^{km}D_m\bar u D_k\bar u + \frac {\varepsilon N(x,u,Du)^2}{(\psi')^2},
\end {equation}
we use the vector $\bar\nu_1$ defined by
\begin {equation} \label {Ebarnudef}
\bar \nu_1= \frac 1{\psi'}\nu_1,
\end {equation}
and we set
\begin {equation} \label {ESdef}
\mathscr S = a^{ij}[c^{km}+ \varepsilon(N^kN^m+NN^{km})]D_{ik}\bar u D_{jm}u.
\end {equation}
It is helpful to notice from \eqref {E10.43} that
\begin {equation} \label {ES}
\mathscr S \ge \frac 12[a^{ij}c^{km}D_{ik}\bar u D_{jm}\bar u + \varepsilon a^{ij}\gamma^k\gamma^m D_{ik}\bar u D_{jm} \bar u].
\end {equation}

Our first step is to compute the gradient of $w_1$.  A simple calculation gives
\begin {equation} \label {ED1w}
D_iw_1 = 2D_{ik}\bar u\bar\nu_1^k +D_i(c^{km})D_k\bar u D_m\bar u + \omega I_0D_i u + 2\frac {\varepsilon}{(\psi')^2}N[N_zD_iu+N_i].
\end {equation}
with
\begin {equation} \label {EI0}
I_0 = \frac {2\varepsilon N(N^kD_ku-N)}{(\psi')^2}.
\end {equation}
Because the exact expression for the second derivatives is quite involved, we jump directly to the main expression of interest, which is $a^{ij}D_{ij}w_1$. A long, tedious, but standard calculation shows that
\begin {equation} \label {ED2w}
\begin {split} 
a^{ij}D_{ij}w_1 &= 2a^{ij}D_{ijk}\bar u \bar \nu_1^k+ 2 \mathscr S +\frac {\omega'}{\psi'} I_0\mathscr E + (\omega^2 A_0  + \omega B_0 +C_0)\mathscr E \\
&\phantom{=\ }{}+\omega (S_0+ I_0a^{ij}D_{ij} u) +S_1+  \frac {2\varepsilon}{(\psi')^2} NN_za^{ij}D_{ij}\bar u,
\end {split}
\end {equation}
with $I_0$ given by \eqref {EI0},
\begin {align*}
A_0 &= 2\varepsilon NN^{km}D_kuD_mu+ 2\varepsilon(N-N^kD_ku)^2, \\
B_0 &= -\frac {4\varepsilon}{(\psi')^2\mathscr E}(N^kD_ku-N)a^{ij}D_iuN_j + \frac {4\varepsilon}{(\psi')^2}(N^kD_ku-N) N_z \\
&\phantom{=}\ {}+\frac {4\varepsilon}{(\psi')^2\mathscr E} NN^k_iD_ku a^{ij}D_ju + \frac {4\varepsilon}{(\psi')^2} NN^k_zD_ku, \\
C_0 &=   \frac {2\varepsilon}{(\psi')^2\mathscr E}a^{ij}(N_iN_j+NN_{ij}) + \frac {4\varepsilon}{(\psi')^2\mathscr E} [N_iN_z+NN_{iz}]a^{ij}D_ju  \\
&\phantom{=}\ {}+ \frac {2\varepsilon[(N_z)^2+NN_{zz}]}{(\psi')^2} +\frac {a^{ij}D_{ij}(c^{km})D_k\bar u D_m\bar u}{\mathscr E}, \\
S_0 &=  \frac {4\varepsilon}{\psi'} (N^kD_ku-N)N^m a^{ij}D_{im}\bar uD_ju +4\varepsilon NN^{km}a^{ij}D_{ik}\bar u D_juD_mu, \\
S_1&= 2a^{ij}D_i(c^{km})D_{jm}\bar u D_k\bar u +\frac {4\varepsilon}{\psi'} (N_zN^k+NN_z^k)a^{ij}D_{ik}\bar uD_ju \\
&\phantom{=\ }{}+ \frac {4\varepsilon}{\psi'} (N^kN_i+NN^k_i)a^{ij}D_{jk} \bar u.
 \end {align*}

We now estimate these terms.  First, we have 
\[
A_0 \ge  2\varepsilon NN^{km}D_kuD_mu,
\]
so \eqref {E10.39} and \eqref {E10.372pp} imply that
\begin {subequations} \label {E0ests}
\begin {gather}
A_0\ge -c w_1, \\
\intertext {where, here and in the remainder of this section, we use $c$ to denote any constant determined only by $\beta_0$, $c_0$, $n$, and $\Omega$. From \eqref {E10.37c}, \eqref {E10.372dot},  \eqref {E10.39}, \eqref {E10.372pz}, and \eqref {E10.372px}, we conclude that}
B_0 \ge -c\left[ \theta_x(v)v\left(\frac \Lambda{\mathscr E}\right)^{1/2}+ \beta_1\varepsilon^{1/2}\right]w_1. \\
\intertext {We use \eqref {E10.37b}, \eqref {E10.37c}, \eqref {E10.37d}, \eqref {E10.39}, \eqref {E10.372pz}, \eqref {E10.372zz}, \eqref {E10.372xz},   and \eqref {E10.372xx} to conclude that}
C_0 \ge -cw_1\left[\beta_1(1+ \varepsilon_x(v)v) \left(\frac {\Lambda}{\mathscr E}\right)^{1/2} +(1+ \varepsilon_x(v)v)^2 \frac \Lambda{\mathscr E} + \beta_1\varepsilon +  \frac {\Lambda}{\mathscr E\varepsilon}\right]. \label {Ec0est} \\
\intertext {The Cauchy-Schwarz inequality, \eqref {E10.372dot}, and \eqref {E10.372pp} imply that}
S_0 \ge - c\varepsilon^{1/2}(\Lambda w_1)^{1/2}\mathscr S^{1/2}. \\
\intertext {Because $D_i(c^{km})= -\gamma ^mD_i\gamma^k-\gamma^kD_i\gamma^m$, we infer from the Cauchy inequality, \eqref {E10.372pp}, \eqref {E10.43}, \eqref {E10.37c}, \eqref {E10.372pz}, and \eqref {E10.39} that}
S_1 \ge- cw_1^{1/2}\left( [\varepsilon^{-1/2}+ \varepsilon_x(v)v]\Lambda + \varepsilon^{1/2}\beta_1\mathscr E\right)^{1/2}\mathscr S^{1/2}. \\
\intertext {Finally, \eqref {E10.372dot} and \eqref {ENv} imply that}
|I_0| \le \frac {c\varepsilon vv_\varepsilon}{(\psi')^2}, \\
\intertext {and hence}
|I_0| \le \frac {c\varepsilon vw_1^{1/2}}{\psi'}, \label {Ei0est1} \\
|I_0| \le c\varepsilon^{1/2}w_1. \label {Ei0est}
\end {gather}
\end {subequations}

Next, we note (see (15.17) from \cite {GT}) that the differential equation \eqref {EPDE} is equivalent to
\begin {equation} \label {E15.17}
\psi' a^{ij}(x,u,Du)D_{ij}\bar u + a(x,u,Du) + \omega \mathscr E(x,u,Du)=0.
\end {equation}
If we apply the operator $\nu_1^kD_k$ to this equation and then add $D\bar u\cdot \bar \nu_1[\omega(r+1)+s]$ times \eqref {E15.17} for functions $r$ and $s$ to be further specified, we obtain (compare with equation (15.22) of \cite {GT}):
\begin {equation} \label {E15.22}
\begin {split}
0 &= [a^{ij}D_{ijk}\bar u +\kappa^iD_{ik}\bar u] \bar\nu_1^k + \frac {\omega'}{\psi'}\mathscr E D\bar u \cdot \nu_1 \\
& \phantom{=}{}+ (\omega^2A +  \omega B +C)\mathscr E D\bar u\cdot \bar\nu_1 + \omega  S_2 +S_3,
\end {split}
\end {equation}
with
\begin {align*}
\kappa^i &= \psi'a^{jk,i}D_{jk}\bar u +a^i+\omega\mathscr E^i, \\
A &=\frac {(\delta_1+r)\mathscr E}{\mathscr E}, \\
B &= \frac {(\delta_1+r)a+ (\delta_2+s)\mathscr E}{\mathscr E}, \\
C &= \frac {(\delta_2+s)a}{\mathscr E} \\
S_2&=[\psi'(\delta_1+r+1)a^{ij}D_{ij}\bar u]D\bar u \cdot \bar\nu_1, \\
S_3 &= [\psi'(\delta_2+s)a^{ij}D_{ij}\bar u]D\bar u\cdot \bar\nu_1,
\end {align*}
and the differential operators $\delta_1$ and $\delta_2$ are defined by
\[
\delta_1f(x,z,p) = p\cdot f_p(x,z,p), \quad
\delta_2f(x,z,p)= f_z(x,z,p)+\frac { f_k(x,z,p)\nu_1^k}{p\cdot\nu_1}.
\]
We defer the estimates for these terms to later sections because these estimates depend on the structure conditions for the differential equation.

By also calculating $\kappa^iD_iw_1$,  we find that
\begin {equation} \label {Ew}
\begin {split}
a^{ij}D_{ij}w_1+ \kappa^iD_iw_1 &= 
 D\bar u\cdot \bar\nu_1\mathscr E\left[\frac {\omega'}{\psi'}(I_1-1) + \omega^2A_2 +\omega B_2+C_2\right] \\
&\phantom {=\ } {}+2\mathscr S  + \omega S_4+ S_5 + \omega I_0a^{ij}D_{ij}u,
\end {split}
\end {equation}
with
\begin {align*}
I_1 &= \frac {I_0}{D\bar u\cdot\bar\nu_1}, \\
A_2 &= (-1+I_1)A +\frac {A_0}{D\bar u \cdot\bar \nu_1}, \\
B_2 &= - B+ \frac {B_0}{D\bar u\cdot\bar\nu_1}+I_1B' -\frac {2\mathscr E^iD_i\gamma^kD_ku Du\cdot\gamma}{Du\cdot\nu_1\mathscr E} \\
&\phantom {=\ }{}+\frac {2\varepsilon NN_zA}{Du\cdot\nu_1}+ \frac {2\varepsilon NN_i\mathscr E^i}{Du\cdot\nu_1\mathscr E}, \\
C_2 &= -C +\frac {C_0}{D\bar u\cdot\bar\nu_1}+ \frac {2\varepsilon Na^iN_i}{Du\cdot\nu_1\mathscr E}- \frac {2a^iD_i\gamma_kD_kuDu\cdot\gamma}{Du\cdot\nu_1\mathscr E} + \frac {2B' \varepsilon NN_z}{Du\cdot\nu_1}, \\
S_4 &=S_0 +(-1+I_1)S_2, \\
S_5 &= S_1-S_3 + \frac {2\varepsilon NN_z(\delta_1+r+1)a^{ij}D_{ij}\bar u}{\psi'}\\
&\phantom{=\ }{} - 2a^{jk,i}D_{jk}\bar u D_i\gamma^mD_m\bar uDu\cdot\gamma - \frac {2\varepsilon NN_ia^{jk,i}D_{jk}\bar u}{\psi'},
\end {align*}
and
\[
B'=  \frac {(\delta_1+r)a}{\mathscr E}.
\]
The remainder of this paper is concerned with deriving a gradient bound under various hypotheses modeled on those of Serrin \cite {Serrin}. With 
\begin {subequations} \label {Elimsup}
\begin {align}
A_\infty &= \limsup_{|p|\to\infty} \sup_{(x,z)\in\Omega\times\mathbb R} A(x,z,p), \\
B_\infty &= \limsup_{|p|\to\infty} \sup_{(x,z)\in\Omega\times\mathbb R} B(x,z,p), \\
C_\infty &= \limsup_{|p|\to\infty} \sup_{(x,z)\in\Omega\times\mathbb R} C(x,z,p),  \label {Elimsupc}
\end {align}
\end {subequations}
and $\nu$ replaced by $Du/|Du|$, Serrin derived a gradient bound in  four cases:  $A_\infty\le 0$, $C_\infty\le0$,  $B_\infty\le -\sqrt{A_\infty C_\infty}$, and the oscillation of $u$ is sufficiently small. These four cases are exactly those for which the differential equation 
\[
\frac {dy}{dt}= A_\infty y^2 + B_\infty y+C_\infty +\eta
\]
has a solution on the range of $u$ for $\eta$ a sufficiently small positive constant.  Unfortunately, there are some important difficulties in trying to translate the full argument in \cite {Serrin} to our situation.
First Serrin uses a decomposition 
\[
a^{ij}= a^{ij}_*+p_ic_j+p_jc_i,
\]
with $[a^{ij}_*]$ a uniformly elliptic matrix and $c$ a convenient vector-valued function.  For the oblique derivative problem, the corresponding decomposition would be
\[
a^{ij}=a^{ij}_* + \nu_1^ic_j+\nu_1^jc_i,
\]
and this decomposition depends in a complicated way on the parameter $\varepsilon$ and on the function $b$.  In addition, our control on the term $A_0$ is not good enough to handle the case $A_\infty=0$, for example.
On the other hand, we can consider several of the critical examples from \cite {Serrin}.  Specifically, we shall consider two important cases:  $C_\infty\le0$ and the small oscillation case.
Moreover, in our situation, it is possible that $B_\infty$ and $C_\infty$ depend on $\varepsilon$.  We will therefore need to take this fact into account.

\section {Global gradient estimates} \label {S1GB}

We now turn to our gradient estimates.  First, we prove them in a neighborhood of $\partial\Omega$, assuming that a bound is already known away from the boundary.  Such bounds are well-known (see, for example, Theorem 3 from \cite {Serrin} or Theorem 15.3 from \cite {GT}).
In this case, our estimate is quite straightforward under suitable structure conditions. To state our result simply, we define
\begin {equation} \label {EA'B'}
B'_\infty=\limsup_{|p|\to\infty} \sup_{(x,z)\in\Omega\times\mathbb R}  |B'|, \quad 
A'_\infty =\limsup_{|p|\to\infty} \sup_{(x,z)\in\Omega\times\mathbb R}  |A|.
\end {equation}
For a positive constant $M_0$, we also write $\Gamma(M_0)$ for the set of all $(x,z,p)\in \Omega\times\mathbb R\times \mathbb R^n$ with $|p| \le M_0$.

Our first estimate assumes that $C_\infty\le0$.

\begin {theorem} \label {TCinfty}
Let $u\in C^2(\clOmega) $ be a solution of \eqref {EODP} with $\partial\Omega\in C^3$ and $b$ satisfying the hypotheses of Theorem \ref {TN} for some $*$-decreasing function $\varepsilon_x$.
 Suppose that there are functions $r$ and $s$, nonnegative constants $\mu_2$ and $M_0$, two decreasing functions $\tilde\mu$ and (for each $\varepsilon\in(0,\varepsilon_0)$) $\tilde\mu_\varepsilon$ with 
\begin {subequations} \label {Etildelimit}
\begin {gather}
\lim_{\sigma\to\infty} \tilde\mu(\sigma)=0, \label {Etildelimitmu} \\ 
\lim_{\sigma\to\infty} \tilde\mu_\varepsilon(\sigma)=0, \label {Etildelimiteps}
\end {gather} 
\end {subequations}
and a $*$-decreasing function $\tilde\mu_*$ such that
\begin {subequations} \label {SC1aij}
\begin {gather}
(\delta_1+r+1)a^{ij}\eta_{ij} \le \frac {\tilde \mu(v)}{v} \mathscr E^{1/2}(a^{ij}\delta^{km}\eta_{ij}\eta_{jm})^{1/2}, \label {SC1aij1}\\
(\delta_2+s)a^{ij}\eta_{ij} \le \frac {\tilde \mu_\varepsilon(v)}{v} \mathscr E^{1/2}(a^{ij}\delta^{km}\eta_{ij}\eta_{jm})^{1/2}, \label {SC1aij2}\\
|a^{ij}_p\eta_{ij}| \le \frac {\tilde\mu_*(v)}{v} \mathscr E^{1/2}(a^{ij}\delta^{km}\eta_{ij}\eta_{jm})^{1/2} \label {SC1aijp}
\end {gather}
\end {subequations}
on $\Gamma(M_0)$ for all matrices $[\eta_{ij}]$ and that
\begin {subequations} \label {Elambdap}
\begin {gather} 
|a| \le \mu_2\mathscr E, \label {EaE} \\
\Lambda \le \tilde\mu_*(v)^2\mathscr E,  \label {ElambdaE1} \\
 |a_p| \le \tilde\mu_*(v)\mathscr E  \label {SC1p}
\end {gather}
\end {subequations}
on $\Gamma(M_0)$.  Suppose also that  $B'_\infty$ is bounded uniformly with respect to $\varepsilon$ and that $C_\infty\le0$ for all $\varepsilon>0$.  If 
\begin {subequations} \label {Edelta1b}
\begin {gather}
\tilde\mu_*(v)\delta_1b \le \tilde\mu(v)b_p\cdot\gamma, \label {Edelta1b1} \\
\tilde\mu_*(v)\delta_2b \le \tilde\mu(v)b_p\cdot\gamma \label {Edelta1b2}
\end {gather}
\end {subequations}
on $\Sigma_0(\tau_0)$, and if
\begin {equation} 
\lim_{\sigma\to\infty} (1+\varepsilon_x(\sigma)\sigma)\tilde\mu_*(s)=0, \label {Etildelimit1}
\end {equation}
then there is a constant $M$, determined only by $\beta_0$, $\beta_1$, $c_0$, $M_0$, $\tau_0$, $n$, $\tilde\mu_*$, $\tilde\mu$, $\Omega$, $\sup_{\{d\ge R_0/4\}}|Du|$, the oscillation of $u$, and the limit behavior in \eqref {Elimsup}, \eqref {EA'B'}, \eqref {Etildelimit}, and \eqref {Etildelimit1}, such that $|Du| \le M$ in $\Omega$.
\end {theorem}
\begin {proof} First, we assume additionally that $u\in C^3(\Omega)$, and we note from the differential equation for $u$ along with \eqref {Ei0est} that
\[
\omega I_0a^{ij}D_{ij}u =-\omega I_0a \ge - c\omega \varepsilon^{1/2}w_1|a|.
\]
(Again, we use $c$ to denote any constant determined only by $\Omega$, $\beta_0$, and $c_0$.)
From \eqref {Ew}, \eqref {SC1aij}, and \eqref {EaE}, we conclude that
\[
a^{ij}D_{ij}w_1+\kappa^iD_iw_1 \ge \frac 32\mathscr S +D\bar u\cdot\bar\nu_1\left( \frac {\omega'}{\psi'}(I_1-1)+\omega^2 A_3 +\omega B_3+C_3\right),
\]
with 
\begin {align*}
A_3 &= A_2-c\left(\varepsilon \frac {\Lambda}{\mathscr E} + \frac {\tilde\mu(v)^2}{\varepsilon}\right), \\
B_3 &= B_2- c\varepsilon\frac {|a|}{\mathscr E}, \\
C_3 &= C_2 -\frac {c}{\varepsilon}\left[( 1 + \varepsilon_x(v)v)^2\frac {\Lambda}{\mathscr E} +\beta_1^2\varepsilon^2 +\tilde\mu_\varepsilon(v) +\tilde\mu_*(v)^2 (1+\varepsilon_x(v)v)^2\right].
\end {align*}  

We now observe that
\[
(\delta_1+r)\mathscr E = (\delta_1+r+1)a^{ij}p_ip_j +\mathscr E,
\]
so \eqref {SC1aij1} with $\eta_{ij}=p_ip_j$ implies that 
\[
A_\infty = A_\infty'=1.
\]
By invoking \eqref {E10.37c}, \eqref {E10.37d}, \eqref {ENv}, \eqref {E0ests}, \eqref {Elambdap} and \eqref {Etildelimit1}, we infer that there are constants $A_{3,\infty}$ and $c$, determined only by $\beta_0$, $c_0$, $n$, and $\Omega$, such that
\begin {align*}
\liminf_{|p|\to\infty} \inf_{(x,z)\in\Omega\times\mathbb R} A_3 &\ge - A_{3,\infty}, \\
\liminf_{|p|\to\infty} \inf_{(x,z)\in\Omega\times\mathbb R} B_3 &\ge -B_\infty- c\mu_2\varepsilon^{1/2}, \\
\liminf_{|p|\to\infty} \inf_{(x,z)\in\Omega\times\mathbb R} C_3 &\ge -c\beta_1^2\varepsilon.
\end {align*}

We now define the function $\chi$ by $\chi=\omega\circ \psi^{-1}$ and write $B_\infty^1$ for the uniform upper bound on $B_\infty$.  Then, as shown on p. 595 of \cite {Serrin}, there is a positive constant $\eta$, determined only by $A_{3,\infty}$ and $B_\infty^1$, such that the differential equation
\[
\chi'(z)= A_{3,\infty}\chi^2 +B_\infty^1 \chi +\eta
\]
has a solution on the range of $u$. With this choice for $\chi$ (which also gives the function $\psi$), we conclude that there are positive constants $M_1$ and $\varepsilon_1$ such that $w_1\ge M_1$ and $\varepsilon<\varepsilon_1$ imply that
\begin {equation} \label {Ew1}
a^{ij}D_{ij}w_1+\kappa^iD_iw_1 \ge \frac 32\mathscr S +\frac 12\eta D\bar u\cdot \bar \nu_1\mathscr E.
\end {equation}
We now fix $\varepsilon\in(0,\varepsilon_1)$.

For $k_1$ a positive constant to be chosen, we now introduce the function 
\[
w_2= w_1+k_1d\int_0^{w_1} \frac 1{\tilde\mu_*(\sqrt \sigma)}\, d\sigma.
\]
Straightforward calculation shows that
\begin {align*}
a^{ij}D_{ij}w_2 +\kappa^iD_iw_2 &= \left(1+\frac {k_1d}{\tilde\mu_*(\sqrt{w_1})}\right)(a^{ij}D_{ij}w_1+\kappa^iD_iw_1)  \\
&\phantom{+\ }{}+ k_1\int_0^{w_1}\frac 1{\tilde\mu_*(\sqrt \sigma)}\, d\sigma (a^{ij}D_{ij}d+\kappa^iD_id) \\
&\phantom{+\ }{} - \frac {k_1\tilde\mu_*'(\sqrt{w_1})d}{\sqrt{w_1}\tilde\mu_*(\sqrt{w_1})^2}a^{ij}D_iw_1D_jw_1 \\
&\phantom{=\ }{}+\frac {2k_1}{\tilde\mu_*(\sqrt{w_1})}a^{ij}D_idD_jw_1.
\end {align*}
Since $\tilde\mu_*'\le0$, we infer that
\[
a^{ij}D_{ij}w_2 +\kappa^iD_iw_2 \ge \frac 32\mathscr S +\frac 12\eta D\bar u\cdot \bar \nu_1\mathscr E +C_4+S_6
\]
with
\begin {align*}
C_4&=k_1\int_0^{w_1} \frac 1{\tilde\mu_*(\sqrt \sigma)}\, d\sigma(a^{ij}D_{ij}d +\omega \mathscr E^iD_id+a^iD_id) \\
&\phantom{=\ } {}-\frac {2k_1}{\tilde\mu_*(\sqrt{w_1})}( a^{ij}D_idD_j(c^{km})D_k\bar u D_m\bar u) \\
&\phantom{=\ }{} - \frac {2k_1}{\tilde\mu_*(\sqrt{w_1})}\left(\left[\omega I_0 +\frac {2\varepsilon}{(\psi')^2} NN_z\right] a^{ij}D_idD_ju + \frac {2\varepsilon} {(\psi')^2} Na^{ij}D_idN_j\right) \\
\intertext {and}
S_6 &=k_1\int_0^{w_1} \frac 1{\tilde\mu_*(\sqrt \sigma)}\, d\sigma\psi' a^{jk,i}D_{jk}\bar uD_id +\frac {2k_1}{\tilde\mu_*(\sqrt{w_1})}a^{ij}D_idD_{jk}\bar u\nu_1^k
\end {align*}
 wherever $w_1\ge M_1$.
 Since $\tilde\mu_*$ is decreasing, it follows that
 \[
 \int_0^{w_1} \frac 1{\tilde\mu_*(\sqrt \sigma)}\, d\sigma\le \frac {w_1}{\tilde\mu_*(\sqrt{w_1})}.
\]
In addition, because $\tilde\mu_*$ is $*$-decreasing, it follows that
\[
\frac 1{\tilde\mu_*(\sqrt\sigma)} \ge \frac {\sqrt\sigma}{\sqrt{w_1}\tilde\mu_*(\sqrt{w_1})}
\]
for $0\le\sigma\le w_1$ and hence
\[
\int_0^{w_1} \frac 1{\tilde\mu_*(\sqrt \sigma)}\, d\sigma\ge \frac {2w_1}{3\tilde\mu_*(\sqrt{w_1})}.
\]
Moreover, because $\tilde\mu_*$ is $*$-decreasing, we conclude that the fraction
\[
\frac {\tilde\mu_*(\sqrt w_1)}{\tilde\mu_*(v)}
\]
is bounded from above and below by positive constants, determined only by $A_\infty'$, $B_\infty'$,  $\beta_0$, $c_0$, $n$, $\varepsilon$, and $\tilde\mu_*$.

Our estimate of $C_4$ uses an estimate of $\mathscr E_p$, which we now derive.
For any vector $\xi$, we have
\[
\mathscr E^i\xi_i = a^{jk,i}p_jp_k\xi_i + 2 a^{ik}p_k\xi_i.
\]
Then \eqref {SC1aijp} implies that
\[
a^{jk,i}p_jp_k\xi_i \le \frac {\tilde\mu_*(v)}v\mathscr E^{1/2}(a^{ij}\delta^{km}p_ip_kp_jp_m)^{1/2}|\xi| = \frac {\tilde\mu_*(v)}v\mathscr E|p||\xi|.
\]
The Cauchy-Schwarz inequality implies that
\[
2a^{ij}p_k\xi \le 2\mathscr E^{1/2}(a^{ij}\xi_i\xi_j)^{1/2} \le 2\mathscr E^{1/2} \Lambda^{1/2}|\xi|.
\]
From \eqref {ElambdaE1} and the inequality $|p|\le v$, it follows that
\[
\mathscr E^i\xi_i \le 3\tilde\mu_*(v) \mathscr E|\xi|,
\]
and hence
\begin {equation} \label {SCEp}
\left| \mathscr E_p\right| \le 3\tilde\mu_*(v) \mathscr E.
\end {equation}

We now apply \eqref {ElambdaE1}, \eqref {SC1p}, and \eqref {SCEp} to conclude that
\[
C_4\ge -D\bar u\cdot\bar\nu_1\mathscr E (ck_1),
\]
and we apply \eqref {SC1aijp} and \eqref {ElambdaE1} to conclude that
\[
S_6 \ge - ck_1(D\bar u\cdot\bar\nu_1)^{1/2}\mathscr E^{1/2}\mathscr S^{1/2}.
\]
(Here, $c$ is determined by all the quantities mentioned in the conclusion of this theorem.)
It follows that
\begin {multline} \label {ELw2}
a^{ij}D_{ij}w_2 + \kappa^iD_iw_2\ge \left(\frac \eta2\left[1+ k_1d\int_0^{w_1} \frac 1{\tilde\mu(\sqrt\sigma)}\ d\sigma\right]- c(k_1^2+k_1)\right) D\bar u\cdot \bar\nu_1\mathscr E \\
+\left(1+k_1d\int_0^{w_1} \frac 1{\tilde\mu(\sqrt\sigma)}\ d\sigma\right)\mathscr S.
\end{multline}
By taking $k_1$ sufficiently small, we conclude that
\[
a^{ij}D_{ij}w_2 + \kappa^iD_iw_2\ge 0
\]
on $E$, the subset of $\Omega_{R_0}$ on which $w_1\ge M_1$.

We also remove the assumption $u\in C^3$ by observing (see also equation (15.13) in \cite {GT}) that $w_2$ is a weak solution of the differential inequality
\[
D_i(a^{ij}D_jw_2) + [\kappa^i-D_j(a^{ij})]D_iw_2 \ge 0
\]
in $E$.
It then follows from Theorem 8.1 in \cite {GT} implies that $w_2$ attains its maximum over $\overline E$ on $\partial E$, which consists of three subsets:
\begin {align*}
E_1 &= \left\{x\in \partial E: d(x)=\frac {R_0}4\right\}, \\
E_2 &= \{x\in \partial E: w_1(x)=M_1\}, \\
E_3 &= \partial E\cap \partial \Omega.
\end {align*}

With $M_3=\sup_{\Omega_{R_0/4}} w_1$, it's straightforward to check that
\[
w_2 \le M_3 + \frac {k_1R_0}4\int_0^{M_3} \frac 1{\tilde\mu_*(\sqrt\sigma)}\, d\sigma
\]
on $E_1$, and we have an upper bound for $M_3$. Moreover,
\begin {equation} \label {Ew2E2}
w_2 \le M_1 + \frac {k_1R_0}4\int_0^{M_1} \frac 1{\tilde\mu_*(\sqrt\sigma)}\, d\sigma
\end {equation}
on $E_2$.

On $E_3$, we compute
\[
b^iD_iw_2= -(\omega\delta_1b +\delta_2b)w_1+b^iD_i(c^{km})D_k\bar uD_m\bar u+k_1\int_0^{w_1}\frac 1{\tilde\mu_*(\sqrt \sigma)}\, d\sigma b_p\cdot\gamma.
\]
We then invoke \eqref {Edelta1b}, \eqref {Etildelimit1}, and the first inequality of \eqref {E10.33} to conclude that there is a  positive constant $M_4$ such that 
\[
b^iD_iw_2>0
\]
at any point of $E_3$ where $w_1\ge M_4$.
Since $w_1=w_2$ on $E_3$, it follows that
\[
w_2\le \max\{M_1,M_4\}
\]
on $E_3$. It follows that $w_2\le c$ on $E$, and we have the upper bound \eqref {Ew2E2} on $\Omega_{R_0/4}\setminus E$ as well. Since $w_2\ge w_1$, we conclude that $w_1\le c$ on $\Omega_{R_0/4}$, so $|Du|\le c$ on $\Omega_{R_0/4}$.  In combination with our assumed upper bound for $|Du|$ on $\Omega\setminus \Omega_{R_0/4}$, we obtain the desired estimate.
\end {proof}

We remark that, although the structure condition \eqref {ElambdaE1} does not appear in the global estimates of Section 3 from \cite {Serrin} or Section 15.2 from \cite {GT} (except in the example of uniformly elliptic equations, where a stronger assumption is made), it is satisfied by many standard equations.  It also appears in the local gradient estimates of Chapter 15 from \cite {GT}, specifically in condition (15.47) there. We also remark that the argument given here does not rely on the precise form of the interior gradient estimate, unlike the argument in \cite {LT} and \cite {NUE}.

As we shall see in our parabolic examples, condition \eqref {SC1aij1} is quite restrictive.  In order to relax it, we need to strengthen the condition on $b_z$ just a little.  Even though this improvement will not be used for our elliptic examples, we include it here for completeness.

\begin {theorem} \label {TCinfty1}
Let $u\in C^2(\clOmega) $ be a solution of \eqref {EODP} with $\partial\Omega\in C^3$ and $b$ satisfying the hypotheses of Theorem \ref {TN} for some $*$-decreasing function $\varepsilon_x$.
 Suppose that there are functions $r$ and $s$, nonnegative constants $M_0$ and $\mu_2$, a decreasing function  $\tilde\mu_\varepsilon$ (for each $\varepsilon\in(0,\varepsilon_0)$) satisfying \eqref {Etildelimiteps}, and a $*$-decreasing function $\tilde\mu_*$ such that \eqref {SC1aij2}, \eqref {SC1aijp}, and 
\begin {equation} \label {SC2aij}
(\delta_1+r+1)a^{ij}\eta_{ij} \le \frac { \mu_2}{v} \mathscr E^{1/2}(a^{ij}\delta^{km}\eta_{ij}\eta_{jm})^{1/2}, 
\end {equation}
 hold on $\Gamma(M_0)$ for all matrices $[\eta_{ij}]$ and that \eqref {Elambdap} holds on $\Gamma(M_0)$.  Suppose also that  $B'_\infty$ is bounded uniformly with respect to $\varepsilon$ and that $C_\infty\le0$ for all $\varepsilon>0$.  If  there is a decreasing function $\tilde\mu$ satisfying \eqref {Etildelimitmu} such that \eqref {Edelta1b} holds on $\Sigma_0(\tau_0)$ and \eqref {Etildelimit1} is satisfied and if
 \begin {equation} \label {Elimbz}
 \lim_{|p|\to\infty} \frac {b_z(x,z,p)}{b_p(x,z,p)\cdot\gamma}=0,
 \end {equation}
then there is a constant $M$, determined only by $\beta_0$, $c_0$, $M_0$, $\tau_0$, $\tilde\mu_*$, $\tilde\mu$, $\Omega$,  the oscillation of $u$,$\sup_{\{d\ge R_0/4\}}|Du|$, the limit behavior in \eqref {Elimsup}, \eqref {EA'B'}, \eqref {Etildelimit1}, \eqref {Etildelimiteps}, and \eqref {Elimbz} such that $|Du| \le M$ in $\Omega$.
\end {theorem}
\begin {proof} 
By using \eqref {Elimbz} in the proof of Theorem \ref {TN}, we conclude that there is a decreasing function $\varepsilon_z$, determined only by the limit behavior in \eqref {Elimbz} with $\lim_{\sigma\to\infty} \varepsilon_z(\sigma)=0$ such that
\begin {align*}
|N_z| &\le 4(1+c_0)\varepsilon_z(v_\varepsilon), \\
|NN_{pz}\cdot\xi| &\le c_1\varepsilon_z(v_\varepsilon) \left( \frac 1{\varepsilon^{1/2}}|\xi'|^2+ |\xi\cdot\gamma|\right)v_\varepsilon \\
\intertext {for any vector $\xi$,}
|NN_{zz}| &\le c_1\beta\varepsilon_z(v_\varepsilon)v_\varepsilon^2, \\
|NN_{xz}| &\le \frac {c_1}{\varepsilon^{1/2}}\varepsilon_z(v_\varepsilon)\left( \frac v{R_0}+ \varepsilon_x(v)v^2\right).
\end {align*}

With these improved estimates, we see from \eqref {Ew} and \eqref {SC1aij} that
\[
a^{ij}D_{ij}w_1+\kappa^iD_iw_1 \ge \frac 32\mathscr S +D\bar u\cdot\bar\nu_1\left( \frac {\omega'}{\psi'}(I_1-1)+\omega^2 A'_3 +\omega B'_3+C'_3\right),
\]
with 
\begin {align*}
A'_3 &= A_2-c\left(\varepsilon \frac {\Lambda}{\mathscr E} + \frac {\mu_2}{\varepsilon}\right), \\
B'_3 &= B_2, \\
C'_3 &= C_2 -\frac {c}{\varepsilon}\left[( 1 + \varepsilon_x(v)v)^2\frac {\Lambda}{\mathscr E} +\beta_1\varepsilon_z(v_\varepsilon)\varepsilon^2 +\tilde\mu_\varepsilon(v) +\tilde\mu_*(v)^2 (1+\varepsilon_x(v)v)^2\right].
\end {align*}  
From the proof of Theorem \ref {TCinfty}, we see that there is a constant $A_{3,\infty}^*$, determined by $\mu_2$ and the same quantities as for $A_{3,\infty}$ such that
\begin {align*} \liminf_{|p|\to\infty} \inf_{(x,z)\in\Omega\times\mathbb R} A'_3 &\ge - \frac {A_{3,\infty}^*}\varepsilon, \\
\liminf_{|p|\to\infty} \inf_{(x,z)\in\Omega\times\mathbb R} B_3 &\ge -B_\infty- c\varepsilon^{1/2}, \\
\liminf_{|p|\to\infty} \inf_{(x,z)\in\Omega\times\mathbb R} C'_3 &\ge 0.
\end {align*}

Just as before, there is a positive constant $\eta$, determined only by $A_{3,\infty}^*$ and $B^1_\infty$, such that the differential equation 
\[
Y'= A_{3,\infty}^*Y^2+B^1_\infty Y+\eta
\]
has a solution on the range of $u$.
If we now take $\chi=Y/\varepsilon$, we have
\[
\chi' = \frac {A_{3,\infty}^*}\varepsilon\,\chi^2 +B^1_\infty\chi +\eta\varepsilon,
\]
and therefore there are positive constants $M_1$ and $\varepsilon_1$ such that the inequality $\varepsilon<\varepsilon_1$ implies that \eqref {Ew1} holds wherever $w_1\ge M_1$.  The remainder of the proof is identical to that of Theorem \ref {TCinfty}.
\end {proof}

In fact, we can weaken condition \eqref {Elimbz} slightly.  It suffices that $\beta_1$ be sufficiently small.

When $C_\infty<0$, we can take $\psi(s)=s$ in the proof of Theorem \ref {TCinfty}.  Hence, a number of hypotheses can be removed or weakened in this case.  For brevity, we just state the result.

\begin {theorem} \label {CCinfty} 
Let $u\in C^2(\clOmega)$ be a solution of \eqref {EODP} with  $\partial\Omega\in C^3$ and $b$ satisfying the hypotheses of Theorem \ref {TN} for some $*$-decreasing function $\varepsilon_x$.  Suppose that there are functions $r$ and $s$, a $*$-decreasing function $\tilde\mu_*$, a decreasing function $\tilde\mu_\varepsilon$ satisfying \eqref {Etildelimiteps}, and nonnegative constants $M_0$ and $\mu_2$ such that \eqref {SC1aij2} and \eqref {SC1aijp} are satisfied on $\Gamma(M_0)$ for all matrices $[\eta_{ij}]$ and that \eqref {ElambdaE1} and \eqref {SC1p}
hold on $\Gamma(M_0)$. Suppose also that $B'_\infty$ is uniformly bounded for $\varepsilon\in (0,\varepsilon_0]$, and that $C_\infty<0$.
If 
\begin {equation} \label {Edelta1beta}
 b_z \le \varepsilon_x(v)vb_p\cdot\gamma
\end {equation}
on $\Sigma_0(\tau_0)$, and if \eqref {Etildelimit1} holds, then there is a constant $M$, determined only by $\beta_0$, $c_0$, $M_0$, $n$, $\tilde\mu_*$, $\mu_2$, $\tilde\mu$, $\tau_0$, $\Omega$, $\sup_{\{d\ge R_0/4\}}|Du|$, and the limit behavior in \eqref {Etildelimiteps}, \eqref {Elimsupc}, and \eqref {Etildelimit1}, such that $|Du| \le M$ in $\Omega$.
\end {theorem}

When the oscillation of $u$ is sufficiently small, then we can derive a gradient bound as long as the quantities $A_\infty$, $A'_\infty$, $B_\infty$, $B'_\infty$, and $C_\infty$ are bounded from above.  In fact, the upper bounds may depend on $\varepsilon$.

\begin {theorem} \label {TGBsmall}
Let $u\in C^2(\clOmega)$ be a solution of \eqref {EODP} with $\partial\Omega\in C^3$ and $b$ satisfying the hypotheses of Theorem \ref {TN} for some $*$-decreasing function $\varepsilon_x$.   Set $\varepsilon = \varepsilon_0/2$, and suppose that there are functions $r$ and $s$, a $*$-decreasing function $\tilde\mu_*$, and nonnegative constants $M_0$ and $\mu_2$ such that conditions \eqref {SC1aijp}, \eqref {SC2aij},  and
\begin {equation} \label {SCaij262}
  (\delta_2+s)a^{ij}\eta_{ij}\le \frac {\mu_2}v \mathscr E^{1/2}(a^{ij}\delta^{km}\eta_{ik}\eta_{jm})^{1/2}
\end {equation}
hold on $\Gamma(M_0)$ for all matrices $[\eta_{ij}]$ and that conditions \eqref {Elambdap} hold on $\Gamma(M_0)$.  Suppose also  that $B'_\infty$ and $C_\infty$ are finite for each $\varepsilon \in (0,\varepsilon_0)$.  If
\begin {equation} \label {Evarep}
\limsup_{\sigma\to\infty} \varepsilon_x(\sigma)\sigma\tilde\mu_*(\sigma)<\infty,
\end {equation} 
and if 
\begin {equation} \label {Edelta1beta1}
\tilde\mu_*(v) \max\{\delta_1b, \delta_2b\} \le  vb_p\cdot\gamma
\end {equation}
on $\Sigma_0(\tau_0)$, then there are constants $M$ and $\omega_0$, determined only by $\beta_0$, $\beta_1$, $\beta_2$, $c_0$, $M_0$,  $\mu_2$, $\tau_0$, $\Omega$, $\sup_{\{d \ge R_0/4\}}|Du|$, and the limit behavior in \eqref {Elimsup}, \eqref {EA'B'}, and \eqref {Evarep}, such that $|Du|\le M$ in $\clOmega$ provided $\osc_{\Omega} u \le \omega_0$.
\end {theorem}
\begin {proof} We first note that the proof of the estimate for $A_\infty'$ from Theorem \ref {TCinfty} shows that $A_\infty'\le c$ by virtue of \eqref {SC2aij}.
 As stated in the hypotheses of this theorem, we take $\varepsilon=\varepsilon_0/2$.  With $k_1$ to be determined, we take $w_2$ as in the proof of Theorem \ref {TCinfty}. Then we can choose $k_1$ so that $b^iD_iw_2>0$ on the subset of $\partial\Omega$ on which $w_1\ge M_0$ and $\omega\le 1$.

This choice of $k_1$ gives a constant $c$ such that
\[
a^{ij}D_{ij}w_2+\kappa^iD_iw_2 \ge D\bar u\cdot \nu_1\mathscr E\left( \frac {\omega'}{\psi'}-c\right)
\]
on the subset of $\Omega_{R_0/4}$ on which $w_1\ge M_0$ and $\omega\le1$.
We now choose $\omega$ so that $\chi'=c+1$ and $\chi(\inf u)=0$.  It follows for $\omega_0 =1/(c+1)$ that $\omega\le1$ and hence by our previous arguments, we obtain the desired gradient bound.
\end {proof}

Note that, if $\tilde\mu_*(s)=K/s$ for some positive constant $K$ (which will be the case in all of our examples), then condition \eqref {Evarep} is automatically satisfied and \eqref {Edelta1beta} and \eqref {Edelta1beta1} follow from the condition
\begin {equation} \label {Ebeta11}
|b_x|+v|b_z| \le \beta_3v^2b_p\cdot\gamma
\end {equation}
for some nonnegative constant $\beta_3$.

\section {Local gradient estimates} \label {Sso}

It is also possible to give local estimates.  To present them in a more compact format, we introduce some further notation.

Specifically, for some $y\in\partial\Omega$ and some $R>0$, we look at solutions of 
\begin {subequations} \label {EODPloc}
\begin {gather}
a^{ij}(x,u,Du)D_{ij}u+a(x,u,Du) =0 \text { in } \Omega \cap B(y,R), \\
b(x,u,Du)=0 \text { on } \partial\Omega\cap B(y,R).
\end {gather}
\end {subequations}
Roughly speaking, we can estimate the gradient of $u$ at $y$ if $\tilde\mu_*$ is a power function with suitable negative exponent.
 We also set
\begin {subequations} \label {ERinfty}
\begin {align}
A'_{\infty,R} &= \limsup_{|p|\to\infty} \sup_{(x,z)\in\Omega\cap B(y,R)\times\mathbb R} |A(x,z,p)|, \\
B_{\infty,R} &= \limsup_{|p|\to\infty} \sup_{(x,z)\in\Omega\cap B(y,R)\times\mathbb R} |B(x,z,p)|, \\
B'_{\infty,R} &= \limsup_{|p|\to\infty} \sup_{(x,z)\in\Omega\cap B(y,R)\times\mathbb R} |B'(x,z,p)|, \\
C_{\infty,R} &= \limsup_{|p|\to\infty} \sup_{(x,z)\in\Omega\cap B(y,R)\times\mathbb R} C(x,z,p), \label {ERinftyc}
\end {align}
\end {subequations}
and we write $\Gamma_{R}(M_0)$ for the set of all $(x,z,p)\in \Gamma(M_0)$ with $|x-y|<R$.
We also write $\Sigma_0(\tau,R)$ for the set of all $(x,z,p)\in \Sigma_0(\tau)$ with $|x-y|<R$.

Our local gradient estimate takes the following form.

\begin {theorem} \label {TCloc} Let $u\in C^2(\clOmega\cap B(y,R))$ be a solution of \eqref {EODPloc} with $\partial\Omega\cap B(y,R)\in C^3$ for some $R\in (0,R_0)$ and $b$ satisfying the hypotheses of Theorem \ref {TN} with $\partial\Omega\cap B(y,R)$ in place of $\partial\Omega$ for some $*$-decreasing function $\varepsilon_x$.  Suppose that there are constants $\theta\in(0,1]$, $M_0>0$, and $\mu_3\ge0$ such that 
\begin {equation} \label {SC1aijploc} 
v^{1+\theta}|a^{ij}_p\eta_{ij}| \le \mu_3 \mathscr E^{1/2}(a^{ij}\delta^{km}\eta_{ik}\eta_{jm})^{1/2}
\end {equation}
on $\Gamma_R(M_0)$ for all matrices $[\eta_{ij}]$, and
\begin {subequations} \label {ElambdaE2}
\begin {gather}
v^{2\theta}\Lambda \le \mu_3^2\mathscr E,\label {ElambdaE2a} \\
v^\theta|a_p| \le \mu_3\mathscr E \label {ElambdaE2b}
\end {gather}
\end {subequations}
on $\Gamma_R(M_0)$.  
\begin {enumerate}
\item  \label {TClocle}
 Suppose also that there are functions $r$ and $s$, a nonnegative constant $\mu_2$, a decreasing function $\tilde\mu$ and, for each $\varepsilon\in(0,\varepsilon_0)$, a decreasing function $\tilde\mu_\varepsilon$ satisfying \eqref {Etildelimit} such that \eqref {SC1aij1} and  \eqref {SC1aij2} hold on $\Gamma_R(M_0)$ for all matrices $[\eta_{ij}]$ and \eqref {EaE} holds on $\Gamma_R(M_0)$,  and suppose that $B'_{\infty,R}$ is bounded uniformly with respect to $\varepsilon$  and that $C_{\infty,R}\le0$. If 
\begin {equation} \label {Ethetalimit}
\lim_{\sigma\to\infty} \sigma^{1-\theta}\varepsilon_x(\sigma)=0
\end {equation}
and if
\begin {equation} \label{Ethetadelta}
\max\{\delta_1b, \delta_2b\} \le v^{\theta}\tilde\mu(v)b_p\cdot\gamma
\end {equation}
on $\Sigma_0(\tau_0,R)$, then there is a constant $M$, determined only by $\beta_0$, $\beta_1$, $c_0$, $M_0$, $n$, $\theta$, $\mu_3$, $R$, $\tau_0$, $\Omega$, the oscillation of $u$ over $\Omega\cap B(y,R)$, and the limit behavior in \eqref {Etildelimit}, \eqref {ERinfty}, and \eqref {Ethetalimit}, such that $|Du(y)| \le M$.
\item  \label {TCloclez}
 Suppose also that there are functions $r$ and $s$, a nonnegative constant $\mu_2$,  and, for each $\varepsilon\in(0,\varepsilon_0)$, a decreasing function $\tilde\mu_\varepsilon$ satisfying \eqref {Etildelimiteps} such that \eqref {SC1aij2} and  \eqref {SC2aij} hold on $\Gamma_R(M_0)$ for all matrices $[\eta_{ij}]$ and \eqref {EaE} holds on $\Gamma_R(M_0)$, and suppose that $B'_{\infty,R}$ is bounded uniformly with respect to $\varepsilon$  and that $C_{\infty,R}\le0$. If \eqref {Ethetalimit} holds, if \eqref {Elimbz} holds, and if \eqref {Ethetadelta} is satisfied on $\Sigma_0(\tau_0,R)$ for some decreasing function $\tilde\mu$ satisfying \eqref{Etildelimitmu}, then there is a constant $M$, determined only by $\beta_0$, $\beta_1$, $c_0$, $M_0$, $n$, $\theta$, $\mu_3$, $R$, $\tau_0$, $\Omega$, the oscillation of $u$ over $\Omega\cap B(y,R)$, and the limit behavior in \eqref {Etildelimit}, \eqref {Elimbz}, \eqref {ERinfty}, and \eqref {Ethetalimit},  such that $|Du(y)| \le M$.
\item \label {TClocless}
 Suppose also there are functions $r$ and $s$ and  a nonnegative constant $\mu_2$ such that \eqref {SC2aij}  holds on $\Gamma_R(M_0)$ for all matrices $[\eta_{ij}]$ and $C_{\infty,R}<0$.  If \eqref {Ethetalimit} holds and if \eqref {Edelta1beta} is satisfied on $\Sigma_0(\tau_0,R)$, then there is a constant $M$, determined only by $\beta_0$, $\beta_1$, $c_0$, $M_0$, $n$, $\mu_2$, $\mu_3$, $R$, $\Omega$, $\tau_0$, and the limit behavior in \eqref {ERinftyc}, \eqref {SC1aijploc}, and \eqref {Ethetalimit}, such that $|Du(y)| \le M$.
\item \label {TClocso}
  Suppose also that there are functions $r$ and $s$ and a nonnegative constant $\mu_2$ such that conditions \eqref {EaE}, \eqref {SC2aij}, and \eqref {SCaij262} hold on $\Gamma_R(M_0)$ for all matrices $[\eta_{ij}]$, and suppose  $B'_{\infty,R}$ and $C_{\infty,R}$ are finite for each $\varepsilon$.  If
\begin {equation} \label {Elimsupvarep}
\limsup_{\sigma\to\infty} \sigma^{1-\theta}\varepsilon_x(\sigma)<\infty,
\end {equation}
and if 
\begin {equation} \label {Edeltabtheta}
\max\{\delta_1b, \delta_2b\} \le \mu_3v^\theta b_p\cdot\gamma
\end {equation}
on $\Sigma_0(\tau_0,R)$, then there are constants $M$ and $\omega_0$, determined only by $\beta_0$, $\beta_1$, $c_0$, $M_0$, $\mu_2$, $\Omega$, $\tau_0$, and the limit behavior in \eqref{EA'B'}, \eqref {ERinfty} and \eqref {Elimsupvarep},  such that $|Du(y)|\le M$  provided $\osc_{\Omega\cap B(y,R)} u \le \omega_0$.
\end {enumerate}
\end{theorem}
\begin {proof} To prove part \eqref {TClocle}, we use the notation from the proof of Theorem \ref {TCinfty} (with $\mu_3s^{-\theta}$ in place of $\tilde\mu_*(s)$); in particular, we take $\varepsilon$, $\psi$, and $M_1$ from the proof of the theorem, and we assume initially that $u\in C^3(\Omega\cap B(y,R))$.  From \eqref {ELw2}, we conclude that there is a positive constant $\eta'$ such that
\[
a^{ij}D_{ij}w_2+\kappa^iD_iw_2 \ge \left(1+ k_2dw_1^{\theta/2}\right)[\mathscr S + \eta'w_1\mathscr E]
\]
for
\[
k_2 = \frac {k_1\theta\mu_3}2
\]
wherever $v\ge M_1$.

Next, the discussion on pages 346 and 347 of \cite {ODPbook} gives  a positive constant $R_1$, determined only by $R_0$ and $\beta_0$, and, for each $R\in (0,R_1)$, a function $\zeta\in C^2(\mathbb R^n)$ such that $\zeta(y)=3/4$, $\zeta\le0$ outside $B(y,R/2)$, and $b_p\cdot D\zeta\ge 0$ at any point of $\partial\Omega\cap B(y,R/2)$ at which $\zeta\ge0$.  Further, there is a constant $c$, determined only by $\beta_0$ and $n$ such that $R|D\zeta|+R^2|D^2\zeta|\le c$. We therefore assume without loss of generality that $R<R_1$.

We now set $q= 1+2/\theta$ and $w= \zeta^qw_2$ and note that
\begin {align*}
D_iw_2 &= \zeta^{-q}D_iw -\zeta^{-q}w_2D_i(\zeta^{-q}) \\
&= \zeta^{-q}D_iw-q\zeta^{-1}w_2D_i\zeta.
\end {align*}
It follows that
\[
a^{ij}D_{ij}w+\kappa_1^iD_iw = \zeta^{1+2/\theta}(a^{ij}D_{ij}w_2+\kappa^iD_iw_2) + C_5w\mathscr E +S_7,
\]
with
\begin {align*}
\kappa_1^i &=\kappa^i-  2q\zeta^{-1}a^{ij}D_j\zeta, \\
C_5 &= \frac 1{\mathscr E} \left( -q^2\zeta^{-2}a^{ij}D_i\zeta D_j\zeta + q\zeta^{-1}a^{ij}D_{ij}\zeta + \zeta^{q-1}[\omega\mathscr E^i+a^i]D_i\zeta\right), \\
 \intertext {and}
S_7 &= q\psi'w_2\zeta^{q-1}a^{jk,i}D_{jk}\bar u.
\end {align*}
The  proof of \eqref {SCEp} shows that $|\mathscr E_p| \le 3\mu_3\mathscr E$, and hence, by invoking \eqref {E10.37c}, \eqref {E10.37d}, \eqref {Ei0est}, \eqref {SC1aijploc},  \eqref {Ethetalimit}, \eqref {ElambdaE2},  and the Cauchy-Schwarz inequality, we conclude that there is a constant $c$ for which
\[
a^{ij}D_{ij}w+\kappa_1^iD_iw \ge  w\mathscr E\left(\tilde\eta - \frac c {\zeta v^{\theta} R}- \frac c{(\zeta v^\theta R)^2}\right).
\]
From the definition of $w_2$, it follows that there is a constant $c$, determined only by the quantities in the conclusion of this theorem (except for $R$), such that
\[
w_2\le cv^2(1+Rv^\theta) \le cv^{2+\theta},
\]
if we further assume (again, without loss of generality) that $R\le1$.
Therefore,
\[
w\le c\zeta^qv^{2+\theta} = c(\zeta v^\theta)^{(2+\theta)/\theta}.
\]
For any positive constant $M_2$, it follows that
\[
\frac 1 {\zeta v^{\theta} R}+ \frac 1{(\zeta v^\theta R)^2} \le c(M_2^{-\theta/(2+\theta)} + M_2^{-2\theta/(2+\theta)})
\]
wherever $w \ge M_2$, and hence, just as in the proof of Theorem 3 of \cite {Serrin} or Theorem 15.3 in \cite {GT}, if we choose $M_2$ sufficiently large, determined also by $R$, then
\[
a^{ij}D_{ij}w + \kappa_1^iD_iw >0
\]
on $E$, the subset of $\Omega\cap B(y,R/2)$ where $w\ge M_2$ and $\zeta>0$. Just as in the proof of Theorem \ref {TCinfty}, it follows that $w\le M_2$ in $E$ or $w$ attains its maximum on $E^0= \partial\Omega\cap E$, even if $u\in C^2(\Omega\cap B(y,R/2))$.

If $w$ attains its maximum on $E^0$, then at that point we have
\[
b^iD_iw= w_2\frac 2{\theta}\zeta^{(2-\theta)/\theta}b^iD_i\zeta + \zeta^{2/\theta}b^iD_iw_2.
\]
Since $b^iD_i\zeta\ge0$, it follows from the proof of Theorem \ref {TCinfty} that $w$ cannot attain its maximum at a point of $E_0$ with $w\ge M_3$ for a suitable constant $M_3$, and hence we obtain the estimate $w\le c$ in $\Omega\cap B(y,R/2)$.  In particular, we obtain an upper bound for $w(y)$ and hence for $|Du(y)|$.

Parts \eqref {TCloclez}, \eqref {TClocless}, and \eqref {TClocso} are proved in a similar fashion.
\end {proof}

Of course,  \eqref {Ebeta11} implies \eqref {Elimsupvarep} and \eqref
{Edeltabtheta} if $\theta=1$.

\section {examples} \label {Se}

\subsection {Capillary-type boundary conditions} \label {Sscbc}
To begin, we look at some boundary functions $b$ which satisfy our structure conditions.  We suppose that there are $C^1$ scalar functions $h$  (defined on $[1,\infty)$) and $\psi$ (defined on $\partial\Omega\times\mathbb R$) such that
\begin {equation} \label {Ebh}
b(x,z,p)=h(v)p\cdot\gamma +\psi(x,z).
\end {equation}
We assume that $h$ is positive and that there is a positive constant $h_0$ such that
\begin {subequations} \label {Ehpsi}
\begin {gather} 
-\frac {h(\sigma)}\sigma \le h'(\sigma) \le h_0h(\sigma)^2 \label {Ehprime} \\
\intertext {for all $\sigma\in(1,\infty)$ and}
\lim_{\sigma\to\infty} \sigma h(\sigma) > \sup_{\partial\Omega\times\mathbb R} |\psi|. \label {Elimpsi}
\end {gather}
\end {subequations}
By virtue of the first inequality in \eqref {Ehprime}, the limit in \eqref {Elimpsi} exists.  Moreover, the latter condition is only a restriction on $\psi$ if the limit is finite.
The capillary problem is the special case $h(\sigma)=1/\sigma$, in which case \eqref {Elimpsi} requires $\sup|\psi|<1$.

We compute
\begin {gather*}
b_p= h(v)\gamma +h'(v) p\cdot\gamma \frac {p}{v}, \\
b_p\cdot \gamma = h(v)+h'(v)\frac {(p\cdot\gamma)^2}v.
\end {gather*}
The first inequality in \eqref {Ehprime} implies that
\[
b_p\cdot\gamma \ge h(v)\left[ 1- \left(\frac {p\cdot\gamma}v\right)^2\right]
\]
and hence $b_p\cdot\gamma>0$.  Moreover,
\[
b_p\cdot\gamma = h(v)\left[ 1-\left (\frac \psi{vh(v)}\right)^2\right]
\]
wherever $b=0$. It follows from \eqref {Elimpsi} (and the monotonicity of the function $\sigma \mapsto \sigma h(\sigma)$) that there are positive constant $\tau_0$ and $\Psi$ for which
\[
1- \left(\frac \psi{vh(v)}\right)^2\ge \Psi
\]
when $v\ge \tau_0$ and hence
\[
b_p\cdot \gamma \ge \Psi h(v)
\]
when $v\ge \tau_0$. 
Moreover, if $\lim_{s\to\infty} sh(s)=\infty$, then we may take $\Psi\in (0,1)$ arbitrary provided $\tau_0$ is sufficiently large.

When $b(x,z,p)=0$, we have
\[
|b_p(x,z,p)| \le h(v)+|h'(v)||p\cdot\gamma| = h(v) + \frac {h'(v)\psi}{h(v)} \le h(v)\left[1 + \frac {\sup|\psi|h'(v)}{h(v)^2}\right],
\]
so we infer the first inequality of \eqref {E10.33} with $\beta_0= (1+ h_0 \sup|\psi|)/\Psi$ by also using the second inequality of \eqref {Ehprime}.
When $b(x,z,p)=0$ and $v\ge \tau_0$, we have
\[
\frac {|p\cdot\gamma|}v = \frac {|\psi(x,z)|}{vh(v)} \le (1-\Psi)^{1/2}.
\]
Simple algebra then yields the second inequality of \eqref {E10.33} with $c_0=( (1-\Psi)/\Psi)^{1/2}$.
Similar computations yield \eqref {E10.35} with
\[
\varepsilon_x(\sigma) = \left(\frac {\sup|\psi_x|}{h(1)\Psi}+ \frac 2{R_0\Psi}\right) \frac 1\sigma,
\]
and $\beta_1= \sup |\psi_z|/h(1)$.
We also observe that \eqref {Elimbz} is satisfied if $\lim_{\sigma\to\infty} \sigma h(\sigma)=\infty$.

Since $\theta>0$, \eqref {Ethetalimit} and \eqref {Elimsupvarep} are immediate.
Furthermore, \eqref {Etildelimit1} follows from \eqref {E10.35} if
\begin {equation} \label {Etildelimit11}
\lim_{\sigma\to\infty} \tilde\mu_*(\sigma)=0
\end {equation}
because $\varepsilon_x(s)s$ is a bounded function of $s$, while \eqref {Evarep} and \eqref {Ebeta11} with $\beta_3=1+\beta_1$ also follow from \eqref {E10.35}.

The other conditions are more delicate.
Since
\[
\delta_1b = (p\cdot\gamma)\left[h(v) + h'(v) \frac {|p|^2}v\right] = -\psi\left(1+ \frac {h'(v)|p|^2}{h(v)v}\right)
\]
and
\[
1+ \frac {h'(v)|p|^2}{h(v)v}\ge 1- \frac {|p|^2}{v^2}>0,
\]
it follows that $\delta_1b\le0$ whenever $\psi\ge0$.
We defer further discussion of the other conditions to the examples of differential equations.

\subsection {Other boundary conditions I} \label {SSOBC1}
It is possible to generalize the previous example slightly by allowing some $x$ and $z$ dependence in the gradient term of $b$.  Specifically, we suppose that there are a positive-definite matrix valued $C^1$ function $\beta^{ij}$ and a positive $C^1([1,\infty))$ function $h$ such that $b$ has the form
\begin {subequations} \label {Ebhtilde}
\begin {gather}
b(x,z,p)= h(\tilde v)\beta^{ij}(x,z)p_i \gamma_j +\psi \\
\intertext {with}
\tilde v= \left( 1+ \beta^{ij}(x,z)p_ip_j\right)^{1/2}.
\end {gather}
\end {subequations}
We further assume that $h$ satisfies the conditions \eqref {Ehprime} for some positive constant $h_0$ and
\begin {equation} \label {Elimpsi1}
\lim_{\sigma\to\infty} \sigma h(\sigma) > \sup \frac {|\psi|}{\lambda_0^{1/2}},
\end {equation}
where $\lambda_0(x,z)$ is the minimum eigenvalue of the matrix $[\beta^{ij}(x,z)]$.  
By imitating the arguments in the previous example, we see that the hypotheses of Theorem \ref {TN} are satisfied with $\beta_0$ and $c_0$ determined by $h_0$, $\sup_{\partial\Omega\times \mathbb R}|\psi|$, the function $\beta^{ij}$, and the quantities in \eqref {Elimpsi1}; $\beta_1$ is determined also by $\sup|\psi_z|$ and $\sup|\beta^{ij}_z|$; and
\[
\varepsilon_x(\sigma) = \frac {K_1}\sigma
\]
for some constant $K_1$ determined by the quantities in \eqref {Elimpsi1}, $h_0$, $\sup |\psi_x|$, and the function $\beta^{ij}$.  The second inequality of \eqref {E10.33} is proved by using the proof of Lemma 2.2(b) from \cite {JDE49}, specifically, the demonstration there that $1-(\nu\cdot\gamma)^2$ is bounded away from zero. It's easy to see that condition (2.2) from \cite {JDE49} is valid with constants $\beta_1, \tau\ge1$ and $\beta_2\in(0,1)$.
Moreover, \eqref {Elimbz} holds if $\lim_{\sigma\to\infty} \sigma h(\sigma)=\infty$.

\subsection {Other boundary conditions II} \label {SSOBC2}

We can also consider boundary conditions with $b$ of the form
\begin {equation} \label {Ebq}
b(x,z,p)= v^{q(x)}p\cdot\gamma +\psi(x,z)
\end {equation}
with $q$ a $C^1$ function satisfying the inequality $q(x)\ge-1$, and $\sup_{\mathscr Q\times\mathbb R} |\psi|<1$, where $\mathscr Q$ is the set on which $q=-1$.  Since the special case $q\equiv-1$ has already been dealt with, we assume here that $q\not\equiv-1$.
  This time, the hypotheses of Theorem \ref {TN} are satisfied for suitable constants $\beta_0$, $\beta_1$ and $c_0$, along with  $\varepsilon_x$ of the form
\[
\varepsilon_x(\sigma) = \frac {K_1\ln \sigma}\sigma 
\]
for $K_1$ a positive constant.  Specifically, the constants $\beta_0$ and $c_0$ are determined by 
\[
\sup_{\mathscr Q\times\mathbb R}|\psi|,
\]
the $C^1$ nature of $q$ at  points of $\mathscr Q$, and the maximum of $q$; $\beta_1$ is determined also by $\sup |\psi_z|$; and $K_1$ is determined also by the $C^1$ nature of $q$ everywhere and $\sup|\psi_x|$. 
This time, \eqref {Etildelimit1} holds if
\[
\lim_{\sigma\to\infty} \tilde\mu_*(\sigma)\ln \sigma=0,
\]
\eqref {Elimbz} holds if $q>-1$, and \eqref {Ethetalimit}, \eqref {Ebeta11}, and \eqref {Elimsupvarep} always hold.  Moreover \eqref {Edelta1beta} holds with $\beta_1$ determined by the $C^1$ nature of $q$ at points of $\mathscr Q$, the maximum of $q$, and $\sup|\psi_z|$.

We infer \eqref {Edelta1b} (for a suitable $\tilde\mu$) if
\begin {equation} \label {Etildeq}
\lim_{\sigma\to\infty} v^{-q}(\sigma)\tilde\mu_*(\sigma)=0,
\end {equation}
while \eqref {Edelta1beta} follows from the weaker condition
\begin {equation} \label {Etildeq1}
\limsup_{\sigma\to\infty} v^{-q}(\sigma)\tilde\mu_*(\sigma)<\infty.
\end {equation}
We infer \eqref {Elimbz} provided $\inf q>-1$.
If $q+\theta>0$, then \eqref {Ethetadelta} holds, while $q+\theta\ge0$ implies  \eqref {Edeltabtheta}.  Finally, we infer \eqref {Evarep} from the limit condition
\begin {equation} \label {Etildemuln}
\lim_{\sigma\to\infty} \tilde\mu_*(\sigma)\ln \sigma= 0.
\end{equation}

\subsection {Non-variational boundary conditions}
All of our boundary conditions so far correspond to the natural boundary condition for a variational problem.  In other words, they have the form
\[
b(x,z,p)= F_p(x,z,p)\cdot\gamma
\]
for some $C^2$ function $F$ which is convex with respect to $p$.  Specifically,
\[
F(x,z,p)=\int_1^v \sigma h(\sigma)\,d\sigma +\psi(x,z)p\cdot\gamma 
\]
for \eqref {Ebh},
\[
F(x,z,p)=\int_1^{\tilde v} \sigma h(\sigma)\,d\sigma + \psi(x,z)p\cdot\gamma
\]
for \eqref {Ebhtilde}, and
\[
F(x,z,p) = \frac {v^{q(x)+2}}{q(x)+2} +\psi(x,z)p\cdot\gamma
\]
for \eqref {Ebq}.  This form allows the possibility of deriving gradient bounds by applying the maximum principle to $F(x,u,Du)$ (assuming suitable structure conditions on $F$, $a^{ij}$, and $a$), an idea going back to Ural$'$tseva \cite {UR0} and Lieberman \cite {JDE49} although both authors only studied the conormal problem in which $a^{ij}= \partial^2F/(\partial p_i\partial p_j)$. We now provide a class of boundary conditions not having this form but for which the paper of the current page provides gradient estimates.

We start with a unit vector valued function $\beta$, defined on $\partial\Omega\times\mathbb R$, such that $\beta\cdot\gamma \ge \beta_*$ for some constant $\beta_*\in(0,1)$.  In particular, $\beta=\gamma$ is such a function for any such $\beta_*$, and if $\beta$ is a unit vector with $\beta\cdot\gamma\ge1$, then $\beta\cdot\gamma\equiv1$ and $\beta=\gamma$. The function $b$ now has the form
\begin {subequations} \label {Ebetastar}
\begin {gather}
b(x,z,p)= h(v)\beta\cdot p +\psi(x,z) \label {Ebetastar1} \\
\intertext {for some positive function $h$ with}
-\beta_* \le \frac {\sigma h'(\sigma)}{h(\sigma)} \le \frac {2\beta_*}{1-\beta_*} \label {Ebetastarh}
\end {gather}
\end {subequations}
for all $\sigma\ge1$.

We first show that $b$ is oblique.  To this end, we compute
\[
b_p\cdot\gamma = h(v)\beta\cdot \gamma + h'(v)v(\beta\cdot \nu()\gamma \cdot \nu).
\]
We also note that
\[
|\beta-\gamma|^2= |\beta|^2+|\gamma|^2- 2\beta\cdot\gamma\ge 2-2\beta_*,
\]
and hence
\begin {align*}
(\beta\cdot p\nu)(\gamma\cdot \nu) &= (\beta\cdot \nu)^2+ (\beta\cdot \nu) \nu\cdot(\beta-\gamma) \ge (\beta\cdot \nu)^2 -|\beta\cdot\nu||\nu|\beta-\gamma| \\
&\ge (\beta\cdot \nu)^2 -|\beta\cdot \nu||\nu|(2-2\beta_*)^{1/2}.
\end {align*}
The Cauchy-Schwarz inequality implies that
\[
|\beta\cdot \nu||\nu(2-2\beta_*)^{1/2} \le (\beta\cdot \nu)^2 - \frac 14|\nu|^2(2-2\beta_*)
\]
and therefore
\begin {equation} \label {Ebpdot}
(\beta\cdot \nu)(\gamma\cdot \nu) \ge - \frac {1-\beta_*}2|\nu|^2.
\end {equation}

If $h'(v)\le0$, then 
\[
h'(v)v(\beta\cdot\nu)(\gamma\cdot\nu)  \ge -\beta_*h(v)
\]
by virtue of the first inequality in \eqref {Ebetastarh}, so $b_p\cdot\gamma>0$ in this case.

On the other hand, if $h'(v)>0$, then
\[
h'(v)v(\beta\cdot \nu)(\gamma\cdot \nu)
\ge -\frac {2\beta_*}{1-\beta} h(v)\frac {1-\beta_*}2|\nu|^2 >-\beta_*h(v).
\]
It follows in either case that $b_p\cdot\gamma>0$.
For large $v$, we can obtain a lower bound for $b_p\cdot\gamma$ in line with the previous lower bounds.  First, we write
\[
h'(v)v(\beta\cdot \nu)(\gamma\cdot\nu) = -\frac {h'(v)v(\gamma\cdot\nu)}{h(v)} \psi(x,z) \ge -\,\frac { 2\beta_*}{1-\beta_*}\frac {|\psi|}v.
\]
The lower bound in \eqref {Ebetastar} implies that $\lim_{\sigma\to\infty} \sigma h(\sigma)=\infty$, so there is a $\tau_0\ge1$, determined only by $v(1)$, $\beta_*$, and $\sup|\psi|$, such that $v\ge\tau_0$ implies that
\[
vh(v) \ge \frac {4\sup|\psi|}{1-\beta_*}
\]
and simple algebra shows that $b_p\cdot\gamma \ge (\beta_*/2)h(v)$ for $v\ge\tau_0$.  
For this reason, we can actually relax the upper bound in \eqref {Ebetastarh} somewhat to
\[
 h'(\sigma) \le \begin {cases} \dfrac {2\beta_*h(\sigma)}{(1-\beta_*)\sigma} &\text { if } \sigma\le\tau_0, \\
h_0h(\sigma)^2 &\text { if } \sigma>\tau_0
\end {cases}
\]
for some $h_0\ge0$.  We leave the details of the analysis for this more extended class of boundary conditions to the reader.
We do point, however, that the hypotheses of Theorem \ref {TN} are satisfied with constants $\beta_1$ and $c_0$ determined by $v(1)$, $\beta_*$, and $\sup|\psi|$ and $\varepsilon_x(\sigma) =K/\sigma$ for some positive constant $K$ determined only by $v(1)$, $\beta_*$, $\sup|\psi|$, and $\Omega$.

Although we could modify Examples \ref {SSOBC1} and \ref {SSOBC2} in a similar fashion, we leave the details of that modification to the interested reader.

We also make a quick comparison between this example and Example \ref {SSOBC1}.  Namely, the expression $\beta^{ij}p_i\gamma_j$ appears in Example \ref {SSOBC1} while the slightly different expression $\beta^ip_i$ appears in this example.  Of course, given a matrix $\beta^{ij}]$ as in Example \ref {SSOBC1}, the vector $\beta$ given by $\beta^i=\beta^{ij}\gamma_j$ makes these two expressions equal. (Although, in order for $\beta$ to be a unit vector, we must define an intermediate vector  $\hat\beta$ by $\hat\beta^i=\beta^{ij}\gamma_j$ and then set $\beta=\hat\beta/|\hat\beta$.  Then $\beta^{ij}p_i\gamma_j$ and $\beta^ip_i$ differ only by a positive factor.)
 On the other hand, given a vector $\beta$ as in this example, we can find a positive-definite, symmetric matrix $[\beta^{ij}]$ making the expressions equal.
To create such a matrix, we assume that $\gamma=(0,\dots,1)$, with the general case being recovered through a simple change of variables.  
We then define
\[
\beta^{ij} = \begin {cases}  \frac 12\beta^n+ \frac 2{\beta^n}\sum_{j=1}^{n-1}(\beta^j)^2 &\text { if } i=j<n, \\
\beta^i &\text { if } j=n, \\
\beta^j &\text { if } i=n, \\
0 &\text {otherwise}.
\end {cases}
\]
It's elementary to check that, for this matrix, $\beta^{ij}\gamma_j=\beta^i$.
Moreover, direct calculation shows that, for any vector $\xi$, we have
\[
\beta^{ij}\xi_i\xi_j = \sum_{i=1}^{n-1}\beta^{ii}(\xi_i)^2 + \beta^n (\xi_n)^2 + 2\sum_{i=1}^{n-1}\beta^i\xi_i\xi_n,
\]
and a simple application of the Cauchy-Schwarz inequality to the last sum in this equation shows that
\[
\beta^{ij}\xi_i\xi_j \ge \frac 12\beta^n|\xi|^2,
\]
so $[\beta^{ij}]$ is positive definite.  It is symmetric by construction, and all of its entries are bounded from above by a constant determined only by $\beta_*$.
Therefore, we could have written this example directly in terms of a matrix $[\beta^{ij}]$ as in Example \ref {SSOBC1}.  We choose not to do so for ease of notation and to point out that this example applies to semilinear boundary conditions, which correspond to $h$ being constant, as well.

Now, we consider some classes of differential equations.

\subsection {The false mean curvature equation}

The model problem for our class of equations is known as the false mean curvature equation. In particular,
\[
a^{ij}= \delta^{ij}+p_ip_j
\]
and $a$ satisfies the conditions
\begin {subequations} \label {SCaFMCE}
\begin {gather}
 |a_x| = o(|p|^5), \\
a_z\le o(|p|^4), \\
|a|+ |p||a_p|=O(|p|^4). \label {SCaFMCEp}
\end {gather}
\end {subequations}
We now take $\tilde\mu_*(s)=K/s$ with $K$ a large constant, determined only by $n$ and the limit behavior in \eqref {SCaFMCEp}. 
Since $\delta_1a^{ij}=2p_ip_j$, conditions \eqref {SC1aij} are satisfied with $r=-3$, $s=0$, $\tilde \mu(\sigma)\ge 4/\sigma$ (the exact choice of $\tilde\mu$ will be made later), and $\tilde\mu_\varepsilon\equiv0$.
With these choices (in particular, with $K$ sufficiently large), \eqref {ElambdaE1}, \eqref {Etildelimit1}, and \eqref {SC1p} follow for any of our examples of boundary conditions.
It's easy to compute $B'_\infty=C_\infty=0$, while $B_\infty$ is bounded by a constant determined by the limit behavior in \eqref {SCaFMCEp}.  
Moreover, local gradient bounds follow from Example 2 on page 585 of \cite {Serrin} (with the constant $\theta$ in that example equal to $1$ and the multipler function $t$ in that example equal to $0$).

When $h$ has the form \eqref {Ebh}, then conditions \eqref {Edelta1beta} and \eqref {Edeltabtheta} are easy to check.  To verify \eqref {Edelta1b}, we consider two possibilities.  First, if 
\[
\lim_{\sigma\to\infty} \sigma h(\sigma)<\infty,
\]
in addition to the bound on $\sup|\psi|$, we assume that
\begin {equation} \label {Ehprime1}
\lim_{\sigma\to\infty}  \frac {h'(\sigma)\sigma}{h(\sigma)} =-1
\end {equation}
and that $\psi_z\le0$.  We then take
\[
\tilde \mu(\sigma) = K_1\left(\frac 1\sigma + 1+\frac {h'(\sigma)\sigma}{h(\sigma)}\right)
\]
for a sufficiently large positive constant $K_1$. Since we have already shown that $\delta_1b\le0$ wherever $\psi\ge0$, we only need to estimate $\delta_1b$ wherever $\psi<0$.  Moreover, it follows from \eqref {Ehprime1} that $h'(\sigma)<0$ for $\sigma$ sufficiently large, in which case, we have
\[
\delta_1b \le -\psi\left[1+ \frac {h'(v)v}{h(v)} +\frac 1{v^2}\right].
\]
From this inequality, \eqref {Edelta1b1} follows easily, and \eqref {Edelta1b2} is easily verified since $b_z=\psi_z\le0$. 

On the other hand, if
\[
\lim_{\sigma\to\infty} \sigma h(\sigma)=\infty,
\]
we assume that
\begin {equation} \label {Ehprime2}
\limsup_{\sigma\to\infty} \frac {h'(\sigma)}{h(\sigma)^2} \le0.
\end {equation}
This time, we take
\[
\tilde\mu(\sigma)= K\left(\frac1\sigma + \frac 1{\sigma h(\sigma)}+\max\{0,\frac {h'(\sigma)}{h(\sigma)^2}\}\right)
\]
to infer \eqref {Edelta1b}.

In particular, we require $\psi_z\le0$ when $h(\sigma) =1/\sigma$ (the capillary problem) or $h(\sigma)=(\arctan \sigma)/\sigma$.  If $h(\sigma)=\sigma^{q-1}$ for some $q>0$ or if $h(\sigma)=\exp(\sigma^\alpha)$ for some $\alpha\ge0$, our gradient estimate holds without any restriction on $\psi_z$.

When $b$ has the form \eqref {Ebhtilde}, we obtain a gradient estimate if $h$ satisfies the same restrictions as before and if $\beta^{ij}_z=0$.

When $b$ has the form \eqref {Ebq} and $q$ is nonconstant, we obtain a gradient estimate provided $\inf q>-1$ because of condition \eqref {Etildeq}. 

Finally, when $b$ has the form \eqref {Ebetastar1} with $h$ satisfying \eqref{Ebetastarh}, we obtain a gradient estimate provided $\beta$ and $\psi$ are $C^1$ with respect to $x$ and $z$.

In fact, the hypotheses of Theorem \ref {TCloc}\eqref {TClocle} are satisfied with $\theta=1$, so we obtain a local gradient bound for all these boundary conditions.

Note that this gradient bound improves the one for operators of this form in \cite {NUE} by relaxing the conditions on $a$. Theorem 4.2 of that work assumes that $|a| = O(|p|^3)$ and that $\delta_1a\le o(|p|^4)$ to obtain a local gradient estimate, while Theorem \ref {TCloc}\eqref {TClocle} allows $|a|=O(|p|^4)$ and $\delta_1a\le O(|p|^4)$.  Theorem 5.1 of \cite {NUE} gives a global gradient estimate under somewhat weaker conditions on $a$, specifically, $|a|= O(|p|^{4-\eta})$ with $\eta>0$ arbitrary, but only under stronger conditions on $b$:
\[
\delta_1b, \ \delta_2b \le O(b_p\cdot\gamma).
\]
In particular, if
\[
b(x,z,p)=v^{q-1}p\cdot\gamma + \psi(x,z,p),
\]
the condition on $\delta_1b$ requires $q=0$ or $q\ge1$. Our results do not fully extend those in \cite {NUE} (see also Theorems 9.8 and 9.9 from \cite {ODPbook}, which are essentially the estimates from \cite {NUE})  because some of the conditions in that source involve the operator $\delta_T$, defined by
\[
\delta_Tf= f_z+\frac 1{|p'|^2}p'\cdot f_x -\frac 1{2|p'|^2} f_p\cdot D(c^{km})p_kp_m;
\]
however, the author is unaware of any special equations which satisfy the conditions in \cite {NUE} are satisfied but not those in this work.

 \subsection {A generalization of the false mean curvature equation} \label {SFMCE}
Our results actually apply to a larger class of differential equations based on the one in \cite {NUE}.  These operators have the form
\begin {equation} \label {EFMCE}
a^{ij}=a^{ij}_*+\tau(x,z,p)p_ip_j 
\end{equation}
assuming that $\tau\ge0$ and $[a^{ij}_*]$ is a symmetric, uniformly elliptic matrix, that is, there is a positive function $\Lambda_*$ and there is a positive constant $\mu_*$ such that
\begin {equation} \label {Estar}
\mu_*\Lambda_*(x,z,p)|\xi|^2\le a^{ij}_*(x,z,p)\xi_i\xi_j\le \Lambda_*(x,z,p)|\xi|^2
\end {equation}
for all $\xi\in\mathbb R^n$ and all $(x,z,p)\in\Omega\times\mathbb R\times\mathbb R^n$ with $|p|\ge M_0$.  We also assume that
\begin {subequations} \label {Estarderiv}
\begin {gather}
\left|\frac {\partial a^{ij}_*}{\partial p}\right| =O((\tau\Lambda_*)^{1/2}), \label {Estarderiva} \\
|p||\delta_1a^{ij}_*| +
|p|\left|\frac {\partial a^{ij}_*}{\partial z}\right| +\left|\frac {\partial a^{ij}_*}{\partial x}\right| =o(|p|^2(\tau\Lambda_*)^{1/2}). \label {Estarderivb}
\end {gather}
\end {subequations}
The function $\tau$ is assumed to satisfy
\begin {subequations} \label {Etau}
\begin {gather}
|\tau_p|=O(\tau/|p|), \label {Etaua} \\
|\tau_x|+|p||\tau_z| =o(\tau|p|).
\end {gather}
\end {subequations}
These conditions on $\tau$ are essentially identical to those in \cite {NUE}. (See conditions (3.7) and (4.4) there.)  The only difference between the two  sets of conditions is that \cite {NUE} assumes an estimate on $\delta_T\tau$ while we rewrite the condition in terms of $\tau_z$ and $\tau_x$ separately.
The function $\tau$ is further assumed to  be appropriately related to $\Lambda_*$.  It was assumed in (4.4) of \cite {NUE} that $\tau=O(\Lambda_*)$.  Here, we shall assume instead
\begin {equation} \label {Elambdastartau}
\Lambda_*=o(\tau|p|^2)
\end {equation}
since, if $\tau=O(\Lambda_*/|p|^2)$, the equation is uniformly elliptic and stronger results are available.  We refer the reader to \cite {LT} and to the next example for details of the stronger results for uniformly elliptic equations. 

We are now ready to check our conditions.  First, if we take $r=-(\delta_1+3)\tau/\tau$, then
\[
(\delta_1+r+1)a^{ij}= (\delta_1+r+1)a^{ij}_*,
\]
so 
\[
(\delta_1+r+1)a^{ij}\eta_{ij} \le o(v)\tau^{1/2}\left(\Lambda_*\sum_{i,j=1}^n(\eta_{ij}^2)\right)\le o(\tau|p|^2)^{1/2}\left(a^{ij}\delta^{km}\eta_{ik}\eta_{jm}\right)^{1/2}.
\]
Since $\tau|p|^2= O(\mathscr E/v^2)$, we infer \eqref {SC1aij1} for some $\tilde\mu_\varepsilon$ determined by the limit behavior in \eqref {Estarderivb} and \eqref {Elambdastartau}. A similar calculation verifies \eqref {SC1aij2} with $s=0$ and  $\tilde\mu_\varepsilon(\sigma) = \bar\mu(\sigma)/\varepsilon^2$ for some decreasing function $\bar\mu$ satisfying
\[
\lim_{\sigma\to\infty}\bar\mu(\sigma)=0
\]
and determined by the limit behavior in \eqref {Estarderivb}, while \eqref {SC1aijp} follows from \eqref {Estarderiva} and \eqref {Etaua} with $\tilde\mu_*(s) = K_1/s$ for some positive constant $K_1$ determined by the limit behavior in those conditions.
If we also assume that
\begin {subequations} \label {Ea2}
\begin {gather}
|a|+|p||a_p| = O(\tau |p|^4), \label {Ea2a} \\
|a_x| =o(\tau |p|^5), \label {Ea2b}\\
a_z\le o(\tau|p|^4), \label {Ea2c}
\end {gather}
\end {subequations}
then \eqref {SC1p} holds with $\tilde\mu_*$ as before and $K_1$ determined also by the limit behavior in \eqref {Ea2a}.  Taking \eqref {Ea2b} and \eqref {Ea2c} into account, we also compute $B'_\infty=C_\infty=0$.

To obtain an interior gradient estimate, we use Theorem 3 from \cite {Serrin}.  First, we define the matrix $\mathscr A'=[a^{ij}_*]$ and take the multipliers $r=0$ and $s=-4$ in the notation of \cite {Serrin}.  It then follows that the quantity $a$ from (21) of \cite {Serrin} is zero, and that the quantities $b$ and $c$ from that equation are bounded from above by nonnegative constants determined only by the limit behavior in \eqref {Estarderiv}, \eqref {Etau}, and \eqref {Ea2}, so the hypotheses of Theorem 2 from \cite {Serrin} are satisfied.  In addition, further computation gives (29) from \cite {Serrin} with $\theta=1$ there.  Theorem 3 from \cite {Serrin} then gives a local gradient estimate for such equations.  
Hence, we obtain a gradient bound here, too, under the same conditions on $b$ as in the previous example.

We point out here that the estimates in \cite {NUE} require stronger hypotheses on $a^{ij}_*$ and on $a$ than the ones here.
For example, the second inequality of (3.1) in \cite {NUE} states that
\[
\left|\frac {\partial a^{ij}_*}{\partial p}\right| =O\left(\frac {\Lambda_*}{|p|}\right)
\]
and \eqref {Elambdastartau} implies that
\[
\frac {\Lambda_*}{|p|}= o((\tau\Lambda_*)^{1/2}),
\]
which means that the second inequality of (3.1) in \cite {NUE} is stronger than \eqref {Estarderiva}. 
Similar comparisons show that $|\partial a^{ij}_*/\partial z|$ and $|\partial a^{ij}_*/\partial x|$ can be larger here than in Theorem 4.2 of \cite {NUE}. 
In addition, Theorem 4.2 of \cite {NUE} requires $|a|=O(\tau|p|^3)$ while Theorem 5.1 of \cite {NUE} requires $|a|=O(\tau|p|^{4-\alpha})$ for some positive $\alpha$ or else $|a|=o(\tau|p|^4)$ and $\tau= O(\lambda_*|p|^{-\alpha})$ for some positive $\alpha$.  Finally, it is also assumed in \cite {NUE} that $\delta_1b \le O(b_p\cdot\gamma)$.

We also note that the choice $a^{ij}_*=  \exp(-v^2)\delta^{ij}$ and $\tau= v^{\alpha}$ (for any real $\alpha$) satisfies \eqref {Estar}, \eqref {Estarderiv}, \eqref {Etau}, and \eqref {Elambdastartau} but not the condition $\tau=O(\Lambda_*)$ from \cite {NUE}.

\subsection {Uniformly elliptic equations} \label {SSUE}
We now suppose that there is a positive constant $\mu^*$ such that
\begin {equation} \label {SC1UE}
\mu^*\Lambda|\xi|^2 \le a^{ij}\xi_i\xi_j
\end {equation}
for all vectors $\xi$. We also suppose that
\begin {subequations} \label {SCUE}
\begin {gather}
|p|^2|a^{ij}_p| + |p||a^{ij}_z| + |a^{ij}_x|= O(\Lambda|p|), \\
|p||a|+|p|^2|a_p|+ |a_x| = O(\Lambda|p|^3), \\
a_z \le O(\Lambda|p|^2).
\end {gather}
\end {subequations}
Finally, we assume that $b$ satisfies conditions \eqref {E10.32}, \eqref {E10.33}, and \eqref {Ebeta11}.

Straightforward calculation shows that conditions  \eqref {SC2aij}, \eqref {SCaij262}, \eqref {SC1aijploc}, and \eqref {ElambdaE2} are satisfied with $\theta=1$, $r=-1$, $s=0$, and the constants $\mu_2$ and $\mu_3$ determined by $n$, $\mu^*$, and the constants in \eqref {SCUE}. In addition, for any $y\in \partial\Omega$ and any $R>0$, the hypotheses of Theorem \ref {TN} are satisfied in $\partial\Omega \cap B(y,R)$ with $\varepsilon_x\equiv \beta_3$ and $\beta_2=\beta_3$. Further, $C_{\infty,R}$ and $B'_{\infty,R}$ are seen to be bounded by constants determined by the same quantities as for $\mu_2$ and $\mu_3$.
 A H\"older continuity estimate for $u$ follows from Theorem 2.3 
of \cite {LT} (or  Lemma 8.4 of \cite {ODPbook}), so we obtain a gradient bound from part \eqref {TClocso} of Theorem \ref {TCloc} if $R$ is sufficiently small.  In conjunction with the usual interior gradient estimate for uniformly elliptic equations Theorem 15.5 from \cite {GT}, we infer a global gradient estimate for solutions of \eqref {EODP} under these hypotheses.
Theorem 3.3 of \cite {LT} gives the exact same estimate as here, but there are two reasons to note the proof given here:  we do not need the special form for the interior gradient estimate proved as Lemma 3.1 in \cite {LT}, and we can use the $N$ function to derive this gradient estimate.

In particular, we obtain a gradient bound for four types of functions $b$.
The first is that $b$ has the form \eqref {Ebh} and $h$ and $\psi$ satisfy \eqref {Ehpsi}, the second is that $b$ has the form \eqref {Ebhtilde} and 
 $h$ and $\psi$ satisfy \eqref {Ehprime} and \eqref {Elimpsi1}, the third is that $b$ has the form \eqref {Ebq} with $q$ a $C^1$ function satisfying the inequality $q(x)\ge-1$, and $\sup_{\mathscr Q\times\mathbb R} |\psi|<1$ at any $x\in \mathscr Q$, where $\mathscr Q$ is the set on which $q=-1$, and the fourth is that $b$ has the form \eqref {Ebetastar1} with $h$ satisfying \eqref {Ebetastarh}.

\section {Global gradient bounds for parabolic problems} \label {SpG}

A large part of our analysis of parabolic problems is essentially identical to that of elliptic problems.  For this reason, we emphasize the few significant different difference between the two cases and sketch the minor modifications needed.

First, for a function $f$ depending on the variables $(x,t)$, we write $Df=\partial f/\partial x$ (exactly as before) and $f_t=\partial f/\partial t$, and we say that a function $f$ is in $C^{3,*}(\Sigma)$ for some open subset $\Sigma$ of $\mathbb R^n$ if the derivatives $D^3f$, $f_t$ and $Df_t$ are continuous in $\Sigma$.

Then, we write $\mathcal P\Omega\in C^{3,*}$ if $S\Omega$ can be written as the intersection of the set $\{(x,t)\in\mathbb R^n: 0<t\le T\}$ for some positive constant $T$ with a level set of a $C^{3,*}$ function $f$ with $|Df|$ bounded away from zero on $S\Omega$.  If the function $f$ is only $C^{2,1}$ (that is, if $D^2f$ and $f_t$ are continuous), then we write $\mathcal P\Omega\in C^{2,1}$.

We use $d^*$ to denote the parabolic distance to $S\Omega$, that is, we define $d^*$ by
\[
d^*(X)= \inf_{\substack {Y\in S\Omega\\s\le t}} \max\{|x-y|, |t-s|^{1/2}\},
\]
and we write $d$ for the spatial distance to $S\Omega$, that is,
\[
d(X)= \inf_{\substack {Y\in S\Omega\\s= t}} |x-y|,
\]
where $Y= (y,s)$. As pointed out in Section 10.3 of \cite {LBook}, there is a positive constant $C$, determined only by $\Omega$, such that $d^*/d$ is trapped between $C$ and $1/C$ in $\Omega$ whenever $\mathcal P\Omega\in C^{2,1}$.  In that reference, it is also shown that, if $\mathcal P\Omega \in C^{3,*}$, then there is a positive constant $R_0$, determined only by $\Omega$, such that $d\in C^{3,*}(\Sigma_{R_0})$, where $\Sigma_{R_0}$ is the subset of $\clOmega$ on which $d<R_0$.

Further, for $r>0$, we write $\Omega_r$ for the subset of $\Omega$ on which $d<r$.  If $\mathcal P\Omega\in C^{3,*}$, then $\gamma=Dd$ is a $C^{2,1}$ unit vector field in $\Sigma_{R_0}$ which extends the unit inner spatial normal to $S\Omega$.  Moreover, for any vector $\xi$, the vector $\tilde\xi$, defined by \eqref {E1.1}, satisfies $|\tilde\xi|\le 2|\xi|/R_0$ in $\Omega_{R_0/2}$.

Since our problem depends on $t$, we need to modify notation slightly before presenting our $N$ function.  For $\tau\ge1$, we write $\Sigma(\tau)$ for the subset of $S\Omega\times\mathbb R\times\mathbb R^n$ on which $v\ge \tau$.  In place of Theorem \ref {TN}, we have the following result.
Since the proof is essentially the same as for Theorem \ref {TN} (except for notational adjustments), we omit it here.  We do mention, though, that the notation $C^{2,1}(\overline{\Omega_{R_0/4}}\times \mathbb R\times\mathbb R^n\times(0,\varepsilon))$ means a function which is once continuously differentiable with respect to  $t$ and twice continuously differentiable with respect to all other variables.

\begin {theorem} \label {TNp}
Let $\mathcal P\Omega\in C^{2,1}$, let $\tau_0\ge1$, and let $b\in C^1(\Sigma(\tau_0)$ with $b_p\cdot\gamma>0$ on $\Sigma_0(\tau_0)$, the subset of $\Sigma(\tau_0)$ on which $b=0$.  Suppose
\begin {equation} \label {E3.1'}
\lim_{t\to\infty} b(X,z,p-\tau \gamma(X)) <0< \lim_{t\to\infty} b(X,z,p+\tau \gamma(X)) 
\end {equation}
for all $(X,z,p)\in \Sigma(\tau_0)$ and that there are positive constants $\beta_0$ and $c_0$ such that conditions \eqref {E10.33} are satisfied on $\Sigma_0(\tau_0)$.  Suppose also that there are a $*$-decreasing function $\varepsilon_x$ and a nonnegative constant $\beta_1$such that \eqref {E10.35} holds.
 Then there is a positive constant $\varepsilon_0(\beta_0,c_0)$, along with a $C^{2,1}(\overline{\Omega_{R_0/4}}\times \mathbb R\times\mathbb R^n\times (0,\varepsilon_0))$ function $N$, such that conditions \eqref {E10.37} hold on $\Omega_{R_0/4}\times\mathbb R\times\mathbb R^n\times(0,\varepsilon_0)$ with $v_\varepsilon$ as in Theorem \ref {TN}.  Moreover, $N=0$ on $\Sigma_0((1+c_0^2)^{1/2}\tau_0)\times\mathbb R\times\mathbb R^n\times(0,\varepsilon_0)$.

If we define $w$ and $\nu_1$ by \eqref {E10.38} for some $\varepsilon\in(0,\varepsilon)$, then \eqref {E10.390} is valid.  Also, there is a nonnegative constant $c_1$, determined only by $\beta_0$, $c_0$, and $n$ such that \eqref {E10.372pp}, \eqref {E10.43}, \eqref {E10.372pz}, and \eqref {E10.372px} all hold for all $\xi\in\mathbb R^n$, and \eqref {E10.372zz} and \eqref {E10.372xz} are satisfied.  Finally, if $\mathcal P\Omega\in C^{3,*}$ and if there is an increasing function $\Lambda_0$ such that
\begin {equation} \label {Ebetat}
|b_t| \le \Lambda_0(v)b_p\cdot\gamma
\end {equation}
on $\Sigma_0(\tau_0)$, then there is a nonnegative constant $c_2$, determined also by $\Omega$, such that \eqref {E10.372xx} is valid and
\begin {equation}
|N_t| \le c_2v+2\Lambda_0(2(1+c_0)v).
\end {equation}
\end {theorem}

We also observe that, if $\Omega$ is a cylinder and if $b$ is time-independent, then the proof of the preceding theorem shows that $N$ is also time-independent.  In particular, all of our elliptic gradient estimates have direct parabolic analogs in this case.  To facilitate our exposition, we include these analogs in the more general results given below.

The first way in which our parabolic gradient estimates differ from the elliptic ones is that the multiplier functions $r$ and $s$ must be chosen as $r=-1$ and $s=0$.  For the reader's convenience, we rewrite the appropriate conditions explicitly for this pair of multipliers.  First, we set
\begin {align*}
A &= \frac {(\delta_1-1)\mathscr E}{\mathscr E}, \\
B &= \frac {(\delta_2E + (\delta_1-1)a}{\mathscr E}, \\
B'  &= \frac {(\delta_1-1)a}{\mathscr E}, \\
C &= \frac {\delta_2a}{\mathscr E},
\end {align*}
and we  introduce the following constants related to the limit behavior these functions.
\begin {subequations} \label {Elimsupp}
\begin {align}
A'_\infty &= \limsup_{|p|\to\infty} \sup_{(X,z)\in \Omega\times\mathbb R} |
A(X,z,p)|, \\
B_\infty &= \limsup_{|p|\to\infty} \sup_{(X,z)\in \Omega\times\mathbb R} B(X,z,p), \\
B'_\infty &= \limsup_{|p|\to\infty} \sup_{(X,z)\in \Omega\times\mathbb R} |B'(X,z,p)|, \label {Elimsuppb'}\\
C_\infty &= \limsup_{|p|\to\infty} \sup_{(X,z)\in \Omega\times\mathbb R} C(X,z,p). \label {Elimsuppc}
\end {align}
\end {subequations}
For ease of notation, we modify the definition of $\Gamma(M_0)$ to deonte the set of all $(X,z,p)\in \Omega\times\mathbb R\times\mathbb R^n$ with $|p|\ge M_0$.

We begin with estimates that are simple analogs of the corresponding elliptic ones.

\begin {theorem} \label {TCp}
Let $u\in C^{2,1}(\clOmega)$ be a solution of \eqref {EPODP} with $\mathcal P\Omega\in C^{3,*}$ and $b$ satisfying the hypotheses of Theorem \ref {TNp} for some $*$-decreasing function $\varepsilon_x$ and some increasing function $\Lambda_0$ such that
\begin {equation} \label {ELambda0}
\Lambda_0((1+c_0)v) \le \mu_4(1+\mathscr E)v
\end {equation}
on $\Gamma(M_0)$ for some nonnegative constant $\mu_4$.
Suppose that there is a $*$-decreasing function $\tilde\mu_*$ such that \eqref {SC1aijp} is satisfied on $\Gamma(M_0)$ for all matrices $[\eta_{ij}]$ and \eqref {SC1p} holds on $\Gamma(M_0)$.
\begin {enumerate}
\item \label {ITCp}
Suppose also that there are decreasing functions $\tilde\mu$ and $\tilde\mu_\varepsilon$ satisfying \eqref {Etildelimit} and a nonnegative constant $\mu_2$ such that
\begin {subequations} \label {pSC1aij}
\begin {gather}
\delta_1a^{ij}\eta_{ij}\le \frac {\tilde \mu(v)}{v} \mathscr E^{1/2}(a^{ij}\delta^{km}\eta_{ij}\eta_{jm})^{1/2}, \label {pSC1aij1}\\
\delta_2a^{ij}\eta_{ij} \le \frac {\tilde \mu_\varepsilon(v)}{v} \mathscr E^{1/2}(a^{ij}\delta^{km}\eta_{ij}\eta_{jm})^{1/2}, \label {pSC1aij2}
\end {gather}
\end {subequations}
hold on $\Gamma(M_0)$ for all matrices $[\eta_{ij}]$, that 
\begin {equation} \label {ELambdav2}
\Lambda|p|^2 \le\mu_2\mathscr E
\end {equation}
on $\Gamma(M_0)$, that  $B'_\infty$ is bounded uniformly with respect to $\varepsilon$ and that $C_\infty\le0$.   If  \eqref {Edelta1b} is satisfied on $\Sigma_0(\tau_0)$ and if $\varepsilon_x$ satisfies \eqref {Etildelimit1}, then there is a constant $M$, determined only by $\beta_0$, $\beta_1$, $c_0$, $n$, $\mu_4$, $\sup_{\{d\ge R_0/4\}} |Du|$, the oscillation of $u$, and the limit behavior in \eqref {Etildelimit}, \eqref {Etildelimit1}, and \eqref {Elimsupp}, such that $|Du| \le M$ in $\Omega$.
\item \label {ITCp1}
Suppose also that there are a decreasing function $\tilde\mu_\varepsilon$ satisfying \eqref {Etildelimiteps} and a nonnegative constant $\mu_2$ such that \eqref {SC1aijp},  \eqref {pSC1aij2}, and
\begin {equation} \label {pSCaij}
\delta_1 a^{ij}\eta_{ij} \le \frac { \mu_2}{v} \mathscr E^{1/2}(a^{ij}\delta^{km}\eta_{ij}\eta_{jm})^{1/2}
\end {equation}
hold on $\Gamma(M_0)$ for all matrices $[\eta_{ij}]$ and \eqref {ELambdav2} holds on $\Gamma(M_0)$, that $B'_\infty$ is finite, and that $C_\infty\le 0$. If $\varepsilon_x$ satisfies \eqref {Etildelimit1} and if $b$ satisfies \eqref {Elimbz},  then there is a positive constant $M$, determined only by $\beta_0$, $\beta_1$, $c_0$, $n$,  $\mu_2$, $\mu_4$, the oscillation of $u$, $\sup_{\{d\ge R_0/4\}} |Du|$,
and the limit behavior in \eqref {Etildelimit}, \eqref {Etildelimit1}, \eqref {Elimbz}, and \eqref {Elimsupp}, such that $|Du| \le M$ in $\Omega$. 
\item \label {ICCp}
Suppose that there are a constant $\mu_2$ and a decreasing function $\tilde\mu_\varepsilon$ satisfying \eqref {Etildelimiteps} such that \eqref {SC1aijp}, \eqref {pSC1aij2}, and  \eqref {pSCaij} hold on $\Gamma(M_0)$ for all matrices $[\eta_{ij}]$ and \eqref {ElambdaE1} holds on $\Gamma(M_0)$, and suppose $B'_\infty$ is finite and $C_\infty<0$.   If \eqref {Edelta1b} and \eqref {Edelta1beta} are satisfied on $\Sigma_0(\tau_0)$ and if $\varepsilon_x$ satisfies \eqref {Etildelimit1}, then there is a constant $M$, determined only by $\beta_0$, $\beta_1$, $c_0$, $n$,  $\mu_4$, $\sup_{\{d\ge R_0/4\}} |Du|$, the oscillation of $u$, and the limit behavior in \eqref {Etildelimit}, \eqref {Etildelimit1}, and \eqref {Elimsupp}, such that $|Du| \le M$ in $\Omega$.
\item \label {ITpGBsmall}
 Suppose also that there are nonnegative constants $M_0$ and $\mu_2$ such that conditions \eqref {pSCaij} and
\begin {equation} \label {SCpaij26}
|\delta_2a^{ij}\eta_{ij}| \le \frac {\mu_2}v \mathscr E^{1/2}(a^{ij}\delta^{km}\eta_{ik}\eta_{jm})^{1/2}
\end {equation}
hold on $\Gamma(M_0)$ for all matrices $[\eta_{ij}]$, that \eqref  {ELambdav2} hold on $\Gamma(M_0)$, and suppose that $B'_\infty$ and $C_\infty$ are finite for each $\varepsilon$. 
If  $\varepsilon_x$ satisfies \eqref {Evarep}, if \eqref {Edelta1beta1} holds on $\Sigma_0(\tau_0)$, and if \eqref {ELambda0} holds, then there are constants $M$ and $\omega_0$, determined only by $\beta_0$, $\beta_1$, $c_0$, $\mu_2$, $\mu_4$, $\Omega$,  $\sup_{\{d\ge R_0/4\}}|Du|$, and  the limit behavior in \eqref {Elimsupp} and \eqref {Evarep} such that $|Du|\le M$ in $\Omega$ provided $\osc_{\Omega} u \le \omega_0$.
\end {enumerate}
\end {theorem}
\begin {proof}
To prove part \eqref {ITCp}, we use the notation from the proof of Theorem \ref {TCinfty}. 
Now, we use a different estimate for $\omega I_0a^{ij}D_{ij}u$ from the one in the elliptic case.  First we write
\[
\omega I_0a^{ij}D_{ij}u = \omega I_0 \psi'a^{ij}D_{ij}\bar u + \omega^2 I_0 \mathscr E.
\]
 It follows from \eqref {ES}, \eqref {Ei0est1}, \eqref {ELambdav2}, and the Cauchy-Schwarz inequality that
\[
\omega I_0\psi'a^{ij}D_{ij}\bar u \ge -c (\mu_2\varepsilon)^{1/2}(w_1\mathscr E)^{1/2}\mathscr S^{1/2},
\]
while \eqref {EI0} yields
\[
\omega^2 I_0\mathscr E \ge -c\varepsilon^{1/2}\omega^2w_1\mathscr E.
\]
 
Using these estimates and assuming initially that $Du\in C^{2,1}(\Omega)$, we see from the proof of Theorem \ref {TCinfty} that, by choosing $\psi$ suitably, there are positive constants $\tilde\eta$ and $M_1$ such that $w_2$ satisfies the differential inequality
\[
-\frac {\partial w_2}{\partial t} +a^{ij}D_{ij}w_2 +\kappa^iD_iw_2\ge \left(\tilde\eta w_1\mathscr E -  \frac {2\varepsilon NN_t}{(\psi')^2}+2c^{km}_t D_k\bar u D_m\bar u\right) \left( 1+ \frac {k_1d}{\tilde\mu_*(\sqrt{w_1})}\right)
\]
on $\Omega'$, the subset of all $\Omega_{R_0/4}$ on which $|Du| \ge M_1$. As in Theorem \ref {TCinfty}, it follows that
\[
w_2\ge \frac {2w_1}{3\tilde\mu_*(\sqrt{w_1})},
\]
and hence (after taking \eqref {ELambda0} into account) there are nonnegative constants $k_2$ and $k_3$, determined only by $c_0$, $\mu_4$, $\sup|\gamma_t|$, and $n$, such that
\[
-\frac {\partial w_2}{\partial t} +a^{ij}D_{ij}w_2 +\kappa^iD_iw_2\ge [\tilde \eta -k_2 \varepsilon] w_2\mathscr E-k_3w_2
\]
on $\Omega'$.

We now decrease $\varepsilon$ so that $\tilde\eta\ge k_3\varepsilon$ to infer that
\[
-\frac {\partial w_2}{\partial t}+a^{ij}D_{ij}w_2+\kappa^iD_iw_2+k_3w_2 \ge 0
\]
on $\Omega'$,  and hence $w_3=e^{-k_3t}w_2$ cannot take its maximum in $\Omega'$. The restriction that $Du\in C^{2,1}$ is removed by using Theorem 6.15 from \cite {LBook} rather than Theorem 8.1 from \cite {GT}.

The analysis of a boundary maximum for $w_3$ is only slightly different from the analysis in the elliptic case.  This time, we note that there is a positive constant $M_2$ such that
\[
b^iD_iw_3>0
\]
at any point of $S\Omega$ where $w_3\ge M_2$, while the size of $w_3$ at  $\mathcal P\Omega_{R_0/4} \setminus S\Omega$ is controlled by the assumed bounds on $|Du|$ there.

The other parts are proved via a similar modification of the corresponding elliptic theorems.
\end {proof}

There is another way to obtain a global gradient bound for parabolic problems using as our model Corollary 11.2 from \cite {LBook}. Specifically, we have the following result, which includes many cases of interest.

\begin {theorem} \label {Tpspecial}
Let $u\in C^{2,1}(\clOmega)$ be a solution of \eqref{EPODP} with $\mathcal P\Omega\in C^{3,*}$ and $b$ satisfying the hypotheses of Theorem \ref {TNp} for some $*$-decreasing function $\varepsilon_x$ and  $\Lambda_0(\sigma)= \mu_5\sigma$ for some nonnegative constant $\mu_5$.
Suppose that there is a nonnegative function $\mu_6$ such that 
\begin {subequations} \label {SCpspecialij}
\begin {gather}
\delta_1a^{ij}\eta_{ij} \le \frac {\mu_5}v \left( a^{ij}\delta^{km}\eta_{ik}\eta_{jm}\right)^{1/2}, \label {SCpspecialaij1}\\
\delta_2a^{ij}\eta_{ij} \le \frac {\mu_6(\varepsilon)}v \left( a^{ij}\delta^{km}\eta_{ik}\eta_{jm}\right)^{1/2}, \label {SCpspecialaij2} \\
|a^{ij}_p\eta_{ij}| \le \frac {\mu_5}v \left( a^{ij}\delta^{km}\eta_{ik}\eta_{jm}\right)^{1/2} \label {SCpspecialaijp}
\end {gather}
\end {subequations}
on $\Gamma(M_0)$ for all matrices $[\eta_{ij}]$ and
\begin {subequations} \label {SCpspecial}
\begin {gather}
|(\delta_1-1)a| \le \mu_5, \label {SCpspeciala1} \\
\delta_2a \le \mu_6(\varepsilon), \label {SCpspeciala2} \\
\mathscr E \le \mu_5 \label {SCpspeicalE}
\end {gather}
\end {subequations}
on $\Gamma(M_0)$.
If
\begin {equation} \label {Elambdaspecial}
(1+\varepsilon_x(v)v)^2\Lambda \le \mu_5
\end {equation}
on $\Gamma(M_0)$ and if \eqref {Edelta1beta} holds on $\Sigma_0(\tau_0)$, then there is a constant $M$, determined only by $n$, $\beta_0$, $\beta_1$, $c_0$, $\mu_5$, $\mu_6$, $M_0$, $\tau_0$, $\Omega$, and $\sup_{\{d\ge R_0/4\}} |Du|$ such that $|Du| \le M$ in $\Omega$.
\end {theorem}
\begin {proof} We first fix $\varepsilon \in (0,\varepsilon_0)$.  With $w_1$ as in the proof of Theorem \ref {TCinfty} and $k_1$ a positive constant to be determined, we set
\[
w_2= w_1+ k_1d\int_0^{w_1} (1+ \varepsilon_x(\sqrt \sigma)\sqrt \sigma)\, d\sigma.
\]
Following the proof of Theorem \ref {TCinfty} and writing $\Omega'$ for the subset of all $X\in \Omega_{R_0/4}$ with $|Du|\ge \max\{M_0,\tau_0\}$, we see that 
\[
b^iD_iw_2 >0
\]
on $S\Omega\cap S\Omega'$ provided $k_1$ is chosen sufficiently large. In addition, there is a constant $c$ such that
\[
-\frac {\partial w_2}{\partial t} +a^{ij}D_{ij}w_2 +\kappa^iD_iw_2 \ge -cw_2
\]
on $\Omega'$.  Therefore $w= \exp(-ct)w_2$ satisfies
\begin {gather*}
-w_t+a^{ij}D_{ij}w +\kappa^iD_iw>0 \text { on } \Omega', \\
b^iD_iw>0 \text { on } S\Omega\cap S\Omega'.
\end {gather*}
These differential inequalities then give the desired bound for $w$, and hence for $|Du|$, by the maximum principle.
\end {proof}

\section {Local gradient estimates for parabolic problems} \label {SLp}

Local estimates are generally straightforward If  $C\le0$ or if the oscillation of $u$ is sufficiently small; however, there is an additional complication if $\Omega$ is not cylindrical.  To keep our discussion of these estimates consistent in the various cases, we start by rewriting the notation from Section \ref {Sso}.

First, for $Y\in \mathbb R^{n+1}$ and positive constants $R$ and $R'$, we write
\[
Q(Y;R,R')= \{X\in \mathbb R^{n+1}: |x-y|<R, s-R'<t<s\}.
\]
Next, for positive constants $M_0$, $R$, $R'$, and $\tau$, we write $\Gamma_{R,R'}(M_0)$ for the set of all $(X,z,p)\in \Gamma(M_0)$ with $X\in Q(Y;R,R')$ and $\Sigma_0(\tau;R,R')$ for the set of all $(X,z,p)\in \Sigma_0(\tau)$ with $X\in Q(Y;R,R')$. 
It will also be convenient to define 
\[
\Omega'(Y;R,R')= \mathscr P(\Omega\cap Q(Y;R,R')) \setminus S\Omega.
\]

Our local problem is then written in the form
\begin {subequations} \label {EPODPloc}
\begin {gather}
-u_t+a^{ij}(X,u,Du)D_{ij}u +a(X,u,Du)=0 \text { in }\Omega\cap Q(Y;R,R'), \\
b(X,u,Du)=0 \text { in } S\Omega \cap Q(Y;R,R'),
\end {gather}
\end {subequations}
and we introduce the numbers 
\begin {subequations} \label {ERinftyp}
\begin {align}
B'_{\infty;R,R'} &= \limsup_{|p|\to\infty} \sup_{(X,z)\in\Omega\cap Q(Y;R,R')\times\mathbb R} |B'(X,z,p)|, \label {ERinftb'p} \\
C_{\infty;R,R'} &= \limsup_{|p|\to\infty} \sup_{(X,z)\in\Omega\cap Q(Y;R,R')\times\mathbb R} C(X,z,p). \label {ERinftycp}
\end {align}
\end {subequations}

Our first local estimate has the following form.

\begin {theorem} \label {TClocp} 
Let $Y\in S\Omega$, let $R$ and $R'$ be positive numbers and suppose $\mathcal P\Omega \cap Q(Y;R,R')\in C^{3,*}$ with $d<R_0/4$ in $\Omega \cap Q(Y;R,R')$
Suppose that the hypotheses of Theorem \ref {TCp} are all satisfied with $\Omega\cap Q(Y;R,R')$ in place of $\Omega$, $S\Omega\cap Q(Y;R,R')$, $\tilde\mu_*(\sigma)= \mu_3\sigma^\theta$ for some positive constants $\mu_3$ and $\theta$ with $\theta<1$, and $\tilde\mu=\varepsilon_x$.
Then $|Du(Y)| \le M$ for some  number $M$ determined by the same quantities as in Theorem \ref {TCp} with $\Omega'(Y;R,R')$ in place of $\{d\ge R_0/4\}$ and also by $\mu_3$ and $\theta$.
\end {theorem}
\begin {proof} We note that the function $\zeta$ from the proof of Theorem \ref {TCloc} (now considered as a function of $X$ which is constant with respect to $t$) satisfies $b_p\cdot D\zeta\ge0$ at any point of $S\Omega\cap Q(Y;R,R_0)$ with $\zeta\ge0$.  The proof of Theorem \ref {TCloc} (with the usual parabolic modifications, as in the proof of Theorem \ref {TCp}) then gives the desired estimate.
\end {proof}

Note that the assumption that $d<R_0/4$ in $\Omega\cap Q(Y;R,R')$ is a restriction on the number $R'$ in a possibly convoluted way.  In particular, if $\Omega$ is cylindrical, which means that $\Omega = S\times (0,T)$ for some fixed subset $S$ of $\mathbb R^n$, then this assumption is automatically satisfied as long as $R<R_0$.

The local estimate corresponding to Theorem \ref {Tpspecial} is somewhat more complicated in that we need to replace the function $\zeta$ by a slightly different function.  Our model is the function $\varphi$ from the proof of Theorem 2.1 in \cite {EH}. This function is defined as
\[
\varphi(X) = 1- \frac {|x-x_0|^2}{R^2} -Kt
\]
for a fixed point $x_0$ and a fixed number $R$ while $K$ is a positive constant at our disposal.

Using this function, we first prove an interior gradient estimate which will be useful in discussing our local gradient bound, which is an analog of Theorem \ref {Tpspecial}. To simplify its statement, we define $\Omega^*$ to be the set of all $(X,z,p)\in B(x_0,R)\times (0,T)\times\mathbb R\times \mathbb R^n$ with $|p|\ge M_0$ with $x_0\in \mathbb R^n$, $R$, $T$, and $M_0$ $R$ given.  We also define the operator $\delta_3$ by
\[
\delta_3f(x,z,p) = f_z(x,z,p)+ \frac {f_k(x,z,p)\delta^{km}p_m}{|p|}.
\]

\begin {theorem} \label {Tpintgrad} Let $u\in C^{2,1}(B(x_0,R)\times (0,T))$ be a solution of
\begin {equation} \label {ETpintgrad}
-u_t +a^{ij}(X,u,Du)D_{ij}u+a(X,u,Du) =0 \text { in } B(x_0,R)\times (0,T)
\end {equation}
for some $x\in \mathbb R^n$ and positive constants $R$ and $T$. Suppose  that there are positive constants $\mu_5$, $\mu_6$, $\mu_7$, $M_0$, and $q$ such that
\begin {subequations} \label {SCTpintgradaij}
\begin {gather}
q>1,\ \mu_7^2<8\left(1- \frac 1q\right), \label {SCTpintgradaijq} \\
\delta_3a^{ij}\eta_{ij} \le \mu_5 \frac {(a^{ij}\delta^{km}\eta_{ik}\eta_{jm})^{1/2}}{|p|}, \label {SCTpintgradaij1}\\
\begin {split}
(|p|^2a^{ij,k} +4a^{ik}\delta^{jm} p_m)\eta_{ij} \xi_k\le \left(\mu_6 |p|^{1-1/q} |\xi| +\mu_7|p|(a^{ij}\xi_i\xi_j)^{1/2}\right)  \\
\times \left(a^{ij}\delta^{km}\eta_{ik}\eta_{jm}\right)^{1/2} 
\end {split} \label {SCTpintgradaij2}
\end {gather}
\end {subequations} 
on $\Omega^*$ for all matrices $[\eta_{ij}]$ and vectors $\xi$. Suppose also that 
\begin {subequations} \label {SCTpintgrada}
\begin {gather}
\Lambda \le \mu_5, \label {SCLambdabound}\\
\delta_3a\le \mu_5, \\
|a_p| \le \mu_5
\end {gather}
\end {subequations}
on $\Omega^*$.
Then there  is a constant $T_0\in(0,T]$, determined only by $\mu_5$, $\mu_6$, $\mu_7$, $q$, and $R$, such that
\begin {equation} \label {ETpintgradest}
\sup_{B(x_0,R/2)\times (0,T_0)}|Du| \le 4\max\{M_0, \sup_{B(x_0,R)\times\{0\}}|Du|\}.
\end {equation}
\end {theorem}
\begin {proof}
First, we set $w=|Du|^2$,  $\mathscr C^2= a^{ij}\delta^{km}D_{ik}uD_{jm}u$, and $\kappa = a^{ij}_pD_{ij}u+a_p$, and we define the operator $\mathscr L$ by
\[
\mathscr Lh =-h_t+a^{ij}D_{ij}h +\kappa^iD_ih.
\] 
A simple calculation shows that
\begin {equation} \label {ELw}
\mathscr Lw=2\mathscr C^2 -2[\delta_3a^{ij}D_{ij}u+\delta_3a]w.
\end {equation}
With $K$ a positive constant to be determined, we now set
\[
\varphi(X)=1- \frac {|x-x_0|^2}{R^2}-Kt,
\]
and we write $\Gamma$ for the subset of $B(x_0,R)\times(0,T)$ on which $\varphi$ is positive. It then follows that $w_3= \varphi^qw$ satisfies the equation
\begin {align*}
\mathscr Lw_3 &= \varphi^q\mathscr Lw +qw\varphi^{q-1}(a^{ij}D_{ij}\varphi-\varphi_t +a^{ij,k}D_{ij}uD_k\varphi +a^kD_k\varphi))  \\&\phantom {=\varphi^q }+q(q-1)w\varphi^{q-2}a^{ij}D_i\varphi D_j\varphi \\
&\phantom {=\varphi^q }+4q\varphi^{q-1}a^{ij}\delta^{km}D_{ik}uD_mu D_j\varphi.
\end {align*}
(Note that we have rewritten $D_iw$ as $2\delta^{km}D_{ik}uD_mu$.) We now use \eqref {ELw} rewrite this equation as
\[
\mathscr Lw_3= 2\varphi^q\mathscr C^2 +J_1w\varphi^{q-1}  +q(q-1)\varphi^{q-2}wa^{ij}D_i\varphi D_j\varphi +J_2 +J_3
\]
with
\begin {align*}
J_1 &=q(a^{ij}D_{ij}\varphi -\varphi_t+a^kD_k\varphi)-2 \delta_3a \varphi, \\
J_2 &=q\varphi^{q-1}( |Du|^2a^{ij,k} + 4a^{ik}\delta^{jm}D_mu)D_{ik}uD_k\varphi, \\
J_3 &=-2\varphi^qw\delta_3a^{ij}D_{ij}u.
\end {align*}
We first use \eqref {SCTpintgradaij2} to see that
\[
J_2\ge -\mu_6q\varphi^{q-1}|Du|^{1-1/q} |D\varphi|\mathscr C - \mu_7q\varphi^{q-1}|Du|\mathscr C(a^{ij}D_i\varphi D_j\varphi)^{1/2}
\]
and \eqref {SCTpintgradaij1} to see that
\[
J_3\ge -2\mu_5|Du|\mathscr C.
\]
By invoking the Cauchy-Schwarz inequality, we conclude, for any positive constants $\theta_1$ and $\theta_2$, that
\[
\mathscr Lw_3 \ge (2-2\theta_1-\theta_2)\mathscr C^2 +J_4w\varphi^{q-1} +J_5qw\varphi^{q-2} a^{ij}D_i\varphi D_j\varphi,
\]
with
\begin {align*}
J_4&= J_1- \frac {\mu_6^2q^2|D\varphi|^2}{4\theta_1w_3^{1/q}} - \frac {\mu_5^2}{\theta_3}\varphi,  \\
J_5 &=q-1- \frac {\mu_7^2q}{4\theta_2}.
\end {align*}
We first choose $\theta_2$ to make $J_5=0$, that is,
\[
\theta_2=\frac {q\mu_7^2}{4(q-1)},
\]
noting that the second inequality of \eqref {SCTpintgradaijq} implies that $\theta_2\in(0,2)$.  Then we choose $\theta_1= 1-(\theta_2)/2$ to infer that
\[
\mathscr Lw_3 \ge J_4w\varphi^{q-1}.
\]
From \eqref {SCTpintgrada} and the explicit form of $\varphi$, we conclude that there is a constant $C$, determined only by $n$, $R$, $\mu_5$ such that
$J_4\ge K-C$ on $\Omega^{**}$, which is the subset of $\Omega^*$ on which $w_3\ge M_0$.  We finally choose $K=C$ and use the maximum principle (Corollary 2.5 from \cite {LBook}) to conclude that $w_3$ must attain its maximum over $\overline {\Omega^{**}}$ somewhere on $\mathscr P\Omega^{**}$ and hence
\[
w_3 \le \max\{M_0^2, \sup_{B(x_0,R)} w_3\}.
\]
The proof is completed by choosing $T_0$ small enough that $\varphi\ge \frac 14$ on $B(x_0,R/2)\times (0,T_0)$.
\end {proof}

Condition \eqref {SCTpintgradaij2} seems rather artificial but our examples will show its value.  We also shall take advantage of a slight variant of this result, which is closer to Theorem 2.1 of \cite {EH}.
In particular, we shall show that it is valid (with suitable constants) if $a^{ij}=g^{ij}$, which is the case studied by Ecker and Huisken in \cite {EH}; in particular, our theorem is a simple generalization of Theorem 2.1 from that work.

\begin {theorem} \label {Tpspecialloc}
Let $u\in C^{2,1}(\clOmega \cap Q(Y;R,R'))$ be a solution of \eqref {EPODPloc} with $\mathcal P(\Omega\cap Q(Y;R,R'))\in C^{3,*}$ and $b$ satisfying the hypotheses of Theorem \ref {TNp}, with $\Omega\cap Q(Y;R,R')$ in place of $\Omega$ and $S\Omega \cap Q(Y;R,R')$ in place of $S\Omega$ for some $*$-increasing function $\varepsilon_x$ and $\Lambda_0(s)= \mu_5s$ for some nonnegative constant $\mu_5$.
Suppose that $d\le R_0$ in $\Omega\cap Q(Y;R,R')$ and that $\Lambda_0$ satisfies \eqref {ELambda0}  on $\Gamma_{R,R'}(M_0)$.
Suppose also that there are a nonnegative function $\mu_6$ and a constant $\theta\in(0,1)$ such that \eqref {SCpspecialaij1}, \eqref {SCpspecialaij2}, and
\begin {equation} \label {SCpspecialaijp1}
|a^{ij}_p\eta_{ij}| \le \frac {\mu_5}{v^{1+\theta/2}} \left( a^{ij}\delta^{km}\eta_{ik}\eta_{jm}\right)^{1/2} 
\end {equation}
hold on $\Gamma_{R,R'}(M_0)$ for all matrices $[\eta_{ij}]$ and \eqref {SCpspecial} hold on $\Gamma_{R,R'}(M_0)$.
If also \eqref {Elambdaspecial} and
\begin {equation} \label {Elambdaspecial2}
(1+\varepsilon_x(v)v)v^\theta\Lambda \le \mu_5
\end {equation}
hold on $\Gamma_{R,R'}(M_0)$ and if \eqref {Edelta1beta} holds on $\Sigma_0(\tau_0)$, then there are a constant $R_2$, determined only by $n$, $\beta_0$, $\beta_1$, $c_0$, $\mu_5$, $\mu_6$, $M_0$, $\tau_0$, $\Omega$, and a constant $M$ determined also by $\sup_{\Omega'(Y;R,R')} |Du|$ such that $|Du(Y)| \le M$ in $\Omega$ provided $R'\le R_2$.
\end {theorem}
\begin {proof}
With $w_2$ as in Theorem \ref {Tpspecial} and $\mathscr L$ as in Theorem \ref {Tpintgrad}, we infer that there are positive constants $\tilde\eta$ and $k_2$ such that
\[
\mathscr Lw_2 \ge (1+k_1d[1+\sqrt{w_1}\varepsilon_x(\sqrt{w_1})])[ \tilde \eta w_2\mathscr E+\mathscr S] -k_2w_2.
\]
Next, we set $\varphi_1(X)= \zeta-Kt$ with $K$ a constant at our disposal and $w=\zeta^qw_2$ with $q=2/\theta$.  A straightforward calculation, taking \eqref {ED1w} into account, yields 
\[
\mathscr Lw = \eta^q\mathscr Lw_2 + C_6 +S_8
\]
with
\begin {align*}
C_6 &= qw_2\varphi_1^{q-1}(a^{ij}D_{ij}\zeta +a^kD_k\zeta +K) \\
 &\phantom{=}{} +q(q-1)\varphi_1^{q-2}w_2a^{ij}D_i\zeta D_j\zeta \\
 &\phantom{=}{}+4( 1+ k_1d[1+\varepsilon_x(\sqrt{w_1})\sqrt{w_1}\,])q\varphi_1^{q-1}NN_za^{ij}D_i\zeta D_ju \\
 &\phantom {=} {}+2( 1+ k_1d[1+\varepsilon_x(\sqrt{w_1})\sqrt{w_1}\,])q\varphi_1^{q-1}Na^{ij}D_i\zeta N_j \\
 &\phantom {=}  {}+2( 1+ k_1d[1+\varepsilon_x(\sqrt{w_1})\sqrt{w_1}\,])q\varphi_1^{q-1}a^{ij}D_i\zeta D_j(c^{km})D_kuD_mu \end {align*}
 and
 \begin {align*}
 S_8 &= 4( 1+ k_1d[1+\varepsilon_x(\sqrt{w_1})\sqrt{w_1}\,]) q\varphi_1^{q-1}a^{ij}D_i\zeta D_{jk}u\nu_1^k \\
 &\phantom {=}{}+ q\varphi_1^{q-1} w_2 a^{ij,k}D_{ij}uD_k\varphi_1.
 \end {align*}
 
 We now estimate the summands in $C_6$.  As a preliminary step, we estimate only portions of each summand.
 From our estimates on $D\zeta$ and $D^2\zeta$, we conclude that
 \[
 a^{ij}D_{ij}\zeta +a^kD_k\zeta \ge - c
 \]
 for some $c$ determined only by $n$, $R$, $R_0$, $\beta_0$, and $\mu_5$.
 We also have $a^{ij}D_i\zeta D_j\zeta \ge0$, and
 \[
 NN_za^{ij}D_i\zeta D_ju \ge - cw_1\Lambda^{1/2} \mathscr E^{1/2}
 \]
 by  \eqref {E10.37c}, \eqref {E10.39}, \eqref{ENv}, the Cauchy-Schwarz inequality and our estimate on $D\zeta$.
 Similarly, \eqref {E10.37d}, \eqref {E10.39}, \eqref {ENv}, the Cauchy-Schwarz inequality and our estimate on $D\zeta$ imply that
 \[
 Na^{ij}D_i\zeta N_j \ge -c w_1(1+\sqrt {w_1}\varepsilon_x(\sqrt w_1)) \Lambda.
 \]
Finally, 
\[
a^{ij}D_i\zeta D_j(c^{km})D_kuD_mu \ge -c\Lambda w_1.
\]
It therefore follows from \eqref {Elambdaspecial} that
\[
C_6 \ge qw_2 \varphi_1^{q-1}[K-c]
\]

Similarly, the Cauchy-Schwarz inequality implies that
\[
a^{ij}D_i\zeta D_{jk}u\nu_1^k \ge -\frac{\varphi_1}2 \mathscr S-\frac c{\varphi_1}w_1\Lambda,
\]
and, along with \eqref {SCpspecialaijp1}, that
\[
a^{ij,k}D_{ij}uD_k\varphi_1 \ge -\frac {\varphi_1}{2w_1}\mathscr S-\frac c{w_2^{\theta/2}\varphi_1}.
\]
Hence, 
\[
S_8 \ge -(1+k_1d[1+\varepsilon_x(\sqrt{w_1})\sqrt{w_1}\,])\varphi_1^q\mathscr S - \frac c{w_2^{\theta/2}\varphi_1}\varphi_1^{q-1}w_2\Lambda(1+\varepsilon_x(v)v).
\]
From \eqref {Elambdaspecial2}, we infer that
\[
\mathscr Lw \ge qw_2\varphi_1^{q-1}(K-c)
\]
wherever $w\ge M_1$ because $w_2^{\theta/2}\varphi_1= w^{\theta/2}$.
Choosing $K$ sufficiently large gives $\mathscr Lw\ge0$ and the proof is completed just as before.
\end {proof}

We note here that Serrin also looked in \cite {Serrin} at estimates that are local with respect to $x$ and $t$. The remark following Theorem 4 of that work points out that only one additional hypothesis is needed in this case. The additional hypothesis is just there are positive constants $\mu_8$ and $\theta_1$ such that
\begin {equation}\label {ESerrin}
\mathscr E \ge \mu_8v^{\theta_1}
\end {equation}
for  $v\ge M_0$.
 (We refer also to Sections 11.4 and 11.5 of \cite {LBook} for a slightly different look at this idea, but we mention that Section 11.4 of \cite {LBook} is based very heavily on \cite {Serrin}.)
In our situation, this condition also leads to estimates which are local in $x$ and $x$, and the estimates are obtained via the same modification of our proofs as in \cite {Serrin}. We include the results for completeness only.

\begin {theorem} \label {TClocp2}
Suppose, in addition to the hypotheses of Theorem \ref {TClocp}, that there are positive constants $\mu_8$ and $\theta_1$ such that \eqref {ESerrin} holds.  Then $|Du(Y)|\le M$ for some constant $M$ determined by the same quantities as Theorem \ref {TClocp} except that, instead of depending on $\sup_{\Omega'(Y;R,R')}|Du|$, it depends on $\mu_8$ and $\theta_1$.
\end {theorem}
\begin {proof} We replace $\zeta(x)^q$ in the proof of Theorem \ref {TClocp}by
\[
\zeta(x)^q\left (1+\frac {t-s}{R'}\right)^{q_1},
\]
with $q_1=2/\theta_1$. We refer the interested reader to the proof of Theorem 4 in \cite {Serrin} or the proof of Theorem 11.3(a) in \cite {LBook} for further information.
 \end {proof}
 
In fact, there is also the possibility of studying estimates that are local with respect to $t$ as well, but we shall not discuss them in detail because they are very similar to the global estimates.

\section {Examples for parabolic problems} \label {Spe}

In this section, we jump immediately to a discussion of differential equations because, except for \eqref {Ebetat} and the connection between $\varepsilon_x$ and the differential equation, all the relevant details of the various boundary conditions have already been established in our elliptic examples.

\subsection {Capillary-type differential equations}
We assume that $a^{ij}$ has the form
\begin {equation} \label {Eparacap}
a^{ij} = h_1(v)g^{ij}
\end {equation}
for some positive function $h_1$.  Different choices for $h_1$ can be used in different theorems, but all of them require some restrictions on this function.
We first compute
\[
a^{ij,k}\eta_{ij} = h_1'(v) g^{ij}\eta_{ij} \nu_k -2h_1v^{-2}g^{ik}\delta^{jm}\eta_{ij}.
\]
A little calculation reveals that we cannot apply Theorem \ref {TCp} because \eqref {SC1aijp}  and \eqref {ELambdav2} fail.
To apply Theorem \ref {Tpspecial}, we note from our computation of $a^{ij,k}\eta_{ij}$ that 
\[
\delta_1a^{ij}\eta_{ij} \le \frac {h_1'(v)|p|}{\sqrt{h_1(v)}}(a^{ij}\delta^{km}\eta_{ik}\eta_{jm})^{1/2} +2 \frac {\sqrt{h_1(v)}}v(a^{ij}\delta^{km}\eta_{ik}\eta_{jm})^{1/2}.
\]
Hence, we obtain \eqref {SCpspecialaij1} if 
\begin {equation} \label {Eh1}
|h_1'(\sigma)|= O(\sqrt{h_1(\sigma)}\sigma^{-2})
\end {equation}
since this inequality implies that $h_1$ is bounded. In addition, it is satisfied for $h_1(\sigma)= \sigma^Q$ provided $Q=0$ or $Q\le-2$.
Of course, \eqref {SCpspecialaij2} holds with $\mu_6\equiv0$.
Finally, we have
\[
|a^{ij}_p\eta_{ij}| \le \left( \frac {h_1'(v)}{\sqrt{h_1(v)}} + \frac {2\sqrt{h_1(v)}}v\right)(a^{ij}\delta^{km}\eta_{ik}\eta_{jm})^{1/2},
\]
so \eqref {SCpspecialaijp} is also satisfied.
We assume that $a$ has the form
\begin {equation} \label {SCaH}
a(X,z,p)=vH(X,z)
\end {equation}
for some $C^1$ function $H$ satisfying $H_z\le0$.  Then conditions \eqref {SCpspecial} are satisfied with $\mu_5$ and $\mu_6$ determined only by upper bounds for $|h_1|$, $|H|$, and $|H_x|$.
Next, we assume that $b$ has the form 
\begin {equation} \label {Epbh}
b(X,z,p)=h(v)p\cdot\gamma + \psi(X,z)
\end {equation}
with $h$ satisfying \eqref {Ehpsi}.  Then the hypotheses of Theorem \ref {TNp} are satisfied with $\beta_0$, $\beta_1$, and $c_0$ determined only by the limit behavior in \eqref {Ehpsi} and bounds for $|\psi_x|$ and $|\psi_z|$.  In addition, we have $\varepsilon_x(\sigma) = c/\sigma$ and $\Lambda_0(\sigma)= c\sigma$ for a constant $c$ determined only by $\Omega$, $h(1)$, and bounds on $|\psi_x|$, and $|\psi_t|$. Hence \eqref {Elambdaspecial2} holds as well.
If we assume finally either that $\psi_z\le0$ or that $h(\sigma)$ is bounded away from zero, then \eqref {Edelta1beta} holds.

If $b$ has the form
\begin {subequations} \label {Bphtilde}
\begin {gather}
b(X,z,p) = h(\tilde v)\beta^{ij}(X,z)p_i\gamma_j +\psi(X,z) \\
\intertext {with}
\tilde v = (1+\beta^{ij}(X,z)p_ip_j)^{1/2}
\end {gather}
\end {subequations}
for some positive definite matrix $[\beta^{ij}]$ and some positive function  $h$ satisfying \eqref {Ehprime} and \eqref {Elimpsi1}, then we can apply the same analysis as in Example \ref  {SSOBC1} to conclude that the conditions of \eqref {TNp} are satisfied and that \eqref {Elambdaspecial2} holds.  If, in addition, $\psi_z\le0$ or $h$ is bounded away from zero and if $\beta^{ij}$ is independent of $z$, then \eqref {Edelta1beta} is valid.

Summarizing, we can apply Theorem \ref {Tpspecial} for $a^{ij}$ having the form \eqref {Eparacap}, $a$ having the form \eqref {SCaH}, and $b$ having the form \eqref {Epbh} with $h$ satisfying \eqref {Ehpsi} provided
\begin {itemize}
\item
$|H_x|$, $|H_z|$, $|\psi_x|$, $|\psi_z|$, and $|\psi_t|$ are uniformly bounded,
\item
$H_z\le0$,
\item $h_1$ satisfies \eqref {Eh1},
\item
$\psi_z\le0$ or $h$ is bounded away from zero.
\end {itemize}
If, instead $b$ has the form \eqref {Bphtilde} with $h$ satisfying \eqref {Ehprime} and \eqref {Elimpsi1} and if $\beta^{ij}$ is independent of $z$, then we can apply Theorem \ref {Tpspecial} under exactly the same hypotheses.
In particular, this theorem applies when $h_1\equiv1$ and when $h_1(\sigma)=\sigma^Q$ with $Q\le-2$. 
Moreover, if $Q<-2$, then we can obtain a local gradient bound via Theorem \ref {Tpspecialloc}.

If $b$ has the form
\begin {equation} \label {Bpq}
b(X,z,p)= v^{q(X)} p\cdot \gamma+\psi(X,z)
\end {equation}
for some $C^1$ function $q$ with $q\ge-1$, then we again have the conditions of Theorem \ref {TNp} satisfied but with
\[
\varepsilon_x(\sigma) = K\frac {\ln\sigma}\sigma
\]
for some positive constant $K$; however, the condition $\Lambda_0(\sigma)=\mu_5\sigma$ requires $q$ to be independent of $t$.
Hence, if $a^{ij}$ has the form \eqref {Eparacap}, $a$ has the form \eqref {SCaH}, and if $b$ has the form \eqref {Bpq}, then we can apply Theorem \ref {Tpspecial} provided
\begin {itemize}
\item
$|H_x|$, $|H_z|$, $|\psi_x|$, $|\psi_z|$, and $|\psi_t|$ are uniformly bounded,
\item
$H_z\le0$,
\item $h_1$ satisfies \eqref {Eh1},
\item
$|h_1(\sigma)|=O(1/\ln\sigma)$,
\item
$\psi_z\le0$ or $q\ge0$,
\item 
$q$ is independent of $t$.
\end {itemize}
In particular, this example provides a gradient estimate if $h_1(v)=v^Q$ with $Q\le-2$ but not if $h_1\equiv1$.

If $b$ has the form
\begin {equation} \label {Epbetastar}
b(X,z,p) = h(v)\beta(X,p)\cdot p + \psi(X,z)
\end {equation}
with $\beta$ a unit vector such that $\beta\cdot \gamma \ge \beta_*$ for some  constant $\beta_*\in (0,1)$ and $h$ a positive function satisfying \eqref{Ebetastarh}, then our analysis yields a gradient bound assuming that
\begin {itemize}
\item
$|H_x|$, $|H_z|$, $|\psi_x|$, $|\psi_z|$, and $|\psi_t|$ are uniformly bounded,
\item
$H_z\le0$,
\item $h_1$ satisfies \eqref {Eh1}.
\end {itemize}

We turn to the conditions on $a^{ij}$ in Theorem \ref {Tpintgrad}.
First, our expression for $a^{ij,k}$ implies that
\[
(|p|^2a^{ij,k}+ 4a^{ij}\delta^{jm})\eta_{ij}\xi_k = h_1'(v)|p|^2g^{ij}\eta_{ij} \nu\cdot \xi +h_1(v)H_1,
\]
with
\begin {align*}
H_1&= |p|^2\left( -\frac {2g^{ik}\delta^{jm}\eta_{ij}p_m\xi_k}{v^2}+4g^{ik}\delta^{jm}\eta_{ij}p_m\xi_k\right) \\
&= \frac {4v^2-2|p|^2}{v^2}g^{ik}\delta^{jm}\eta_{ij}p_m\xi_k \\
&=\left(2- \frac 2{v^2}\right)g^{ik}\delta^{jm}\eta_{ij}p_m\xi_k.
\end {align*}
Since $\delta_3a^{ij}\equiv0$, we infer \eqref {SCTpintgradaij} with $\mu_7=2$ and suitable $\mu_6$ provided 
\[
|h_1'(\sigma)|=O(\sqrt{h_1(\sigma)}\sigma^{-1-1/q})
\]
for some $q>2$.  It's easy to see that $h_1(\sigma)=\sigma^Q$ satisfies these hypotheses for $Q=0$ if we take $q>2$ arbitrarily. For $Q<0$, these hypotheses are satisfied if we take $q>\max\{2, -2/Q\}$.
Moreover, \eqref {SCLambdabound} is satisfied for $h_1(\sigma) =\sigma^Q$ for $Q\le0$.
Hence, Theorem \ref {Tpintgrad} applies under the conditions given about which satisfy the hypotheses of Theorem \ref {Tpspecial}.

In particular, we can obtain a gradient bound for the problem
\begin {gather*}
-u_t+g^{ij}D_{ij}u +vH(X,u)=0 \text { in } \Omega, \\
\gamma \cdot Du+\psi(X,u)=0 \text { on } S\Omega, \\
u=u_0 \text { on } B\Omega
\end {gather*}
provided $|H_x|$, $|H_z|$, $|\psi_x|$, $|\psi_z|$, $|\psi_t|$, and $|Du_0|$ are uniformly bounded and $H_z\le0$.  This problem is similar to the one studied by Mizuno and Takasao in \cite {Mizuno}.  As previously pointed out, we improve their results by allowing nonzero boundary data but we need to assume that $H$ is differentiable.

\subsection {False mean curvature equations}
We now consider $a^{ij}$ having the form \eqref {EFMCE} with $\tau\ge0$ and $a^{ij}_*$ satisfying slightly more general conditions than in Example \ref {SFMCE}.  We assume that that there are a positive function $\Lambda_*$ and a positive constant $\mu_*$ such that \eqref {Estar} holds for all $(x,z,p)\in \Omega\times\mathbb R\times\mathbb R^n$ and all $\xi\in\mathbb R^n$, and we assume that \eqref {Estarderiva} holds.  In place of \eqref {Estarderivb}, we assume that
\begin {equation} \label {Estarderivb1}
|p| \left| \frac {\partial a^{ij}_*}{\partial z}\right| + \left| \frac {\partial a^{ij}_*}{\partial x}\right| =o (|p|^2(\tau\Lambda_*)^{1/2}).
\end {equation}
We also assume that $\tau$ satisfies \eqref {Etau} and
\begin {equation} \label {Elambdastartau1}
\Lambda_*= O(\tau|p|^2),
\end {equation}
while $a$ satisfies \eqref {Ea2}.

Then \eqref {SC1aijp}, \eqref {Elambdap}, and \eqref {pSC1aij2} were verified in Example \ref {SFMCE} with $\tilde\mu_*(s) = K_1/s$ and a suitable $\tilde\mu_\varepsilon$ satisfying \eqref {Etildelimiteps}.  To check \eqref {pSCaij}, we first use \eqref {Estarderiva} to conclude that there is a constant $c$ such that
\[
|\delta_1a^{ij}_*| \le c(\tau\Lambda^*)^{1/2}.
\]
Then, for any matrix $[\eta_{ij}]$ and this constant $c$, we have
\begin {align*}
\delta_1a^{ij}_*\eta_{ij} &\le nc\tau^{1/2}|p|(\Lambda_*\sum_{i,j=1}^n|\eta_{ij}|^2)^{1/2} \\
&\le \frac {nc}{\mu_*|p|} (\tau|p|^4)^{1/2}(a^{ij}_*\delta^{jk}\eta_{ik}\eta_{jm})^{1/2} \\
&\le \frac {nc}{\mu_*|p|} (\tau|p|^4)^{1/2}(a^{ij}\delta^{jk}\eta_{ik}\eta_{jm})^{1/2}.
\end {align*} 
Moreover, direct computation shows that
\[
\delta_1(\tau p_ip_j) =(\delta_1+2)\tau p_ip_j,
\]
so similar reasoning, using \eqref {Etaua}, shows that
\[
\delta_1(\tau p_ip_j)\eta_{ij} \le \frac c{|p|}\mathscr E^{1/2}(a^{ij}\delta^{jk}\eta_{ik}\eta_{jm})^{1/2},
\]
and therefore \eqref {pSCaij} is also satisfied.

Since $\lim_{\sigma\to\infty}\varepsilon_x(\sigma)=0$ in all our examples, \eqref {Etildelimit1} is also satisfied.

To verify \eqref {Edelta1b} and \eqref {Elimbz}, we assume more about $b$.
 If $b$ has the form \eqref {Ebh} with $h$ satisfying \eqref {Ehpsi} for some nonnegative constant $h_0$, we consider two cases. First, if
 \[
 \lim_{\sigma\to\infty} \sigma h(\sigma)<\infty,
 \]
 we assume \eqref {Ehprime1} and $\psi_z=0$ (although this condition can be relaxed slightly to assuming that $\psi_z$ is nonpositive and sufficiently small).  Then, as already shown, \eqref {Edelta1b} and \eqref {Elimbz} hold.
 On the other hand, if
 \[
 \lim_{\sigma\to\infty}\sigma h(\sigma)=0,
 \]
 then we assume \eqref {Ehprime2}.  Again, we have already shown that \eqref {Edelta1b} holds in this case and that \eqref {Elimbz} holds.
 
 If $b$ has the form \eqref {Bphtilde} with $h$ satisfying \eqref {Ehpsi} and \eqref {Elimpsi1}, then we obtain a gradient bound under the same additional restrictions on $h$ and $\psi_z$ as for the previous form for $b$.
 
If $b$ has the form \eqref {Bpq}, then we assume that $\inf q>-1$.  From this assumption, we immediately infer \eqref {Etildeq} and hence \eqref {Edelta1b} as well as \eqref {Elimbz}.

If $b$ has the form \eqref {Epbetastar} with $h$ satisfying \eqref {Ebetastarh}, we also infer a global gradient estimate.
 
 Furthermore, if we assume also that there is a positive constant $\theta_1$ such that 
 \[
 \tau \ge v^{\theta_1-4},
 \]
 then we obtain a gradient estimate which is local in $x$ and $t$.
 
Let us now assume that $a^{ij}$ has the form \eqref {EFMCE} with $a^{ij}_*$ and $\tau$ satisfying a different set of conditions.  First, we suppose that there are positive functions $\lambda_*$ and $\Lambda_*$ such that
\[
\lambda_*|\xi|^2\le a^{ij}_*\xi_i\xi_j \le \Lambda_*|\xi|^2, 
 \]
 for all $\xi\in\mathbb R^n$. We also assume that
 \[
 \Lambda_*= O(v^{-2}),
 \]
 and that
\[
 |p|^2|\left| \frac {\partial a^{ij}_*}{\partial p}\right| +
 |p|\left| \frac {\partial a^{ij}_*}{\partial z}\right| +
 \left| \frac {\partial a^{ij}_*}{\partial x}\right|= O(\sqrt{\lambda_*}),
 \]
 while $\tau$ satisfies
 \[
 \tau = O(v^{-4}), \
 |p|^2|\tau_p|+|p||\tau_z| +|\tau_x|=O(\sqrt\tau).
 \]
Finally, we assume that $a$ satisfies \eqref {SCpspeciala1} and \eqref {SCpspeciala2}.  This will be the case if there is a $C^1$ function $a_0$ such that $a(x,z,p)=va_0(x,z)$.
By imitating our previous arguments, we see that all the hypotheses of Theorem \ref {Tpspecial} as long as $b$ satisfies the hypotheses of Theorem \ref {TNp} with $\varepsilon_x$ constant.  It follows that we have a gradient estimate for this class of parabolic differential equations and for all of the boundary conditions in our examples. In particular, we can get a gradient estimate for the problem
\begin {align*}
-u_t+ v^{-4}(\delta^{ij}+D_iuD_ju)D_{ij}u &=0 \text { in } \Omega, \\
\frac {Du\cdot\gamma}v+\psi(x,z) &=0 \text { on }S\Omega, \\
u &=u_0 \text { on } B\Omega
\end {align*}
as long as $\psi_z$, $\psi_x$, and $Du_0$ are bounded.

 \subsection {Uniformly parabolic equations}
 
 Our final example is the parabolic analog of Example \ref {SSUE}.  We assume that there is a positive constant $\mu_*$ such that \eqref {SC1UE} is satisfied.  In place of \eqref {SCUE}, we assume that
\begin {subequations} \label {SCUP}
 \begin {gather}
|p| |a^{ij}_p| = O(\Lambda), \\
|p||a^{ij}_z|+|a^{ij}_x| = o(\Lambda), \label {SCUPaijo} \\
|a| +|p||a_p| = O(\Lambda|p|^2), \\
|a_x| = o(\Lambda|p|^3), \label {SCUPaxo} \\
a_z\le o(\Lambda|p|^2). \label {SCUPazo}
 \end {gather}
 \end {subequations}
Then conditions \eqref {Etildelimiteps}, \eqref {Elambdap}, \eqref {Etildelimit1}, \eqref  {SC1aijp}, \eqref {pSC1aij2}, and \eqref {pSCaij}, and  are all satisfied with $\tilde\mu_\varepsilon$ determined by $\varepsilon$ and the limit behavior in \eqref {SCUPaijo}, \eqref {SCUPaxo}, and \eqref {SCUPazo}, and $\tilde\mu_*(\sigma)=K_1/\sigma$ for a sufficiently large constant $K_1$. 
We therefore obtain a gradient bound for all of our boundary conditions under the restrictions mentioned in the previous example.

There is one special case when we can replace \eqref {SCUP} by \eqref {SCUE}.  If we assume \eqref {SC1UE} holds with $\Lambda$ bounded from above and below by positive constants and if $|a|=O(|p|^2)$ (which is the same as $|a|=O(\Lambda|p|^2)$ in this situation), then a modulus of continuity is known. Theorem 3 of \cite {Nazarov} gives a slightly weaker result, but standard arguments can be used to obtain this result from that theorem.  (See also Lemma 13.14 of \cite {LBook}, which states the result with a very minimal proof.)  It follows that, if $a^{ij}$ and $a$ also satisfy \eqref {SC1UE}, then we obtain a gradient estimate for all of our examples of boundary conditions.  This result was first proved by Ural$'$tseva in \cite {UR3} as Theorem 1 (but assuming that $b$ is twice differentiable with second derivatives satisfying certain structure conditions) by very different means.  It was also proved as Theorem 13.13 in \cite {LBook} but, just as in the elliptic case, the proof there takes advantage of the exact form of the interior gradient estimate.

If we assume further that $\Lambda=O(|p|^{-2})$, then the discussion of the previous example (with $\tau\equiv0$) gives a gradient estimate under the hypotheses \eqref {SC1UE},
\begin {gather*}
|p|^2|a^{ij}_p| +|p| |a^{ij}_z|+ |a^{ij}_x| =O(\sqrt\Lambda), \\
|p|^2|a_p|+|p||a| +|a_x|= O(|p|^3), \\
a_z\le O(|p|^2),
\end {gather*}
and any of our examples of boundary conditions.

\begin {bibsection}
\begin {biblist}

\bib{Dong}{article}{
   author={Dong, G. C.},
   title={Initial and nonlinear oblique boundary value problems for fully
   nonlinear parabolic equations},
   journal={J. Partial Differential Equations Ser. A},
   volume={1},
   date={1988},
   number={2},
   pages={12--42},
}

\bib{EH}{article}{
   author={Ecker, K.},
   author = {Huisken, G.},
   title={Interior estimates for hypersurfaces moving by mean curvature},
   journal={Invent. Math.},
   volume={105},
   date={1991},
   pages={547--569},
}

\bib {GT}{book}{
author = {Gilbarg, D.},
author = {Trudinger, N. S.},
title = {Elliptic Partial Differential Equations of Second Order},
publisher = {Springer-Verlag},
place = {Berlin},
date = {2001},
note = {reprint of 1998 edition},
}

\bib {H} {article}{
author = {Huisken, G.},
title = {Non-parametric mean curvature evolution with boundary conditions},
journal = {J. Differential Equations},
volume = {77},
date = {1989},
pages = {369\ndash385},
}

\bib {JDE49} {article}{
author = {Lieberman, G. M.},
title = {The conormal derivative problem for elliptic equations of variational type},
journal = {J. Differential Equations},
volume = {49},
date = {1982},
pages = {218\ndash257},
}

\bib {PJM} {article}{
author = {Lieberman, G. M.},
title = {Regularized distance and its applications},
journal = {Pac. J. Math.},
volume = {117},
date = {1985},
pages = {329\ndash352},
}

\bib {NUE} {article}{
author = {Lieberman, G. M.},
title = {Gradient estimates for solutions of nonuniformly elliptic oblique derivative problems},
journal = {Nonl. Anal.},
volume = {11},
date = {1987},
pages = {49\ndash61},
}

\bib {maxcap} {article}{
author = {Lieberman, G. M.},
title = {Gradient estimates for capillary-type problems via the maximum principle},
journal = {Comm. Partial Differential Equations},
volume = {13},
date = {1988},
pages = {33\ndash59},
}

\bib {IUMJ88} {article}{
author = {Lieberman, G. M.},
title = {The conormal derivative problem for non-uniformly parabolic equations},
journal = {Indiana U. Math. J.},
volume = {37},
date = {1988},
pages = {23\ndash72},
}

\bib {CPDE16} {article}{
author = {Lieberman, G. M.},
title = {The natural generalization of the natural conditions of Ladyzhenskaya and Ural$'$tseva for elliptic equations},
journal = {Comm. Partial Differential Equations},
volume = {16},
date = {1991},
pages = {311\ndash361},
}

\bib {LBook}{book}{
author = {Lieberman, G. M.},
title = {Second Order Parabolic Differential Equations},
publisher={World Scientific Publishing Co., Inc., River Edge, NJ},
date={1996},
}

\bib {NA119} {article}{
author = {Lieberman, G. M.},
title = {Gradient estimates for singular fully nonlinear elliptic equations},
journal = {Nonlin. Anal.},
volume = {119},
date = {2015},
pages = {382\ndash397},
}

\bib {ODPbook}{book}{
author = {Lieberman, G. M.},
title = {Oblique Derivative Problems for Elliptic Equations},
publisher = {World Scientific},
address = {Hackensack, N. J.},
year = {2013},
}

\bib {maxcap2} {article}{
author = {Lieberman, G. M.},
title = {Gradient estimates for capillary-type problems via the maximum principle, A second look},
note = {to appear}
}

\bib {LT} {article}{
author = {Lieberman, G. M.},
author = {Trudinger, N. S.},
title = {Nonlinear oblique boundary value problems for nonlinear elliptic equations},
journal = {Trans. Amer. Math. Soc.},
volume = {295},
date = {1986},
pages = {509\ndash546},
}

\bib {XuMa} {article}{
author= {Ma, X. N.},
author = {Xu, J. J.},
title = {Gradient estimates of mean curvature equations with Neumann boundary conditions},
journal = {Adv. Math.},
volume = {290},
date = {2016},
pages = {1010\ndash1036},
}

\bib {Mizuno} {article}{
author= {Mizuno, M.},
author = {Takasao, K.},
title = {Gradient estimates for mean curvature flow with Neumann boundary conditions},
journal = {Nonlinear Diff. Eqn. Appl.},
volume = {24},
number ={4},
date = {2017},
pages = {Art. 32,24},
}

\bib {Nazarov}{article}{
   author={Nazarov, A. I.},
   title={H\"older estimates for bounded solutions of problems with an oblique
   derivative for parabolic equations of nondivergence structure},
   language={Russian},
   note={Translated in J. Soviet Math.\ {\bf 64} (1993), no.\ 6,
   1247--1252},
   conference={
      title={Nonlinear equations and variational inequalities. Linear
      operators and spectral theory (Russian)},
   },
   book={
      series={Probl. Mat. Anal.},
      volume={11},
      publisher={Leningrad. Univ., Leningrad},
   },
   date={1990},
   pages={37--46, 250},
}

\bib {Serrin}{incollection}{
author = {Serrin, J.},
title = {Gradient estimates for solutions of nonlinear elliptic and parabolic equations},
booktitle = {Contributions to Nonlinear Functional Analysis},
publisher = {Academic Press},
place = {New York},
date = {1971},
pages = {565--601},
}

\bib {UR0} {article} {
author = {Ural'ceva, N. N.},
year = {1971},
language = {Russian},
translation = {
journal = {Proc. Steklov Math. Inst.},
volume = {116},
year = {1971},
pages = {227--237},
},
title = {Nonlinear boundary value problems for equations of minimal-surface type},
journal = {Trudy Mat. Inst. Steklov.},
volume = {116},
pages = {217--226},
}

\bib{UR3}{incollection}{ 
   author={Uraltseva, N. N.},
   title={Gradient estimates for solutions of nonlinear parabolic oblique
   boundary problem},
   conference={
      title={Geometry and nonlinear partial differential equations},
      address={Fayetteville, AR},
      date={1990},
   },
   book={
      series={Contemp. Math.},
      volume={127},
      publisher={Amer. Math. Soc., Providence, RI},
   },
   date={1992},
   pages={119--130},
}
	
\end {biblist}
\end {bibsection}
\end {document}